\documentclass{siamart0516}

\usepackage[utf8]{inputenc}
\usepackage[T1]{fontenc}
\usepackage{iflang}
\usepackage{csquotes}

\usepackage[scaled = 0.92]{helvet}
\usepackage{lmodern}
\usepackage{mathrsfs} 

\usepackage{tikz,pgfplots}
\usepackage{tikzscale} 
\usetikzlibrary{decorations}
\usetikzlibrary{decorations.markings}
\usetikzlibrary{decorations.pathreplacing}
\usetikzlibrary{shapes,arrows,patterns}
\usetikzlibrary{arrows.meta}
\usetikzlibrary{patterns}
\usetikzlibrary{positioning,calc}
\usetikzlibrary{circuits.ee.IEC}
\usetikzlibrary{intersections}
\pgfplotsset{compat=newest}
\usetikzlibrary{fit,matrix}
\usepackage{pgfkeys}

\newlength\figureheight
\newlength\figurewidth
\usetikzlibrary{external} 
\usepgfplotslibrary{external} 

\usepackage[colorinlistoftodos,prependcaption]{todonotes}

\usepackage{empheq}
\usepackage{enumerate}
\usepackage{fancyhdr}
\usepackage[headheight = 14 pt]{geometry}
\usepackage{setspace}

\usepackage{graphicx}
\usepackage{import} 
\usepackage{pdfpages}

\usepackage[font=small]{subfig}
\usepackage{caption}
\usepackage{epstopdf}

\usepackage{amssymb}
\usepackage{siunitx}
\usepackage{xfrac}
\usepackage{cancel}
\usepackage{nicefrac}

\usepackage{longtable}
\usepackage{multirow}

\usepackage[square, numbers]{natbib}
\usepackage{acro}
\usepackage{thmtools}

\usepackage{hyperref}

\newcommand{\smo}{\smashoperator}
\usepackage{upgreek} 
\usepackage{stmaryrd} 

\DeclareMathOperator{\flambda}{\uplambda}
\renewcommand{\Re}{\mathrm{Re}}

\renewcommand{\vec}[1]{\boldsymbol{#1}} 
\newcommand{\mat}[1]{\boldsymbol{#1}} 
\newcommand{\matn}[1]{\boldsymbol{\mathrm{#1}}} 
\newcommand{\pd}{\partial}                      
\newcommand{\abs}[1]{| #1 |}                    

\newcommand{\comment}[1]{ }  

\makeatletter
\renewcommand{\cleardoublepage}{\clearpage\if@twoside \ifodd\c@page\else
	\hbox{}
	\vspace*{\fill}
	\thispagestyle{empty}
	\newpage
	\if@twocolumn\hbox{}\newpage\fi\fi\fi}
\makeatother   

\newcommand{\D}{\displaystyle}

\providecommand*{\eu}{\ensuremath{\mathrm{e}}}

\makeatletter
\newsavebox\myboxA
\newsavebox\myboxB
\newlength\mylenA

\newcommand*\xoverline[2][0.75]{%
	\sbox{\myboxA}{$\m@th#2$}%
	\setbox\myboxB\null
	\ht\myboxB=\ht\myboxA%
	\dp\myboxB=\dp\myboxA%
	\wd\myboxB=#1\wd\myboxA
	\sbox\myboxB{$\m@th\overline{\copy\myboxB}$}
	\setlength\mylenA{\the\wd\myboxA}
	\addtolength\mylenA{-\the\wd\myboxB}%
	\ifdim\wd\myboxB<\wd\myboxA%
	\rlap{\hskip 0.5\mylenA\usebox\myboxB}{\usebox\myboxA}%
	\else
	\hskip -0.5\mylenA\rlap{\usebox\myboxA}{\hskip 0.5\mylenA\usebox\myboxB}%
	\fi}
\makeatother

\newdimen\savedrulethickness

\makeatletter
\providecommand*{\dd}%
{\@ifnextchar^{\DIfF}{\DIfF^{}}}
\def\DIfF^#1{%
	\mathop{\mathrm{\mathstrut d}}%
	\nolimits^{#1}\gobblespace}
\def\gobblespace{%
	\futurelet\diffarg\opspace}
\def\opspace{%
	\let\DiffSpace\!%
	\ifx\diffarg(%
	\let\DiffSpace\relax
	\else
	\ifx\diffarg[%
	\let\DiffSpace\relax
	\else
	\ifx\diffarg\{%
	\let\DiffSpace\relax
	\fi\fi\fi\DiffSpace}

\usepackage{pdflscape} 
\makeatletter
\global\let\orig@begin@landscape=\landscape%
\global\let\orig@end@landscape=\endlandscape%
\gdef\@true{1}
\gdef\@false{0}
\gdef\landscape{%
	\global\let\within@landscape=\@true%
	\orig@begin@landscape%
}%
\gdef\endlandscape{%
	\orig@end@landscape%
	\global\let\within@landscape=\@false%
}%
\@ifpackageloaded{pdflscape}{%
	\gdef\pdf@landscape@rotate{\PLS@Rotate}%
}{
\gdef\pdf@landscape@rotate#1{}%
}
\let\latex@outputpage\@outputpage
\def\@outputpage{
	\ifx\within@landscape\@true%
	\if@twoside%
	\ifodd\c@page%
	\gdef\LS@rot{\setbox\@outputbox\vbox{%
			\pdf@landscape@rotate{-90}%
			\hbox{\rotatebox{90}{\hbox{\rotatebox{180}{\box\@outputbox}}}}}%
	}%
	\else%
	\gdef\LS@rot{\setbox\@outputbox\vbox{%
			\pdf@landscape@rotate{+90}%
			\hbox{\rotatebox{90}{\hbox{\rotatebox{0}{\box\@outputbox}}}}}%
	}%
	\fi%
	\else%
	\gdef\LS@rot{\setbox\@outputbox\vbox{%
			\pdf@landscape@rotate{+90}%
			\hbox{\rotatebox{90}{\hbox{\rotatebox{0}{\box\@outputbox}}}}}%
	}%
	\fi%
	\fi%
	\latex@outputpage%
}
\makeatother

\definecolor{tum blue}{HTML}{0065BD}
\definecolor{tum blue 1}{HTML}{98C6EA}
\definecolor{tum blue 2}{HTML}{64A0C8}
\definecolor{tum blue 3}{HTML}{0073CF}
\definecolor{tum blue 4}{HTML}{005293}
\definecolor{tum blue 5}{HTML}{003359}

\definecolor{tum green}{HTML}{A2AD00}
\definecolor{tum orange}{HTML}{E37222}
\definecolor{tum ivory}{HTML}{DAD7CB}

\definecolor{tum dia violet}{HTML}{69085A}
\definecolor{tum dia dark blue}{HTML}{0F1B5F}
\definecolor{tum dia turquoise}{HTML}{00778A}
\definecolor{tum dia dark green}{HTML}{007C30}
\definecolor{tum dia light green}{HTML}{679A1D}
\definecolor{tum dia light yellow}{HTML}{FFDC00}
\definecolor{tum dia dark yellow}{HTML}{F9BA00}
\definecolor{tum dia dark orange}{HTML}{D64C13}
\definecolor{tum dia red}{HTML}{C4071B}
\definecolor{tum dia dark red}{HTML}{9C0D16}

\newsiamremark{remark}{Remark}

\usepackage{lipsum}
\usepackage{amsfonts}
\usepackage{graphicx}
\usepackage{epstopdf}
\usepackage{algorithmic}
\ifpdf
  \DeclareGraphicsExtensions{.eps,.pdf,.png,.jpg}
\else
  \DeclareGraphicsExtensions{.eps}
\fi

\newcommand{\TheTitle}{Impulse response of bilinear systems based on \\ Volterra series representation} 
\newcommand{\TheShortTitle}{Impulse response of bilinear systems}
\newcommand{\TheAuthors}{M. Cruz Varona and R. Gebhart}

\headers{\TheShortTitle}{\TheAuthors}

\title{{\TheTitle}\thanks{Submitted to the editors March 29, 2018.}}

\author{
  Maria Cruz Varona\footnotemark[2]
  \and
  Raphael Gebhart\thanks{Chair of Automatic Control (Prof. Dr.-Ing. B. Lohmann), Department of Mechanical Engineering, Technical University of Munich, Boltzmannstr. 15, 85748 Garching, Germany (\email{maria.cruz@tum.de},
    \email{raphael-gebhart@t-online.de}).}
}

\usepackage{amsopn}

\ifpdf
\hypersetup{
  pdftitle={\TheTitle},
  pdfauthor={\TheAuthors}
}
\fi

\begin{document}

\maketitle

\begin{abstract}
%
This paper focuses on the systems theory of bilinear dynamical systems using the Volterra series representation. The main contributions are threefold. First, we gain an input-output representation in the frequency domain, where the Laplace transform of the kernels can indeed be interpreted as transfer functions. Then, we derive the response of bilinear systems to a nascent delta function in time domain, i.e. the \emph{impulse response of bilinear systems}. Finally, we study the relationships between this novel impulse response and the well-known Volterra kernels and adjust those to be compatible to the impulse response.
\end{abstract}

\begin{keywords}
  bilinear systems, Volterra series, integral kernels, multidimensional Laplace transform, impulse response, systems theory
\end{keywords}

\begin{AMS}
	93C10, 93C15, 93A15
\end{AMS}

\section{Introduction}
Bilinear dynamical systems are weakly nonlinear systems that are linear separately in state and in input, but not jointly linear in both. They represent an interface between fully nonlinear and linear dynamical systems. In the 1970s, bilinear systems began to receive attention, since many physical phenomena can be described by them, e.g. biological, physiological and economical processes, as well as applications such as catalysts in chemical reactions, nuclear fission and stochastic problems \cite{mohler1970natural,mohler1973bilinear,mohler1980overview,mohler1991nonlinear}. The description of the input-output behavior, the analysis of system-theoretic concepts and the minimal realization theory of bilinear state-space systems were first investigated in \cite{bruni1971on,isidori1973direct,isidori1973realization,bruni1974bilinear} and further studied in \cite{rugh1981nonlinear,elliott2009bilinear,isidori2013nonlinear}. Recently, a renewed interest in bilinear systems has aroused in the context of model order reduction. The close relation of this system class to linear state-space systems as well as the previous development of similar system-theoretic concepts (e.g. Gramians, transfer functions, $\mathcal{H}_2$-norm, etc.) based on the Volterra series representation has enabled the generalization of well-known linear reduction methods to the bilinear setting \cite{zhang2002h2,bai2006projection,breiten2010krylov,breiten2012interpolation,breiten2013interpolatory,flagg2012interpolation,flagg2015multipoint,ahmad2017implicit}.

Consider a multiple-input multiple-output (MIMO) bilinear system of the form:
\begin{equation}\label{eq:bilinear-FOM}
	\begin{aligned}
		\dot{\vec x}(t) &= \mat A \vec x(t) + \sum_{j = 1}^{m} \mat N_j u_j(t) \vec x(t) + \mat B \vec u(t)\,, \quad \vec x(0) = \vec x_0\,, \\
		\vec y(t) &= \mat C \vec x(t)\,,
	\end{aligned}
\end{equation}
where the matrices~$\mat A, \mat N_j \in \mathbb{R}^{n \times n}$ for $j=1,\ldots,m$, $\mat B \in \mathbb{R}^{n \times m}$ and $\mat C \in \mathbb{R}^{p \times n}$. The vectors~$\vec x(t) \in \mathbb{R}^n$, $\vec u(t) \in \mathbb{R}^{m}$ and $\vec y(t) \in \mathbb{R}^{p}$ denote the state, inputs and outputs of the system, respectively. For the sake of brevity, throughout the paper we will derive the results for the single-input single-output (SISO) case\footnote{For SISO ($m\!=\!1$, $p\!=\!1$), replace:~$\sum_{j = 1}^{m} \mat N_j u_j(t) \rightarrow \mat N u(t)$,~$\mat B\vec u(t) \!\!=\!\! \sum_{j = 1}^{m}\vec b_j u_j(t) \rightarrow \vec bu(t)$,~$\mat C \rightarrow \vec c^{\mathsf T}$,~$\vec y(t) \rightarrow y(t)$.}, but we will also state the formulas for the MIMO case. In case that the bilinear system is given in implicit form with a regular matrix~$\mat E \in \mathbb{R}^{n \times n}$, i.e.~$\det \mat E \not= 0$, then the theoretical statements still hold true by just replacing~$\mat A,\mat N_j, \mat B \rightarrow \mat E^{-1}\mat A, \mat E^{-1}\mat N_j, \mat E^{-1}\mat B$. Note, however, that from a computational point of view the inversion of~$\mat E$ becomes prohibited in the large-scale setting (e.g. $n=10^4$ and bigger) and should be circumvented in algorithms.

In this paper, we extensively investigate the systems theory for bilinear dynamical systems making use of the Volterra series representation. Our main contribution is the derivation of the response of a bilinear system to an impulse input in time domain, i.e. the \emph{impulse response of bilinear systems}, together with the adjustment of the Volterra kernels, especially along lines of equal time arguments, such that the kernels and the impulse response are compatible to each other. The motivation for this endeavor lies in the fact that, in the linear setting, the impulse response constitutes a very important characteristic of the system, with which the output response for arbitrary inputs can be computed. In the bilinear setting, however, the input-output behavior is characterized by means of an infinite series of multidimensional kernels that are normally not defined at certain discontinuities \cite[pp. 14-16]{rugh1981nonlinear}. Since values at discontinuities are not taken into account through the multidimensional (inverse) Laplace transform and these precise values will turn out to be crucial for the response computation, we will derive the impulse response of bilinear systems directly in time domain and adjust the Volterra kernels afterwards to obtain consistent answers. 

The paper is organized as follows. In the next section, we will survey the main fundamentals concerning bilinear systems theory. This will include the characterization of the input-output behavior in the time domain as well as the multidimensional Laplace transform of the Volterra kernels. In \cref{sec:inp-out-freq-dom} we will focus on the multidimensional Laplace transform of the output to gain an input-output representation in the frequency domain, where the Laplace transform of the kernels can be interpreted as transfer functions. The derivation of the impulse response of bilinear systems to a nascent delta function is given in \cref{sec:imp-resp}, and in \cref{sec:5} the multidimensional Volterra kernels are adjusted to be compatible and consistent with the previously derived impulse response. Finally, we conclude the paper with some final remarks and conclusions in \cref{sec:conclusions}.

\section{Background on bilinear systems theory} 
In this section, we review some system theoretic concepts of bilinear systems. The groundwork for the existing bilinear system theory is laid by the \emph{Volterra series expansion} for general nonlinear dynamical systems, where the solution $\vec x(t)$ is constructed as an infinite sum of multivariable convolution integrals by applying the Picard fixed-point iteration. Consequently, the nonlinear system can be alternatively interpreted as an infinite sequence of interconnected (cascaded) subsystems. Pursuing these considerations, the bilinear system \eqref{eq:bilinear-FOM} can be represented by the following infinite series of coupled linear subsystems\footnote{A state-equation for each subsystem can be derived by applying an input of the form~$\alpha u(t)$ and assuming the response~$\vec x(t) = \alpha \vec x_1(t) + \alpha^2 \vec x_2(t) + \ldots$.}:
\begin{equation}\label{eq:bilinear-subsystems}
\begin{aligned}
	\dot{\vec x}_1(t) &= \mat A \vec x_1(t) + \vec b u(t)\, , \quad & \vec x_1(0) &= \vec x_0\, , \\
	\dot{\vec x}_k(t) &= \mat A \vec x_k(t) + \mat N u(t) \vec x_{k-1}(t)\, , \quad & \vec x_k(0) &= \vec 0\, , \quad k \geq 2\, .
\end{aligned}
\end{equation}
In \cite{bruni1971on} it is proven that, for bounded inputs $u(t)$, the sequence of solutions $\vec x_k(t)$ of the linear subsystems converges to the solution $\vec x(t)$ of the bilinear system for $k \to \infty$: $\vec x(t) = \sum_{k=1}^{\infty} \vec x_k(t)$.\\
In the following, we make use of the subsystems to obtain the input-output representation for bilinear systems.

\subsection{Input-output representation in the time domain} \label{subsec:i-o_time-dom}
First, we aim at obtaining a descriptive equation for the output response of the bilinear system to given inputs in the time domain. Depending on the used description form, the input-output representation will be given in terms of the so-called \emph{triangular} and \emph{regular kernels}. Both description forms can be equivalently transformed into each other by change of variables.

\subsubsection{Triangular kernel representation}
Making use of \eqref{eq:bilinear-subsystems}, the solution~$\vec x_k(t)$ of the~$k$-th subsystem is in general given by:
\begin{equation}
\begin{aligned}
\vec x_1(t) &= \int_{\tau = 0}^t \eu^{\mat A(t - \tau)} \vec b u(\tau)\, \dd \tau + \eu^{\mat At}\vec x_0\, , \\
\vec x_k(t) &= \int_{\tau = 0}^t \eu^{\mat A(t - \tau)} \mat N u(\tau) \vec x_{k-1}(\tau) \, \dd \tau\,, \quad k \geq 2\, .
\end{aligned}
\end{equation}
Thus, the solution~$\vec x_2(t)$ of the second subsystem is:
\begin{equation}\label{eq:sy_lo_su1}
\begin{aligned}
\vec x_2(t) &= \int_{\tau_1 = 0}^t \eu^{\mat A(t - \tau_1)} \mat N u(\tau_1) \vec x_{1}(\tau_1) \, \dd \tau_1 \\
&= \int_{\tau_1 = 0}^t \eu^{\mat A(t - \tau_1)} \mat N u(\tau_1) \Bigg(\:\underbrace{
	\int_{\tau_2 = 0}^{\tau_1} \eu^{\mat A(\tau_1 - \tau_2)} \vec b u(\tau_2)\, \dd \tau_2 + \eu^{\mat A\tau_1}\vec x_0}_{\vec x_1(\tau_1)} \Bigg)
\, \dd \tau_1 \\
&= \int_{\tau_1 = 0}^t\int_{\tau_2 = 0}^{\tau_1}  \eu^{\mat A(t - \tau_1)} \mat N u(\tau_1) 
\eu^{\mat A(\tau_1 - \tau_2)} \vec b u(\tau_2) \, \dd \tau_2\dd \tau_1 \\
&+ 
\int_{\tau_1 = 0}^t  \eu^{\mat A(t - \tau_1)} \mat N u(\tau_1) \eu^{\mat A\tau_1}\vec x_0
\, \dd \tau_1\, .
\end{aligned}
\end{equation}
After $k$ successive substitutions, the solution~$\vec x_k(t)$ of the~$k$-th subsystem is:
\begin{equation}
\begin{aligned}
\vec x_k(t) &= \int_{\tau_1 = 0}^{t}\!\!\!\!\cdots\int_{\tau_k=0}^{\tau_{k-1}}\eu^{\mat A(t-\tau_1)}\mat N u(\tau_1)\eu^{\mat A(\tau_1 - \tau_2)}\mat N u(\tau_2) \cdots \mat N u(\tau_{k-1})\eu^{\mat A(\tau_{k-1} - \tau_k)}\vec b u(\tau_k)\,\dd \tau_k \cdots \dd \tau_1 \\ 
&+ \int_{\tau_1 = 0}^{t}\!\!\!\!\cdots\int_{\tau_{k-1}=0}^{\tau_{k-2}}\eu^{\mat A(t-\tau_1)}\mat N u(\tau_1)  \eu^{\mat A(\tau_1 - \tau_2)}\mat N u(\tau_2) \cdots
\mat N u(\tau_{k-1})\eu^{\mat A\tau_{k-1}}\vec x_0\,\dd \tau_{k-1} \cdots \dd \tau_1\, .
\end{aligned}
\end{equation}
In the sequel, we will substitute the integration variables often. To simplify the substitutions, we will use infinite limits and introduce heaviside step functions~$\sigma(t)$\footnote{$\sigma(t) := \left\{ \begin{array}{ll} 1 & t > 0 \\ 0 & t < 0 \end{array} \right.$. The value of the heaviside step function at $t=0$ will be discussed later.} to take the integration domain into account. The solution~$\vec x_k(t)$ of the~$k$-th subsystem for zero initial condition $\vec x_0 \!=\! \vec 0$\footnote{Without loss of generality we set~$\vec x_0 = \vec 0$, since we are especially interested in the input-output relationship, and not in the initial condition term.} is then:
\begin{equation} \label{eq:tri-1}
\begin{aligned}
\vec x_k(t) &= \int_{\tau_1=-\infty}^{\infty}\!\!\!\!\cdots\int_{\tau_k=-\infty}^{\infty}\underbrace{\eu^{\mat A(t-\tau_1)}\sigma(t-\tau_1)}_{:=\eu^{\mat A (t-\tau_1)}}\mat N \underbrace{u(\tau_1)\sigma(\tau_1)}_{:=u(\tau_1)}\eu^{\mat A(\tau_1 - \tau_2)}\sigma(\tau_1-\tau_2)\mat N u(\tau_2)\sigma(\tau_2) \cdots \\
& \qquad\qquad\qquad \times \mat N u(\tau_{k-1})\sigma(\tau_{k-1})\eu^{\mat A(\tau_{k-1} - \tau_k)}\sigma(\tau_{k-1} - \tau_k)\vec b u(\tau_k)\sigma(\tau_k)\,\dd \tau_k \cdots \dd \tau_1.
\end{aligned}
\end{equation}
The heaviside step functions do not need to be explicitly written down further, when we assume one-sided input signals~$u(t):=u(t)\sigma(t)$ and define one-sided matrix exponentials~$\eu^{\mat At} := \eu^{\mat A t}\sigma(t)$. Since the magnitude of the Jacobian determinant $\left|\det \left(\pd \tau_i/\pd \tilde{\tau}_j\right)\right| = 1$ for all sequal substitutions\footnote{The sign of the Jacobian determinant does not matter, as it cancels out with the integration limits:$\:\:\int_{a}^{b}(\cdot)\,\dd \tau \rightarrow \int_{b}^{a}(\cdot)\,(-\dd \tilde{\tau})$.}, we only need to substitute the variables of the input signal~$u(t)$ and matrix exponentials~$\eu^{\mat A t}$.

Although in the literature \cite{rugh1981nonlinear,flagg2012interpolation} the triangular kernels are usually defined directly out of equation \eqref{eq:tri-1}, here a change of integration variables is first performed, in order to obtain a triangular input-output representation where the kernels appear without reflected arguments $-\tau_1, \ldots, -\tau_k$. Making the change of variables
\begin{equation*}
\begin{aligned}
\left.\begin{array}{r@{\:=\:}lr@{\:=\:}lr@{\:=\:}l}
\tilde{\tau}_k & t - \tau_1, \quad &\tilde{\tau}_{k-1} & t - \tau_2,\:\: \ldots, &\tilde{\tau}_1 & t - \tau_k, \\[0.3em]
\tau_1 & t - \tilde{\tau}_k, \quad &\tau_2 & t - \tilde{\tau}_{k-1},\:\: \ldots, &\tau_k & t - \tilde{\tau}_1,
\end{array}\right\}
\quad \text{with} \quad
\left\{\begin{array}{rl}
\tau_1 - \tau_2 & \!\! = \ \tilde{\tau}_{k-1} - \tilde{\tau}_k, \\[-0.4em]
&\vdots \\[-0.2em]
\tau_{k-1} - \tau_k & \!\! = \ \tilde{\tau}_1 - \tilde{\tau}_2,
\end{array}\right.
\end{aligned}
\end{equation*}
and letting again~$\tilde{\tau}_k \rightarrow \tau_k$,~$\ldots$,~$\tilde{\tau}_1 \rightarrow \tau_1$, then the triangular input-output representation (SISO) is given by:
\\
\fbox{\parbox{14.7cm}{
\begin{equation}
\begin{aligned}\label{eq:y-tri-SISO}	
y_k(t) &=  \int_{\tau_1 = -\infty}^{\infty}\!\!\!\!\cdots\int_{\tau_k=-\infty}^{\infty} g^{\triangle}_{k}(\tau_1,\ldots,\tau_k)u(t-\tau_k)\cdots u(t-\tau_1)\,\dd \tau_k \cdots \dd \tau_1\, ,
\end{aligned}
\end{equation}
with the triangular kernels
\begin{equation}
\begin{aligned}\label{eq:tri-SISO-kernels}
g^{\triangle}_{k}(t_1,\ldots,t_k) \!= \!\left\{\!
\begin{array}{ll}
\vec c^{\mathsf T}\eu^{\mat A t_k}\mat N\eu^{\mat A(t_{k-1} - t_k)}\mat N \cdots \mat N\eu^{\mat A(t_1 - t_2)}\vec b, & 0 < t_k < \ldots < t_1 \\
\text{not yet defined}, & \text{on the surface} \\
0, & \text{else.} 
\end{array} \right.
\end{aligned} 
\end{equation}
}}\\[1ex]
For a better understanding of the definition domain of the triangular kernels, Figure~\ref{fig:defdtriangle_2} exemplary shows the triangular kernel of the second subsystem~$g_2^\triangle(t_1,t_2)$. Its domain is the triangle ($2$-simplex)~$t_1 > t_2 > 0$. Its surface is the line ($1$-simplex)~$t_1 = t_2 > 0$. We do not care about the bottom line~$t_2 = 0$, since we are not interested in the value of the impulse response at~$t=0$, which will become clear later. In Figure~\ref{fig:defdtriangle_3} the triangular kernel of the third subsystem~$g_3^\triangle(t_1,t_2,t_3)$ is illustrated. Its domain is the tetrahedron ($3$-simplex)~$t_1 > t_2 > t_3 > 0$. Its surface consists of the triangles ($2$-simplices)~$t_1 = t_2 > t_3 > 0$ and~$t_1 > t_2 = t_3 > 0$ and their intersection line ($1$-simplex)~$t_1 = t_2 = t_3 > 0$. Again, we do not care about the bottom triangle of the tetrahedron, where~$t_3 = 0$. This can be generalized to the $k$-th subsystem.
\begin{figure}[htp!]
	\begin{center}
		\subfloat[second subsystem]{\centering
			\begin{tikzpicture}
[xscale=1/(5 - (-0.5))*0.0352778*0.44*\textwidth,
yscale=1/(5 - (-0.5))*0.0352778*0.44*\textwidth,
xtick/.style={below,fill=white,font=\tiny},
ytick/.style={left=0.2,fill=white,font=\tiny},
xlabel/.style={below right,font=\scriptsize},
ylabel/.style={left,font=\footnotesize},
beziernode/.style={circle, inner sep=0pt,minimum size=5pt,fill=black}
]
\def\tickin{0.1}
\def\tickout{0.1}
\def\delta{0.5}
\def\xmax{5}
\def\xmin{0}
\def\ymax{0.61803*\xmax} 
\def\ymin{0}

\draw [->] (0,\ymin) -- (0,\ymax) node [ylabel] {$t_2$};
\draw [->] (\xmin,0) -- (\xmax,0) node [xlabel] {$t_1$};
\draw [dashed] (\xmin,\xmin) -- (\ymax,\ymax);
\draw[arrows={-latex}] (0.5*\xmin + 0.5*\xmax,0.5*\xmin + 0.5*\xmax) --++ (0.5,-0.5) node[xshift = 20, yshift = 15,fill=white, font=\scriptsize, align = center] {$\sigma(t_1 - t_2)$,\\$t_1 > t_2$};
\draw[arrows={-latex}] (0.8*\xmin + 0.8*\xmax,0) --++ (0,0.707) node[xshift = 20, yshift = -10,fill=white, font=\scriptsize, align = center] {$\sigma(t_2)$,\\$t_2 > 0$};
\node[beziernode] at (0.2*\ymax,0.8*\ymax) {};
\draw [ultra thin]  (0.2*\ymax,0.8*\ymax)   .. controls +(-0.2,+0.4) and ++(+0.2,-0.4) .. (-0.05*\ymax,0.5*\ymax)  node[anchor=south,font=\scriptsize, rotate = 90] {$g_2^{\triangle}(t_1,t_2) = 0$};
\draw [ultra thin] (0.5*\ymax,0.5*\ymax) .. controls +(+0.05,+0.05) and ++(+0.05,-0.05) .. (0.5*\ymax,0.6*\ymax)  node[anchor=south,font=\scriptsize, rotate = 45] {$g_2^{\triangle}(t_1,t_2)$ not yet defined};
\node[beziernode] at (0.9*\ymax,0.3*\ymax) {};
\draw [ultra thin] (0.9*\ymax,0.3*\ymax)  .. controls +(+0.2,-0.4) and ++(-0.2,+0.4) .. (0.75*\ymax,-0.05*\ymax)  node[anchor=north,font=\scriptsize] {$g^{\triangle}_2(t_1,t_2) = \vec c^{\mathsf T} \eu^{\mat A t_2}\mat N \eu^{\mat A (t_1 - t_2)}\vec b$};

\end{tikzpicture}%
			\label{fig:defdtriangle_2}}%
		\subfloat[third subsystem]{\centering
			\definecolor{mycolor1}{rgb}{0.97690,0.98390,0.08050}%
\definecolor{mycolor2}{rgb}{0.24220,0.15040,0.66030}%
\definecolor{mycolor3}{rgb}{0.09640,0.75000,0.71204}%

\definecolor{mycolor4}{rgb}{0.01, 0.68, 0.33}%

\begin{tikzpicture}[
helpline/.style={dashed},
]

\begin{axis}[
width=0.44*\textwidth,
tick label style={font=\tiny, fill=white},
label style={font=\color{white!15!black}, font=\tiny},
axis on top,
scale only axis,
tick align=outside,
axis lines=center,
y dir=reverse,
name=myplot,
view={22}{35.6},
enlargelimits=.1,
ymin=0, ymax=1, xmin=0, xmax=1, zmin=0, zmax=1,
xlabel=$t_1$, ylabel=$t_2$, zlabel=$t_3$,
ticks=none,
every axis x label/.append style={at=(ticklabel* cs:1.03)},
every axis y label/.append style={at=(ticklabel* cs:0.2)},
every axis z label/.append style={at=(ticklabel* cs:1.05)}
]

\addplot3[area legend, draw=black, fill=mycolor2, fill opacity=0.5]
table[row sep=crcr] {%
	x	y	z\\
	0	0	0\\
	1	0	0\\
	1	1	1\\
}--cycle;

\addplot3[area legend, draw=black, fill=mycolor1, fill opacity=0.5]
table[row sep=crcr] {%
	x	y	z\\
	0	0	0\\
	1	1	0\\
	1	1	1\\
}--cycle;

\draw [helpline] (1,0,0) -- (1,0.08,0);
\draw [helpline] (1,0.3,0) -- (1,1,0);

\draw [helpline] (0,1,0) -- (1,1,0);

\draw [helpline] (0,0,1) -- (1,0,1);

\draw [helpline] (0,0,1) -- (0,1,1);
\draw [helpline] (0,1,1) -- (1,1,1);
\draw [helpline] (1,0,1) -- (1,1,1);

\draw [helpline] (1,0,0) -- (1,0,1);
\draw [helpline] (0,1,0) -- (0,1,1);

\draw [draw=mycolor4, ultra thick] (0,0,0) -- (1,1,1);

\draw [ultra thin]  (0.5,0.5,0.5)   .. controls +(-0.1,-0.1,0.1) and ++(+0.1,+0.1,-0.1) .. (0.5,0.5,0.75)  node[anchor=south,font=\scriptsize,rotate=-10] {$t_1 \!= \!t_2 \!= \!t_3 \!> \!0$};

\draw [ultra thin]  (0.4,0.5,0.25)   .. controls +(-0.2,+0.2,-0.1) and ++(+0.2,-0.2,0.1) .. (0.3,0.7,0)  node[anchor=north,font=\scriptsize,rotate=-10] {$t_1 \!=\! t_2 \!>\! t_3 \!>\! 0$};

\draw [ultra thin]  (0.85,0.25,0.25)   .. controls +(-0.1,+0.15,-0.1) and ++(+0.05,-0.1,0.1) .. (0.87,0.04,0)  node[anchor=north,font=\scriptsize,rotate=-10] {$t_1\! >\! t_2 \!=\! t_3 \!>\! 0$};

\end{axis}

\end{tikzpicture}		
			\label{fig:defdtriangle_3}}		
	\end{center}
	\caption{(a) Domain of the triangular kernel $g^{\triangle}_2(t_1,t_2)$ of the second subsystem. The kernel is not yet defined on the line $t_1=t_2 > 0$. (b) Domain of the triangular kernel $g^{\triangle}_3(t_1,t_2,t_3)$ of the third subsystem. The kernel is not yet defined on the triangles $t_1 =  t_2 > t_3 > 0$ and~$t_1 > t_2 = t_3 > 0$ and its intersection line~$t_1 = t_2 = t_3 > 0$.}
	\label{fig:defdtriangle}
	\vspace{-2.5em}
\end{figure}

\subsubsection*{MIMO case}
All the concepts we have seen so far can be generalized for MIMO bilinear systems. To this end, replace in the following~$\vec c^{\mathsf T} \rightarrow \mat C$,~$\mat N u(t) \rightarrow \sum_{j = 1}^{m}\mat N_j u_j(t)$ and~$\vec b u(t) \rightarrow \sum_{j = 1}^{m}\vec b_j u_j(t)$. Thus, equations~\eqref{eq:y-tri-SISO} and~\eqref{eq:tri-SISO-kernels} become:
\begin{equation}
\begin{aligned}
\lefteqn{\vec y_k(t)
= \int_{\tau_1 = -\infty}^{\infty} \!\!\!\! \cdots \int_{\tau_k=-\infty}^{\infty} \mat C\eu^{\mat A \tau_k}
\left(\sum_{j = 1}^{m}\mat N_ju_j(t-\tau_k)\right)\eu^{\mat A(\tau_{k-1} - \tau_k)}\left(\sum_{j = 1}^{m}\mat N_ju_j(t-\tau_{k-1})\right) \cdots} \\
& \qquad\qquad\qquad\qquad \times \left(\sum_{j = 1}^{m}\mat N_j u_j(t-\tau_2)\right) \eu^{\mat A(\tau_1 - \tau_2)}\left(\sum_{j = 1}^{m}\vec b_ju_j(t-\tau_1)\right)\,\dd \tau_k \cdots \dd \tau_1 \\
& = \sum_{j_1 = 1}^{m}\cdots\sum_{j_k = 1}^{m}\int_{\tau_1 = -\infty}^{\infty} \!\!\!\! \cdots\int_{\tau_k=-\infty}^{\infty}\mat C\eu^{\mat A \tau_k}
\mat N_{j_k}u_{j_k}(t-\tau_k)\eu^{\mat A(\tau_{k-1} - \tau_k)}\mat N_{j_{k-1}}u_{j_{k-1}}(t-\tau_{k-1}) \cdots\\
& \qquad\qquad\qquad\qquad \times \mat N_{j_2}u_{j_2}(t-\tau_2) \eu^{\mat A(\tau_1 - \tau_2)}\vec b_{j_1}u_{j_1}(t-\tau_1)\,\dd \tau_k \cdots \dd \tau_1\, .
\end{aligned}
\end{equation}
The triangular input-output representation in the MIMO case is therefore:
\\
\fbox{\parbox{14.7cm}{	
\begin{equation}
	\begin{aligned}	
	\vec y(t) = \sum_{k=1}^{\infty} \vec y_k(t), && && && \vec y_k(t) = \sum_{j_1 = 1}^{m}\cdots \sum_{j_k = 1}^{m}\vec y^{(j_1,\ldots,j_k)}_k(t)\, ,
	\end{aligned}
\end{equation}
with the output components
\begin{equation}
	\begin{aligned}\label{eq:y-tri-MIMO}
	\vec y^{(j_1,\ldots,j_k)}_k(t) = \int_{\tau_1 = -\infty}^{\infty} \cdots\int_{\tau_k=-\infty}^{\infty} \vec g^{(j_1,\ldots,j_k)}_{k, \triangle}(\tau_1,\ldots,\tau_k)u_{j_k}(t-\tau_k)\cdots u_{j_1}(t-\tau_1)\,\dd \tau_k \cdots \dd \tau_1\, ,
	\end{aligned}
\end{equation}
and the triangular MIMO kernels 
\begin{equation}
	\begin{aligned} \label{eq:tri-MIMO-kernels}
	\vec g^{(j_1,\ldots,j_k)}_{k,\triangle}(t_1,\ldots,t_k) = \left\{
	\begin{array}{ll}
	\mat C \eu^{\mat A t_k}\mat N_{j_k}\eu^{\mat A(t_{k-1} - t_k)}\mat N_{j_{k-1}}\cdots \mat N_{j_2}\eu^{\mat A(t_1 - t_2)}\vec b_{j_1}, & 0 < t_k < \ldots < t_1 \\
	\text{not yet defined}, & \text{on the surface} \\ 
	\vec 0, & \text{else.} 
	\end{array} \right.
	\end{aligned} 
\end{equation}
}}
\begin{remark}[Kronecker notation for MIMO kernels]
	In the literature, the input-output representation and the kernels for the MIMO case are usually expressed by means of Kronecker products. In that sense, the kernels are matrix-valued functions of dimension $\vec g^\triangle_k(t_1,\ldots,t_k) \in \mathbb{R}^{p\times m^k}$. Note, however, that each kernel $\mat g^\triangle_k(t_1,\ldots,t_k) \in \mathbb{R}^{p\times m^k}$ includes all $m^k$ combinations of the kernels~$\vec g^{(j_1,\ldots,j_k)}_{k,\triangle}(t_1,\ldots,t_k) \in \mathbb{R}^{p}$ defined here. Although the presentation with Kronecker products might be more compact, we prefer the sum-notation, as the generalization from SISO to MIMO is straightforward by summing over all SISO-like terms~$\vec g^{(j_1,\ldots,j_k)}_{k,\triangle}(t_1,\ldots,t_k)$, whereas Kronecker products can be cumbersome. For a detailed presentation with Kronecker notation the reader is referred to \cite[\S4.2]{breiten2013interpolatory}, \cite[\S2]{flagg2015multipoint}, \cite[\S2.1.2.]{gebhart2017impulse}.
\end{remark} 

\subsubsection{Regular kernel representation} \label{subsubsec:reg-kernel-rep}
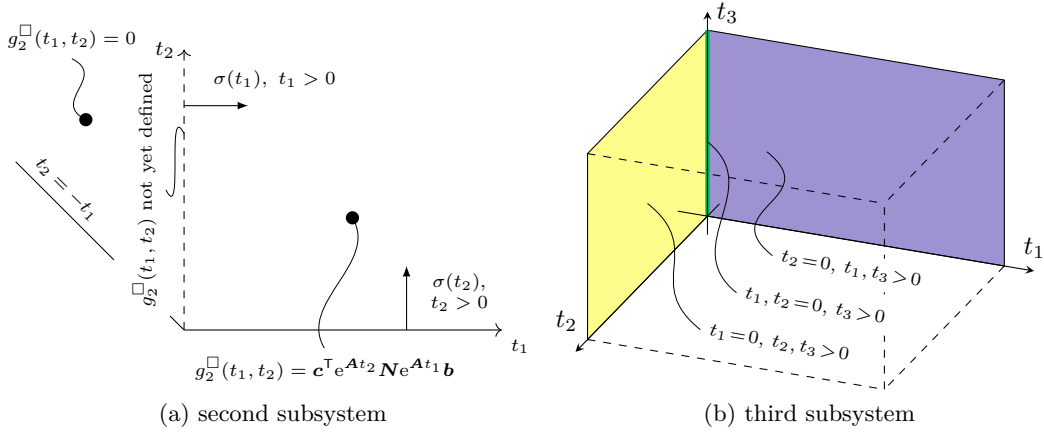
\begin{figure}
	\begin{center}
	\subfloat[second subsystem]{\centering
		\begin{tikzpicture}
[xscale=1/(5 - (-0.5))*0.0352778*0.44*\textwidth,
yscale=1/(5 - (-0.5))*0.0352778*0.44*\textwidth,
xtick/.style={below,fill=white,font=\tiny},
ytick/.style={left=0.2,fill=white,font=\tiny},
xlabel/.style={below right,font=\scriptsize},
ylabel/.style={left,font=\footnotesize},
beziernode/.style={circle, inner sep=0pt,minimum size=5pt,fill=black}
]
\def\tickin{0.1}
\def\tickout{0.1}
\def\delta{0.5}
\def\xmax{3.5}
\def\xmin{0}
\def\ymax{0.61803*5} 
\def\ymin{0}

\draw [->,dashed] (0,\ymin) -- (0,\ymax) node [ylabel] {$t_2$};
\draw [->] (\xmin,0) -- (\xmax,0) node [xlabel] {$t_1$};
\draw [-] (-0.6*\ymax,0.6*\ymax) -- (-0.25*\ymax,0.25*\ymax) node[pos=0.4, anchor=south,font=\scriptsize, rotate = -45] {$t_2 = -t_1$};
\draw [-] (-0.05*\ymax,0.05*\ymax) -- (0,0);
\draw[arrows={-latex}] (0,0.8*\ymax) --++ (+0.707,0) node[xshift = 10, yshift = 10, font=\scriptsize] {$\sigma(t_1), \ t_1 > 0$};
\draw[arrows={-latex}] (0.7*\xmax,0) --++ (0,0.707) node[xshift = 20, yshift = -10,fill=white, font=\scriptsize, align = center] {$\sigma(t_2)$,\\$t_2 > 0$};

\draw [ultra thin]  (0,0.7*\ymax)   .. controls +(-0.2,+0.4) and ++(+0.2,-0.4) .. (-0.2,0.5*\ymax)  node[anchor=south,font=\scriptsize, rotate = 90] {$g_2^{\square}(t_1,t_2)$ not yet defined};

\node[beziernode] at (-0.35*\ymax,0.75*\ymax) {};
\draw [ultra thin] (-0.35*\ymax,0.75*\ymax)  .. controls +(-0.4,+0.2) and ++(+0.4,-0.2) .. (-0.4*\ymax,-0.1+\ymax)  node[anchor=south,font=\scriptsize] {$g_2^{\square}(t_1,t_2) = 0$};
\node[beziernode] at (0.6*\ymax,0.4*\ymax) {};
\draw [ultra thin] (0.6*\ymax,0.4*\ymax)  .. controls +(+0.3,-0.5) and ++(-0.3,+0.5) .. (0.5*\ymax,-0.2)  node[anchor=north,font=\scriptsize] {$g^{\square}_2(t_1,t_2) = \vec c^{\mathsf T} \eu^{\mat A t_2}\mat N \eu^{\mat A t_1}\vec b$};

\end{tikzpicture}%
		\label{fig:defdsquare_2}}%
	\subfloat[third subsystem]{\centering
		\definecolor{mycolor1}{rgb}{0.97690,0.98390,0.08050}%
\definecolor{mycolor2}{rgb}{0.24220,0.15040,0.66030}%
\definecolor{mycolor3}{rgb}{0.09640,0.75000,0.71204}%

\definecolor{mycolor4}{rgb}{0.01, 0.68, 0.33}%

\begin{tikzpicture}[
helpline/.style={dashed},
]

\begin{axis}[
width=0.44*\textwidth,
tick label style={font=\tiny, fill=white},
label style={font=\color{white!15!black}, font=\scriptsize},
axis on top,
scale only axis,
axis lines=center,
y dir=reverse,
name=myplot,
view={22}{35.6},
enlargelimits=.1,
ymin=0, ymax=1, xmin=0, xmax=1, zmin=0, zmax=1,
xlabel=$t_1$, ylabel=$t_2$, zlabel=$t_3$,
ticks=none,
every axis x label/.append style={at=(ticklabel* cs:1.03)},
every axis y label/.append style={at=(ticklabel* cs:0.2)},
every axis z label/.append style={at=(ticklabel* cs:1.05)}
]

\addplot3[area legend, draw=black, fill=mycolor1, fill opacity=0.5]
table[row sep=crcr] {%
	x	y	z\\
	0	0	0\\
	0	1	0\\
	0	1	1\\
	0	0	1\\
}--cycle;

\addplot3[area legend, draw=black, fill=mycolor2, fill opacity=0.5]
table[row sep=crcr] {%
	x	y	z\\
	0	0	0\\
	1	0	0\\
	1	0	1\\
	0	0	1\\
}--cycle;

\draw [helpline] (1,0,0) -- (1,1,0);
\draw [helpline] (0,1,0) -- (1,1,0);
\draw [helpline] (1,0,1) -- (1,1,1);
\draw [helpline] (0,1,1) -- (1,1,1);
\draw [helpline] (1,1,0) -- (1,1,1);
\draw [draw=mycolor4, ultra thick] (0,0,0) -- (0,0,1);

\draw [ultra thin]  (0.2,0,0.4)   .. controls +(+0.2,0,-0.2) and ++(-0.2,0,0.2) .. (0.3,0.2,0)  node[anchor=west,rotate=-10,font=\scriptsize,fill=white] {$t_2\!=\!0,\, t_1, t_3\!>\!0$};

\draw [ultra thin]  (0,0,0.4)   .. controls +(+0.2,0,-0.2) and ++(-0.2,0,0.2) .. (0.3,0.5,0)  node[anchor=west,rotate=-10,font=\scriptsize,fill=white] {$t_1,t_2\!=\!0,\, t_3\!>\!0$};

\draw [ultra thin]  (0,0.5,0.4)   .. controls +(+0.2,0,-0.2) and ++(-0.2,0,0.2) .. (0.3,0.8,0)  node[anchor=west,rotate=-10,font=\scriptsize,fill=white] {$t_1\!=\!0,\, t_2, t_3\!>\!0$};

\end{axis}

\end{tikzpicture}		
		\label{fig:defdsquare_3}}		
	\end{center}
		\caption{(a) Domain of the regular kernel $g^{\square}_2(t_1,t_2)$ of the second subsystem. The kernel is not yet defined on the line $t_1=0, t_2 > 0$. (b) Domain of the regular kernel $g^{\square}_3(t_1,t_2,t_3)$ of the third subsystem. The kernel is not yet defined on the rectangles $t_1 = 0, t_2, t_3 > 0$ and~$t_2=0, t_1, t_3 > 0$ and its intersection line~$t_1 = t_2 = 0, t_3 > 0$.}
		\label{fig:defdsquare11}
\end{figure}

To simplify the terms~$\eu^{\mat A(t_{k-1} - t_{k})}$ arising in the triangular representation~\eqref{eq:tri-SISO-kernels}, we now give the input-output behavior of the bilinear system in terms of the regular kernels. For this purpose, we perform the change of variables
\begin{equation*}
\begin{aligned}
\left.\begin{array}{r@{\:=\:}lr@{\:=\:}lr@{\:=\:}l}
\tilde{\tau}_k & \tau_k, \quad &\tilde{\tau}_{k-1} & \tau_{k-1} - \tau_{k}, \:\: \ldots, &\tilde{\tau}_1 & \tau_1 - \tau_2, \\[0.3em]
\tau_k & \tilde{\tau}_k, \quad &\tau_{k-1} & \tau_{k} + \tilde{\tau}_{k-1}, \:\: \ldots, &\tau_1 & \tau_2 + \tilde{\tau}_1,
\end{array}\right\}
\quad \text{with} \quad
\left\{\begin{array}{r@{\:=\:}l}
\tau_{k-1} & \tau_k + \tilde{\tau}_{k-1} = \tilde{\tau}_k +  \tilde{\tau}_{k-1}, \\[-0.4em]
\multicolumn{2}{l}{\quad\quad \, \, \vdots} \\[-0.2em]
\tau_1 & \tau_2 + \tilde{\tau}_1 = \tilde{\tau}_k + \ldots + \tilde{\tau}_1,
\end{array}\right.
\end{aligned}
\end{equation*}
and let again~$\tilde{\tau}_k \rightarrow \tau_k,\ldots,\tilde{\tau}_1 \rightarrow \tau_1$. Then, the regular input-output representation (SISO) is:
\\
\fbox{\parbox{14.7cm}{	
\begin{equation}
	\begin{aligned} \label{eq:y-reg-SISO}
	y_k(t) &= \int_{\tau_1 = -\infty}^{\infty}\!\!\!\!\cdots\int_{\tau_k = -\infty}^{\infty} g^{\square}_{k}(\tau_1,\ldots,\tau_k)u(t - \tau_k)\cdots
	u(t - \tau_k - \ldots - \tau_1)\, \dd \tau_k \cdots \dd \tau_1\, ,
	\end{aligned}
\end{equation}
with the regular kernels
\begin{equation}
	\begin{aligned}\label{eq:reg-SISO-kernels}
	g^{\square}_{k}(t_1,\ldots,t_k) = \left\{
	\begin{array}{ll}
	\vec c^{\mathsf T} \eu^{\mat A t_k}\mat N\eu^{\mat At_{k-1}}\mat N\cdots \mat N \eu^{\mat A t_1}\vec b, & t_1,\ldots,t_k > 0 \\
	\text{not yet defined}, & \text{on the surface}\\
	0, & \text{else.} 
	\end{array} \right.
	\end{aligned} 
\end{equation}
}}\\[1ex]
The regular kernel of the second subsystem~$g_2^\square(t_1,t_2)$ is depicted in Figure~\ref{fig:defdsquare_2}. Its domain is the rectangle ($2$-dimensional hypercube)~$t_1, t_2 > 0$. Its surface  is the line ($1$-dimensional hypercube)~$t_1 = 0, t_2 > 0$. We do not care about the bottom line~$t_2 = 0$, since we are not interested in the value of the impulse response at~$t=0$, which will become clear later. In Figure~\ref{fig:defdsquare_3} the domain of the regular kernel of the third subsystem~$g_3^\square(t_1,t_2,t_3)$ is the cube ($3$-hypercube)~$t_1,t_2,t_3 > 0$. Its surface consists of the rectangles ($2$-hypercubes)~$t_1 = 0,t_2, t_3 > 0$ and~$t_2 = 0, t_1,t_3 > 0$ and their intersection line ($1$-hypercube)~$t_1 = t_2 = 0, t_3 > 0$. Again, we do not care about the bottom rectangle of the cube, where~$t_3 = 0$. This can be generalized to the $k$-th subsystem.
%
%

\subsubsection*{MIMO case}
Analog to the triangular kernels, the regular input-output representation in the MIMO case is given by:
\\
\fbox{\parbox{14.7cm}{	
\begin{equation}
	\begin{aligned}	
	\vec y(t) = \sum_{k=1}^{\infty} \vec y_k(t), && && && \vec y_k(t) = \sum_{j_1 = 1}^{m}\cdots \sum_{j_k = 1}^{m}\vec y^{(j_1,\ldots,j_k)}_k(t)\, ,
	\end{aligned}
\end{equation}
with the output components
\begin{equation}
	\begin{aligned}
	\vec y^{(j_1,\ldots,j_k)}_k(t) = \int_{\tau_1 = -\infty}^{\infty}\!\!\!\!\!\!\!\!\cdots\int_{\tau_k = -\infty}^{\infty}\!\!\!\! \vec g_{k,\square}^{(j_1,\ldots,j_k)}(\tau_1,\ldots,\tau_k)u_{j_k}(t - \!\tau_k)\cdots
	u_{j_1}(t - \!\tau_k - \ldots - \!\tau_1)\, \dd \tau_k \cdots \dd \tau_1
	\end{aligned}
\end{equation}
and the regular MIMO kernels
\begin{equation}
	\begin{aligned} \label{eq:reg-MIMO-kernels}
	\vec g^{(j_1,\ldots,j_k)}_{k,\square}(t_1,\ldots,t_k) = \left\{
	\begin{array}{ll}
	\mat C \eu^{\mat A t_k}\mat N_{j_k}\eu^{\mat At_{k-1}}\mat N_{j_{k-1}}\cdots \mat N_{j_2}\eu^{\mat A t_1}\vec b_{j_1}, & t_1,\ldots,t_k > 0 \\
	\text{not yet defined}, & \text{on the surface} \\
	0, & \text{else.} 
	\end{array} \right.
	\end{aligned} 
\end{equation}
}}\\[0.1em]

\subsection{Multidimensional Laplace transform of the kernels}
Once the Volterra kernels have been derived, we now want to analyze them in the frequency domain. 
Therefore, and similar as for single-variable kernels arising in linear systems, the Laplace transform for multivariable functions is first introduced in the following definition.

\begin{definition}[Multidimensional Laplace transform]
	A multivariable function~$f(t_1,\ldots,t_k)$ with $f\!:\! \mathbb{R}^k_{\geq 0} \rightarrow \mathbb{C}^k$ is \emph{Laplace transformable} (in terms of the $k$-dimensional Laplace transform), if the integral
	\begin{equation}
	\begin{aligned}
		F(s_1,\ldots,s_k) &:= \mathcal{L}_k\{f(t_1,\ldots,t_k) \}(s_1,\ldots,s_k) \\[0.3em]
		&:= \int_{t_1 = 0}^{\infty} \ldots \int_{t_k = 0}^{\infty} f(t_1,\ldots,t_k)\eu^{-s_1t_1}\cdots\eu^{-s_kt_k}\,\dd t_k\cdots\dd t_1\, ,\\[0.7em]
		\vec s &= (s_1 \:\:\ldots\:\: s_k)^{\mathsf T} \in H_{\gamma_1,\ldots,\gamma_k} := H_{\vec \gamma} = \left.\left\{ \vec s \in \mathbb{C}^k \right|\:\: \Re (s_i) > \gamma_i, \:\: i = 1,\ldots,k \right\}\, ,
		\end{aligned}
	\end{equation}
	converges for values~$\vec s~\in \mathbb{C}^k$ on the $k$-dimensional complex half-space~$H_{\vec \gamma}$ of the $k$-dimensional complex vector space~$\mathbb{C}^k$. The function~$F(s_1,\ldots,s_k)$ with~$F: H_{\vec \gamma}\rightarrow \mathbb{C}^k$ is called \emph{$k$-dimensional Laplace transform} of~$f(t_1,\ldots,t_k)$. 
\end{definition}
\begin{remark}[Existence and uniqueness of the Laplace transform]
	The existence of the Laplace transform is guaranteed, if~$f(t_1,\ldots,t_k)$ is piecewise continuous on~$\mathbb{R}^k$ and of exponential order in each variable $t_1,\ldots,t_k$~\cite[\S2.1]{rugh1981nonlinear}. The~$k$-dimensional Laplace transform is unique except at finitely many -- at most $(k-1)$-dimensional -- discontinuities that are not taken into account through the integral.
\end{remark}


\subsubsection{Laplace transform of triangular kernels}
Based on equation~\eqref{eq:tri-SISO-kernels}, the~$k$-dimensional Laplace transform of the triangular kernels~$g^{\triangle}_k(t_1,\ldots,t_k)$ is given by:
\begin{equation}
\begin{aligned}
G^{\triangle}_k(s_1,\ldots,s_k) &:= \mathcal{L}_k \{g^{\triangle}_k(t_1,\ldots,t_k) \} (s_1,\ldots,s_k) \\[0.3em]
&= \int_{t_1=-\infty}^{\infty} \!\!\!\! \cdots\int_{t_k=-\infty}^{\infty} 
\underbrace{\vec c^{\mathsf T} \eu^{\mat A t_k}\mat N\cdots\mat N\eu^{\mat A (t_2 - t_3)}\mat N\eu^{\mat A (t_1 - t_2)}\vec b}_{g^{\triangle}_k(t_1,\ldots,t_k)} \\[0.3em]
&\qquad\qquad\qquad\qquad \times \eu^{-s_1t_1}\cdots\eu^{-s_{k-1}t_{k-1}}\eu^{-s_kt_k}\sigma(t_1)\cdots\sigma(t_k) \,\dd t_k \cdots \dd t_1\, .
\end{aligned}
\end{equation}
Performing the change of variables:
\begin{equation*}
\begin{aligned}
\left.\begin{array}{r@{\:=\:}lr@{\:=\:}lr@{\:=\:}l}
\!\!\!\! \tilde{t}_k & t_k, \quad &\tilde{t}_{k-1} & t_{k-1} - t_k, \:\: \ldots, &\tilde{t}_1 & t_1 - t_2, \\[0.3em]
\!\!\!\! t_k & \tilde{t}_k, \quad &t_{k-1} & \tilde{t}_{k-1} + t_k, \:\: \ldots, &t_1 & \tilde{t}_1 + t_2,
\end{array}\!\right\}
\:\: \text{with} \:\:
\left\{\!\begin{array}{r@{\:=\:}l}
t_{k-1} & \tilde{t}_{k-1} + t_k, \\
t_{k-2} & \tilde{t}_{k-2} + t_{k-1} = \tilde{t}_{k-2} + \tilde{t}_{k-1} + t_k, \\[-0.4em]
\multicolumn{2}{l}{\quad\quad \, \, \vdots} \\[-0.2em]
t_1 &\tilde{t}_1 + t_2 = \tilde{t}_1 +  \ldots + \tilde{t}_{k-1} +  t_k,
\end{array}\right.
\end{aligned}
\end{equation*}
and letting again~$\tilde{t}_k \rightarrow t_k$,~$\ldots$,~$\tilde{t}_1 \rightarrow t_1$ yields the~$k$-dimensional Laplace transform of the triangular kernels:
\\
\fbox{\parbox{14.7cm}{
\begin{equation}
\begin{aligned}
\lefteqn{G^{\triangle}_k(s_1,\ldots,s_k) = \int_{t_1=-\infty}^{\infty}\!\!\!\!\cdots\int_{t_k=-\infty}^{\infty}\vec c^{\mathsf T} \eu^{\mat A t_k}\mat N\cdots\mat N\eu^{\mat A t_{2}}\mat N\eu^{\mat A t_1}\vec b} \\[0.3em]
&\qquad\qquad\qquad\qquad \times\eu^{-s_1(t_1 + \ldots + t_k)}\eu^{-s_2(t_2 + \ldots + t_k)}\cdots\eu^{-s_kt_k}\sigma(t_{1} + \ldots + t_k)\cdots\sigma(t_k) \,\dd t_k\cdots \dd t_1 \\[0.3em]
&= \vec c^{\mathsf T} \int_{t_k = 0}^{\infty}\eu^{\big(\mat A - (s_1 + \ldots + s_k)\matn I\big)t_k}\, \dd t_k\; \mat N\cdots\mat N \int_{t_2 = 0}^{\infty}\eu^{\big(\mat A - (s_1 + s_2)\matn I\big)t_2}\, \dd t_2 \mat N \int_{t_1 = 0}^{\infty}\eu^{\big(\mat A - s_1 \matn I\big)t_1}\, \dd t_1\; \vec b \\[0.3em]
&= \vec c^{\mathsf T} \big( (s_1 + \ldots + s_k)\matn I - \mat A \big)^{-1}\mat N \cdots\mat N\big((s_1 + s_2)\matn I - \mat A \big)^{-1}\mat N\big(s_1 \matn I - \mat A \big)^{-1}\vec b\, .
\end{aligned}
\end{equation}
}}\\[1ex]
The~$k$-dimensional Laplace transform~$G^\triangle_k(s_1,\ldots,s_k)$ of the triangular kernels exists in the region~$\Re (s_1)> \mathrm{max}\: \Re (\flambda(\mat A))$,~$\ldots$,~$\Re (s_1 + \ldots + s_k) > \mathrm{max}\: \Re (\flambda(\mat A))$, with the real part~$\Re \flambda(\mat A)$ of the eigenvalues~$\flambda(\mat A)$ of the dynamic matrix~$\mat A$.
\subsubsection*{MIMO case}
Starting from the triangular MIMO kernels stated in~\eqref{eq:tri-MIMO-kernels}, the~$(j_1,\ldots,j_k)$-th transfer function~$\vec G^{(j_1,\ldots,j_k)}_{k,\triangle}(s_1,\ldots,s_k) \in \mathbb{C}^p$ can be given as:
\\
\fbox{\parbox{14.7cm}{
\begin{equation}
\begin{aligned}
\lefteqn{\vec G_{k,\triangle}^{(j_1,\ldots,j_k)}(s_1,\ldots,s_k)} \\
&= \mat C \big( (s_1 + \ldots + s_k)\matn I - \mat A \big)^{-1}\mat N_{j_k} \cdots \mat N_{j_3}\big((s_1 + s_2)\matn I - \mat A \big)^{-1}\mat N_{j_2}\big(s_1 \matn I - \mat A \big)^{-1}\vec b_{j_1}.
\end{aligned}
\end{equation}
}}
\begin{remark}[Kronecker notation for MIMO transfer functions]
	Similar as with the kernels, in the literature the MIMO transfer functions are usually given in Kronecker notation. Note, however, that each transfer function matrix~$\mat G^{\triangle}_k(s_1,\ldots,s_k) \in \mathbb{R}^{p\times m^k}$ includes all $m^k$ combinations of the~$(j_1,\ldots,j_k)$-th transfer functions~$\vec G^{(j_1,\ldots,j_k)}_{k,\triangle}(s_1,\ldots,s_k) \in \mathbb{R}^{p}$ defined here. 
\end{remark}
Please note that the triangular transfer functions~$G^{\triangle}_k(s_1,\ldots,s_k)$ yield factors~$\big( (s_1 + \ldots + s_k)\matn I - \mat A \big)^{-1}$, where the frequency variables~$s_1,\ldots,s_k$ are summed up and do not appear as independent, separable variables. In other words, each singularity term~$\big( (s_1 + \ldots + s_k)\matn I - \mat A \big)^{-1}$ of the transfer function is not expressed by means of a single variable~$s_k$, but rather by several variables. This complicates the application of the residue calculus and the inverse Laplace transform in each variable independently. For that reason, the triangular transfer functions~$G^{\triangle}_k(s_1,\ldots,s_k)$ are rarely used as a starting point for model order reduction, but the regular transfer functions (see section~\ref{subsubsec:reg-TFs}) have established instead.
\vspace{0.3em}
\subsubsection{Laplace transform of regular kernels} \label{subsubsec:reg-TFs}
Next, the regular kernels~$g^{\square}_k(t_1,\ldots,t_k)$ are transformed into the frequency domain. Based on equation~\eqref{eq:reg-SISO-kernels},
the multidimensional Laplace transform of the regular kernels is given by:\nopagebreak
\\
\nopagebreak
\fbox{\parbox{14.7cm}{
\begin{equation}
\begin{aligned} \label{eq:tf-reg-SISO}
G^{\square}_k(s_1,\ldots,s_k) &:= \mathcal{L}_k \{g^{\square}_k(t_1,\ldots,t_k) \} (s_1,\ldots,s_k) \\[0.3em]
&= \int_{t_1 = 0}^{\infty} \cdots \int_{t_k = 0}^{\infty}
\underbrace{\vec c^{\mathsf T} \eu^{\mat A t_k}\mat N\cdots\mat N\eu^{\mat A t_1}\vec b}_{g^{\square}_k(t_1,\ldots,t_k)} \, \eu^{-s_1t_1}\cdots\eu^{-s_kt_k} 
\,\dd t_k \cdots \dd t_1 \\[0.3em]
&=\vec c^{\mathsf T} \int_{t_k = 0}^{\infty}\eu^{(\mat A - s_k \matn I)t_k}\,\dd t_k\, \mat N \cdots \mat N \int_{t_1 = 0}^{\infty}\eu^{(\mat A - s_1 \matn I)t_1}\,\dd t_1\, \vec b \\[0.3em]
&= \vec c^{\mathsf T} (s_k \matn I - \mat A)^{-1}\mat N\cdots \mat N(s_1 \matn I - \mat A)^{-1}\vec b \, .
\end{aligned}
\end{equation}
}}\\[1ex]
The~$k$-dimensional Laplace transform~$G^\square_k(s_1,\ldots,s_k)$ of the regular kernels exists in the region $\Re (s_1), \ldots, \Re (s_k) > \mathrm{max}\: \Re (\flambda(\mat A))$.

The main feature of the regular transfer functions~$G^{\square}_k(s_1,\ldots,s_k)$ -- in comparison to the triangular ones -- is, that they are expressed in terms of a product of \emph{single-variable} factors~$(s_k \matn I - \mat A)^{-1}$. This permits to write the transfer functions as (see~\cite[p. 73, \S4.2]{rugh1981nonlinear}, \cite[\S2.2]{flagg2015multipoint})
\begin{equation*}
G^{\square}_k(s_1,s_2,\ldots,s_k) = \frac{P(s_1,s_2,\ldots,s_k)}{Q(s_1,s_2,\ldots,s_k)} = \frac{P(s_1,s_2,\ldots,s_k)}{Q(s_1) Q(s_2) \cdots Q(s_k)}\, ,
\end{equation*}
where the denominator~$Q(s_1,s_2,\ldots,s_k)$ is a $k$-variate polynomial with maximum degree $kn$ that can be factored into the product of simple factors $Q(s_{\ell})$ with~$\ell=1,\ldots,k$. Each $Q(s_{\ell})=\det(s_{\ell} \matn I - \mat A)$ is a polynomial of degree $n$ in the single variable $s_{\ell}$, having each $n$ distinct roots at $\lambda_1, \lambda_2, \ldots, \lambda_n \in \mathbb{C}$. The fact that the polar sets of $G^{\square}_k(s_1,\ldots,s_k)$ are separable into $(k-1)$-dimensional hyperplanes permits the straightforward application of the residue calculus in each variable separately. Furthermore, the representation~\eqref{eq:tf-reg-SISO} as a product of factors~$(s_k \matn I - \mat A)^{-1}$ allows to perform independent single-variable inverse Laplace transforms when applying the multidimensional inverse Laplace transform~\cite[pp. 59-60, 73]{rugh1981nonlinear}, \cite[\S2.4]{flagg2012interpolation}.
\vspace{0.3em}
\subsubsection*{MIMO case}
Corresponding to the regular MIMO kernels from~\eqref{eq:reg-MIMO-kernels}, the~$(j_1,\ldots,j_k)$-th transfer function~$\vec G^{(j_1,\ldots,j_k)}_{k,\square}(s_1,\ldots,s_k) \in \mathbb{C}^p$ is:
\\
\fbox{\parbox{14.7cm}{
\begin{equation}
\begin{aligned}
\vec G_{k,\square}^{(j_1,\ldots,j_k)}(s_1,\ldots,s_k) &= \mat C (s_k \matn I - \mat A)^{-1}\mat N_{j_k} \cdots \mat N_{j_3}(s_2 \matn I - \mat A)^{-1}\mat N_{j_2} 
(s_1 \matn I - \mat A)^{-1}\vec b_{j_1} \in \mathbb{C}^{p}\, .
\end{aligned}
\end{equation}
}}\\
\newpage
\section{Input-output representation in the frequency domain}\label{sec:inp-out-freq-dom}
In the previous section, we derived the Laplace transform of the triangular kernel~$g_k^\triangle(t_1,\ldots,t_k)$ and regular kernel $g^\square_k(t_1,\ldots,t_k)$: $G_k^\triangle(s_1,\ldots,s_k)$ and~$G_k^\square(s_1,\ldots,s_k)$. Instead of merely transforming the kernels, the output $y_k(t)$ of the bilinear system -- i.e. equations~\eqref{eq:y-tri-SISO} and~\eqref{eq:y-reg-SISO} -- should be rather transformed. Thus, our aim is now to obtain an input-output representation in the frequency domain similar to the one known for linear systems~$Y(s) = G(s) \, U(s)$ using the convolution property of the multidimensional Laplace transform. Then,~$G^{\triangle}_k(s_1,\ldots,s_k)$ and~$G^{\square}_k(s_1,\ldots,s_k)$ can be interpreted as multidimensional triangular and regular transfer functions.
%
\vspace{-0.5em}
\subsection{Triangular transfer function representation}
First, the triangular input-output representation in the frequency domain is derived. To make the derivation more clear, we start with the second subsystem and generalize then the results for the~$k$-th subsystem. 

To apply the two-dimensional Laplace transform to~\eqref{eq:y-tri-SISO}, let
\vspace{-0.2em}
\begin{equation}
\begin{aligned}
\label{eq:sy_mu_mu2_1}
y^{\triangle}_2(t_1,t_2)
&:= \int_{\tau_1 = -\infty}^{\infty}\int_{\tau_2= - \infty}^{\infty} g^{\triangle}_2(\tau_1,\tau_2)u(t_2 - \tau_2)u(t_1 - \tau_1)\,\dd\tau_2\dd\tau_1\, ,\\[0.2em]
&:= \int_{\tau_1 = -\infty}^{\infty}\int_{\tau_2= - \infty}^{\infty} g^{\square}_2(\tau_1,\tau_2)u(t_2 - \tau_2)u(t_1 - \tau_2 - \tau_1)\,\dd\tau_2\dd\tau_1\, 
\end{aligned}
\end{equation}
\vspace{-0.2em}
be the triangular \emph{auxiliary} output of the second subsystem, which depends on the two time variables~$t_1$ and $t_2$. Note that the auxiliary output~$y^{\triangle}_2(t_1,t_2)$ can be equivalently described by either the triangular or regular kernels, since they can be transformed into each other according to  section~\ref{subsubsec:reg-kernel-rep}. The output of the second subsystem can be then given in terms of the triangular auxiliary output:~$y_2(t) = y^{\triangle}_2(t_1 = t, t_2 = t)$.

The Laplace transform~$Y^{\triangle}_2(s_1,s_2)$ of the auxiliary output~$y^{\triangle}_2(t_1,t_2)$ is obtained by applying the convolution property in a similar manner as it is generally done for the output~$y_1(t_1)$ of the first subsystem. In a first step, all the integral limits are lowered and raised to~$(-\infty, \infty)$ through heaviside step functions~$\sigma(t)$:
\vspace{-0.2em}
\begin{equation}
\begin{aligned}
\label{eq:sy_mu_mu2_2}
\lefteqn{Y^{\triangle}_2(s_1,s_2) := \mathcal{L}_2\{y^{\triangle}_2(t_1,t_2)\}(s_1,s_2)} \\
&= \int \smo{\int_{t_1,t_2,\tau_1,\tau_2 = -\infty}} \int \int\limits^{\infty}\vec c^{\mathsf T} \eu^{\mat A\tau_2}\mat N \eu^{\mat A \tau_1}\vec b u(t_2 - \tau_2)u(t_1 - \tau_2 - \tau_1)\eu^{-s_1t_1}\eu^{-s_2t_2} \, \sigma(t_1)\sigma(t_2) \, \dd \tau_2\dd\tau_1\dd t_2 \dd t_1\, . \!\!\!\!\!\!\!
\end{aligned}
\end{equation}
\vspace{-0.2em}
The substitution of~$t_2, t_1$ by~$\tilde{\tau}_2,\tilde{\tau}_1$:
\begin{equation*}
\begin{aligned}
\left.\begin{array}{r@{\:=\:}l}
\tilde{\tau}_2 & t_2 - \tau_2, \\ 
\tilde{\tau}_1 & t_1 - \tau_2 - \tau_1,
\end{array}\right\}
\qquad \text{with} \qquad
\left\{\begin{array}{r@{\:=\:}l}
t_2 & \tilde{\tau}_2 + \tau_2, \\
t_1 & \tilde{\tau}_1 + \tau_2 + \tau_1,
\end{array}\right.	
\end{aligned}
\end{equation*}
finally yields the Laplace transform~$Y^{\triangle}_2(s_1,s_2)$ of the triangular auxiliary output~$y^{\triangle}_2(t_1,t_2)$\footnote{Note that~$\sigma(\tilde{\tau}_1 + \tau_2 + \tau_1)\sigma(\tilde{\tau}_2 + \tau_2) = 1$ for~$\tau_1,\tau_2,\tilde{\tau}_1,\tilde{\tau}_2 > 0$.}:
\begin{equation}
\begin{aligned} \label{eq:Y2-tri}
Y^{\triangle}_2(s_1,s_2)
&= \int \smo{\int_{\tilde{\tau}_1,\tilde{\tau}_2,\tau_1,\tau_2 = -\infty}} \int \int\limits^{\infty}\vec c^{\mathsf T} \eu^{\mat A\tau_2}\mat N \eu^{\mat A \tau_1}\vec b u(\tilde{\tau}_2)u(\tilde{\tau}_1)\eu^{-s_1(\tilde{\tau}_1 + \tau_2 + \tau_1)}\eu^{-s_2(\tilde{\tau}_2 + \tau_2)} \\
&\qquad\times\sigma(\tilde{\tau}_1 + \tau_2 + \tau_1)\sigma(\tilde{\tau}_2 + \tau_2) \,\dd \tau_2\dd\tau_1\dd \tilde{\tau}_2 \dd \tilde{\tau}_1\, , \\
&= \int_{\tau_1 = 0}^{\infty}\int_{\tau_2 = 0}^{\infty}\vec c^{\mathsf T} \eu^{(\mat A - (s_1 + s_2)\matn I)\tau_2}\mat N \eu^{(\mat A - s_1\matn I)\tau_1}\vec b \,\dd \tau_2 \dd \tau_1\\
&\qquad\qquad\times
\underbrace{\int_{\tilde{\tau}_2 = 0}^{\infty}u(\tilde{\tau}_2)\eu^{-s_2\tilde{\tau}_2}\,\dd \tilde{\tau}_2}_{U(s_2)}
\underbrace{\int_{\tilde{\tau}_1 =0}^{\infty}u(\tilde{\tau}_1)\eu^{-s_1\tilde{\tau}_1}\,\dd \tilde{\tau}_1}_{U(s_1)}\, , \\
&= \vec c^{\mathsf T}\big((s_1 + s_2)\matn I - \mat A\big)^{-1}\mat N\big(s_1 \matn I - \mat A\big)^{-1}\vec b\: U(s_2)U(s_1)\, .
\end{aligned}
\end{equation}

To apply the~$k$-dimensional Laplace transform, let
\\
\fbox{\parbox{14.7cm}{
\begin{equation}
	\begin{aligned} \label{eq:y_k-tri-aux-SISO}
		y^{\triangle}_k(t_1,\ldots,t_k)
		&:= \int_{\tau_1 = -\infty}^{\infty}\!\!\!\!\cdots\int_{\tau_k= - \infty}^{\infty}
		g^{\triangle}_k(\tau_1,\ldots,\tau_k)u(t_k - \tau_k)\cdots u(t_1 - \tau_1)\,\dd\tau_k\cdots\dd\tau_1\, , \\[0.2em]
		&:= \int_{\tau_1 = -\infty}^{\infty}\!\!\!\!\cdots\int_{\tau_k= - \infty}^{\infty}
		g^{\square}_k(\tau_1,\ldots,\tau_k)u(t_k - \tau_k)u(t_{k-1} - \tau_k - \tau_{k-1})\cdots \\
		&\qquad \qquad \qquad \qquad \qquad \times u(t_1 - \tau_k - \ldots - \tau_1)
		\,\dd\tau_k\cdots\dd\tau_1\, ,
	\end{aligned}
\end{equation}
}}\\[1ex]
be the triangular auxiliary output of the~$k$-th subsystem, which depends on the time variables $t_1,\ldots,t_k$. The output of the~$k$-th subsystem can be then given in terms of the triangular auxiliary output:~$y_k(t) = y^{\triangle}_k(t_1 = t, \ldots, t_k = t)$.

Again, the convolution property of the multivariable Laplace transform is applied and infinite integral limits are used:
\begin{equation}
\begin{aligned}
\lefteqn{Y^{\triangle}_k(s_1,\ldots,s_k) := \mathcal{L}_k\{y^{\triangle}_k(t_1,\ldots,t_k)\}(s_1,\ldots,s_k)} \\
&=
\smo{\int_{\ \ t_1,\cdots,t_k = -\infty}} \cdots \int\limits^{\infty}  \int \cdots \smo{\int_{\!\!\!\!\!\!\! \tau_1,\cdots,\tau_k = -\infty}}^{\infty} \ \
\vec c^{\mathsf T} \eu^{\mat A\tau_k}\mat N \cdots \mat N\eu^{\mat A \tau_1}\vec b \, u(t_k - \tau_k)\cdots u(t_1 - \tau_k - \ldots - \tau_1) \\
&\qquad \qquad \qquad \qquad \qquad \times\eu^{-s_1t_1}\cdots \eu^{-s_kt_k} \, \sigma(t_1)\cdots\sigma(t_k)
\,\dd \tau_k\cdots\dd\tau_1\dd t_k\cdots\dd t_1\, .
\end{aligned}
\end{equation}
The substitution of~$t_k,\ldots,t_1$ by~$\tilde{\tau}_k,\ldots,\tilde{\tau}_1$:
\begin{equation}
\begin{aligned}
\left.\begin{array}{r@{\:=\:}l}
\tilde{\tau}_k & t_k - \tau_k, \\ 
\tilde{\tau}_{k-1} & t_{k-1} - \tau_k - \tau_{k-1}, \\  
\multicolumn{2}{l}{\quad\quad \, \, \vdots} \\
\tilde{\tau}_{1} & t_{1} - \tau_k - \ldots - \tau_{1}, \\
\end{array}\right\}
\qquad \text{with} \qquad
\left\{\begin{array}{r@{\:=\:}l}
t_k & \tilde{\tau}_k + \tau_k, \\
t_{k-1} & \tilde{\tau}_{k-1} + \tau_k + \tau_{k-1}, \\
\multicolumn{2}{l}{\quad\quad \, \, \vdots} \\
t_{1} & \tilde{\tau}_1 + \tau_k + \ldots + \tau_{1},
\end{array}\right.		
\end{aligned}
\end{equation}
finally yields\footnote{Note that~$\sigma(\tilde{\tau}_1 + \tau_k + \ldots + \tau_1)\cdots\sigma(\tilde{\tau}_k + \tau_k) = 1$ for~$\tau_1,\ldots,\tau_k,\tilde{\tau}_1,\ldots,\tilde{\tau}_k > 0$.}:
\begin{equation}
\begin{aligned} 
\lefteqn{Y^{\triangle}_k(s_1,\ldots,s_k) = \smo{\int_{\ \ \tilde{\tau}_1,\cdots,\tilde{\tau}_k = -\infty}} \cdots \int\limits^{\infty} \ \int \cdots \smo{\int_{\!\!\!\!\!\!\!\!\!\!\!\! \tau_1,\cdots,\tau_k = -\infty}}^{\infty} \ \ \ \vec c^{\mathsf T} \eu^{\mat A\tau_k}\mat N \cdots \mat N\eu^{\mat A \tau_1}\vec b  u(\tilde{\tau}_k)\cdots u(\tilde{\tau}_1)} \\[0.2em]
&\quad \times\eu^{-s_1(\tilde{\tau}_1 + \tau_k + \ldots + \tau_1)}\cdots \eu^{-s_k(\tilde{\tau}_k + \tau_k)} \, \sigma(\tilde{\tau}_1 + \tau_k + \ldots + \tau_1)\cdots\sigma(\tilde{\tau}_k + \tau_k)
\,\dd \tau_k\cdots\dd\tau_1 \dd \tilde{\tau}_k \cdots \dd \tilde{\tau}_1\, , \!\! \\[0.2em]
&= \int_{\tau_1 = 0}^{\infty} \cdots \int_{\tau_k = 0}^{\infty} 
\vec c^{\mathsf T} \eu^{\big(\mat A - (s_1 + \ldots + s_k)\matn I\big)\tau_k}\mat N \cdots \mat N\eu^{(\mat A - s_1\matn I)\tau_1}\vec b \,\dd \tau_k \cdots \dd \tau_1 \\[0.2em]
&\qquad\qquad\qquad \times
\underbrace{\int_{\tilde{\tau}_k = 0}^{\infty}u(\tilde{\tau}_k)\eu^{-s_k\tilde{\tau}_k}\,\dd \tilde{\tau}_k}_{U(s_k)} \: \cdots \:
\underbrace{\int_{\tilde{\tau}_1 = 0}^{\infty}u(\tilde{\tau}_1)\eu^{-s_1\tilde{\tau}_1}\,\dd \tilde{\tau}_1}_{U(s_1)}\, .
\end{aligned}
\end{equation}
Thus, the Laplace transform~$Y^{\triangle}_k(s_1,\ldots,s_k)$ of the triangular output~$y^{\triangle}_k(t_1,\ldots,t_k)$ is:
\\[1ex]
\fbox{\parbox{14.7cm}{
\begin{equation}
\begin{aligned} \label{eq:Y-tri-SISO}
Y^{\triangle}_k(s_1,\ldots,s_k) = \underbrace{\vec c^{\mathsf T}\big((s_1 + \ldots + s_k)\matn I - \mat A\big)^{-1}\mat N \cdots \mat N \big(s_1 \matn I - \mat A\big)^{-1}\vec b}_{G_k^{\triangle}(s_1,\ldots,s_k)}\: U(s_k)\cdots U(s_1)
\end{aligned}
\end{equation}
}}\\[1ex]
\begin{remark}[Response computation from the triangular frequency domain representation]
The triangular input-output representation from equation~\eqref{eq:Y-tri-SISO} is given in terms of the triangular transfer functions~$G_k^{\triangle}(s_1,\ldots,s_k)$. If the transfer functions of all subsystems and the Laplace transform~$U(s)$ of the input are known, then~$y^{\triangle}_k(t_1,\ldots,t_k)$ can be computed via the multivariable inverse Laplace transform of~$Y^{\triangle}_k(s_1,\ldots,s_k)$. The actual (single-variable) output of the bilinear system is finally calculated from $y_k(t) = y^\triangle_k(t_1 = t,\ldots,t_k = t)$.

Alternatively, one could compute~$Y_k(s)$ via the associated Laplace transform of~$Y^{\triangle}_k(s_1,\ldots,s_k)$ and then perform single-variable inverse Laplace transforms on~$Y_k(s)$ to calculate the outputs~$y_k(t)$ \cite[\S2.3]{rugh1981nonlinear}. The associated Laplace transform avoids to perform the multivariable inverse Laplace transform of~$Y^{\triangle}_k(s_1,\ldots,s_k)$ and implicitly sets~$t_1 = \cdots = t_k$ for $y_k(t)$.

\end{remark}

\subsubsection*{MIMO case}
Similar to the SISO case, the triangular auxiliary output~$\vec y^{\triangle}_k(t_1, \ldots, t_k)$ can be given in terms of both the triangular and regular MIMO kernels:
\begin{equation}
\begin{aligned} \label{eq:y(t1,...,tk)-tri-MIMO}
\lefteqn{\vec y^{\triangle}_k(t_1, \ldots, t_k) = \sum_{j_1 = 1}^{m} \! \cdots \! \sum_{j_k = 1}^{m} \, \int \cdots \smo{\int_{\!\!\!\!\!\!\!\! \tau_1,\cdots,\tau_k = -\infty}}^{\infty} \vec g^{(j_1,\ldots,j_k)}_{k, \triangle}(\tau_1,\ldots,\tau_k)u_{j_k}(t_k-\tau_k)\cdots u_{j_1}(t_1-\tau_1)\,\dd \tau_k \cdots \dd \tau_1\, ,} \\
&= \sum_{j_1 = 1}^{m} \! \cdots \! \sum_{j_k = 1}^{m} \, \int \cdots \smo{\int_{\!\!\!\!\!\!\!\! \tau_1,\cdots,\tau_k = -\infty}}^{\infty} \vec g^{(j_1,\ldots,j_k)}_{k, \square}(\tau_1,\ldots,\tau_k)u_{j_k}(t_k-\tau_k)\cdots u_{j_1}(t_1-\tau_k-\ldots-\tau_1)\,\dd \tau_k \cdots \dd \tau_1\, .
\end{aligned}
\end{equation}
Due to the linearity property of the Laplace transform,~$\vec Y^{\triangle}_k(s_1, \ldots, s_k) \in \mathbb{C}^p$ is given in sum-notation by:\\[1ex]
\fbox{\parbox{14.7cm}{
\begin{equation}
\begin{aligned} \label{eq:Y_k-tri-MIMO}
\lefteqn{\vec Y^{\triangle}_k(s_1, \ldots, s_k) =  \sum_{j_1 = 1}^{m}\cdots \sum_{j_k = 1}^{m} \vec Y^{(j_1,\ldots,j_k)}_{k,\triangle}(s_1, \ldots, s_k)} \\
&= \sum_{j_1 = 1}^{m} \cdots \sum_{j_k = 1}^{m} \underbrace{\mat C\big((s_1 + \ldots + s_k)\matn I - \mat A\big)^{-1} \mat N_{j_k} \cdots \mat N_{j_2} \big(s_1 \matn I - \mat A\big)^{-1}\vec b_{j_1}}_{\vec G_{k,\triangle}^{(j_1,\ldots,j_k)}(s_1,\ldots,s_k)} \: U_{j_k}(s_k) \cdots U_{j_1}(s_1)\, .
\end{aligned}
\end{equation}
}}

\subsection{Regular transfer function representation}
Now our aim is to derive an input-output representation in the frequency domain by means of the regular transfer functions~$G^{\square}_k(s_1,\ldots,s_k)$ instead of the triangular transfer functions~$G^{\triangle}_k(s_1,\ldots,s_k)$ as in~\eqref{eq:Y-tri-SISO}. The regular representation will be derived reversely starting from the triangular representation, since we do not know yet, how the regular auxiliary outputs~$y^{\square}_k(t_1,\ldots,t_k)$ look like.

The Laplace transform~$Y^{\triangle}_2(s_1,s_2)$ of the triangular output~$y^{\triangle}_2(t_1,t_2)$ was given by means of~$G_2^{\triangle}(s_1,s_2)$ in equation~\eqref{eq:Y2-tri}.
To come up with a formula in terms of~$G_2^{\square}(s_1,s_2)$, we first perform the substitution~$\tilde{s}_1 = s_1$,~$\tilde{s}_2 = s_1 + s_2$ with $s_2 = \tilde{s}_2 - \tilde{s}_1$ and let again~$\tilde{s}_1 \rightarrow s_1$,~$\tilde{s}_2 \rightarrow s_2$:
\begin{equation}
\begin{aligned}
Y^{\square}_2(s_1,s_2) = \underbrace{\vec c^{\mathsf T}\big(s_2\matn I - \mat A\big)^{-1}\mat N (s_1\matn I - \mat A)^{-1}\vec b}_{G_2^{\square}(s_1,s_2)} \, U(s_1)U(s_2 - s_1)\, .
\end{aligned}
\end{equation}
In order to derive the regular output~$y^{\square}_2(t_1,t_2)$, the convolution is then applied backwards:
\begin{equation}
\begin{aligned}
\lefteqn{Y^{\square}_2(s_1,s_2) = \underbrace{\int_{\tau_2 = 0}^{\infty}\int_{\tau_1 = 0}^{\infty}\vec c^{\mathsf T}\eu^{\mat A \tau_2}\mat N \eu^{\mat A \tau_1}\vec b \eu^{-s_2 \tau_2}\eu^{-s_1\tau_1}\,\dd\tau_1\dd\tau_2}_{G_2^{\square}(s_1,s_2)}} \\
&\qquad\qquad\qquad\times
\underbrace{\int_{\tilde{\tau}_1=0}^{\infty}u(\tilde{\tau}_1)\eu^{-s_1\tilde{\tau}_1}\,\dd\tilde{\tau}_1}_{U(s_1)}
\underbrace{\int_{\tilde{\tau}_2=0}^{\infty}u(\tilde{\tau}_2)\eu^{-(s_2 - s_1)\tilde{\tau}_2}\,\dd\tilde{\tau}_2}_{U(s_2 - s_1)} \\
&= \int \smo{\int_{\tilde{\tau}_1,\tilde{\tau}_2,\tau_1,\tau_2 = -\infty}} \int \int\limits^{\infty}\underbrace{\vec c^{\mathsf T} \eu^{\mat A\tau_2}\mat N \eu^{\mat A \tau_1}\vec b}_{g_2^{\square}(\tau_1,\tau_2)} u(\tilde{\tau}_1)u(\tilde{\tau}_2)\eu^{-s_1\scriptsize{\overbrace{(\tau_1 + \tilde{\tau}_1 - \tilde{\tau}_2)}^{t_1}}}\eu^{-s_2\scriptsize{\overbrace{(\tau_2 + \tilde{\tau}_2)}^{t_2}}}\sigma(\ast)\,\dd \tau_2\dd \tau_1\dd \tilde{\tau}_2\dd \tilde{\tau}_1\, .
\end{aligned}
\end{equation}
The heaviside step function~$\sigma(\ast)$ will later define the integral limits for the Laplace transform and should be such that~$\sigma(\ast) = 1$ for~$\tau_1,\tau_2,\tilde{\tau_1},\tilde{\tau}_2 > 0$. To finally come up with well-known~$\eu^{-s_kt_k}$-terms, the change of variables:
\begin{equation*}
\begin{aligned}
	\left.\begin{array}{r@{\:=\:}l}
		t_2 & \tau_2 + \tilde{\tau}_2, \\ 
		t_1 & \tau_1 + \tilde{\tau}_1 - \tilde{\tau}_2,
	\end{array}\right\}
	\qquad \text{with} \qquad
	\left\{\begin{array}{r@{\:=\:}l}
		\tilde{\tau}_2 & t_2 - \tau_2, \\
		\tilde{\tau}_1 & \tilde{\tau}_2 + t_1 - \tau_1 = t_2 + t_1 - \tau_2 - \tau_1,
	\end{array}\right.
\end{aligned}
\end{equation*}
yields
\begin{equation}
\begin{aligned}  \label{eq:Y2-reg}
\lefteqn{Y^{\square}_2(s_1,s_2) = \int_{t_1 = - \infty}^{\infty}\int_{t_2 = - \infty}^{\infty}}\\ 
&\times
\underbrace{\int_{\tau_1 = - \infty}^{\infty} \int_{\tau_2 = - \infty}^{\infty} \!\! g_2^{\square}(\tau_1,\tau_2)u(t_2 + t_1 - \tau_2 - \tau_1)u(t_2 - \tau_2)\,\dd \tau_2\dd \tau_1}_{:=y_2^{\square}(t_1,t_2)}
\eu^{-s_1t_1}\eu^{-s_2t_2}\sigma(\ast)\dd t_2\dd t_1\, , \\
&= \int_{t_1 = -\infty}^{\infty}\int_{t_2 = - \infty}^{\infty} y_2^{\square}(t_1,t_2)\eu^{-s_1t_1}\eu^{-s_2t_2}\sigma(\ast)\dd t_2\dd t_1\, .
\end{aligned}
\end{equation}
Comparing the definition of the regular auxiliary output~$y^{\square}_2(t_1,t_2)$ from~\eqref{eq:Y2-reg} with equation~\eqref{eq:y-reg-SISO}, it follows that
\begin{equation}
\begin{aligned}
y_2(t) := y^{\square}_2(t_1 = 0,t_2 = t)\, .
\end{aligned}
\end{equation}
For the usual two-dimensional Laplace transform, the heaviside step function would normally be~$\sigma(\ast) = \sigma(t_1)\sigma(t_2)$. A substitution of~$t_1$, $t_2$ by~$\tilde{\tau}_1$,~$\tilde{\tau}_2$ however leads to
\begin{equation}
\begin{aligned}
\sigma(t_1)\sigma(t_2) = \underbrace{\sigma(\tau_1 + \tilde{\tau}_1 - \tilde{\tau}_2)}_{\lightning}\underbrace{\sigma(\tau_2 + \tilde{\tau}_2)}_{\checkmark} \not> 0, \text{ for } \tau_1,\tau_2,\tilde{\tau}_1,\tilde{\tau}_2 > 0\, .
\end{aligned}
\end{equation}
On the contrary, the following ansatz for~$\sigma(\ast)$ fulfills the requirements:
\begin{equation}
\begin{aligned}
\sigma(\ast) = \sigma(t_2 + t_1)\sigma(t_2) = \underbrace{\sigma(\tau_2 + \tau_1 + \tilde{\tau}_1 )}_{\checkmark}\underbrace{\sigma(\tau_2 + \tilde{\tau}_2)}_{\checkmark} > 0, \text{ for } \tau_1,\tau_2,\tilde{\tau}_1,\tilde{\tau}_2 > 0\, .
\end{aligned}
\end{equation} 
Therefore, we define the Laplace transform of~$y^{\square}_2(t_1,t_2)$ as:
\begin{equation}
\begin{aligned}
Y^{\square}_2(s_1,s_2) &:= \int_{t_1 = -\infty}^{\infty}\int_{t_2 = - \infty}^{\infty}y_2^{\square}(t_1,t_2)\eu^{-s_1t_1}\eu^{-s_2t_2}\sigma(t_2 + t_1)\sigma(t_2)\dd t_2\dd t_1 \\[0.2em]
&= \int_{t_2 = 0}^{\infty}\int_{t_1 = - t_2}^{\infty}y_2^{\square}(t_1,t_2)\eu^{-s_1t_1}\eu^{-s_2t_2}\dd t_1\dd t_2\, .
\end{aligned}
\end{equation}

Proceeding in a similar way, the Laplace transform~$Y^{\triangle}_k(s_1,\ldots,s_k)$ of the triangular output $y^{\triangle}_k(t_1,\ldots,t_k)$ of the~$k$-th subsystem was given by means of~$G_k^{\triangle}(s_1,\ldots,s_k)$ in equation~\eqref{eq:Y-tri-SISO}.
To come up with a formula in terms of~$G_k^{\square}(s_1,\ldots,s_k)$, we perform the substitution:
\begin{equation}
\begin{aligned}
\left.\begin{array}{r@{\:=\:}l}
\tilde{s}_1 & s_1, \\ 
\tilde{s}_{2} & s_1 + s_2,\\  
\multicolumn{2}{l}{\quad \, \, \, \vdots} \\
\tilde{s}_{k} & s_1 + \ldots + s_k,
\end{array}\right\}
\qquad \text{with} \qquad
\left\{\begin{array}{r@{\:=\:}l}
s_1 & \tilde{s}_1, \\
s_{2} & \tilde{s}_2 - \tilde{s}_1, \\
\multicolumn{2}{l}{\quad \, \, \, \vdots} \\
s_{k} & \tilde{s}_k - \tilde{s}_{k-1},
\end{array}\right.		
\end{aligned}
\end{equation} 
and let again~$\tilde{s}_1 \rightarrow s_1$,~$\ldots$,~$\tilde{s}_k \rightarrow s_k$:\\[0.2em]
\fbox{\parbox{14.7cm}{
\begin{equation}
\begin{aligned} \label{eq:Y-reg-SISO}
Y^{\square}_k(s_1,\ldots,s_k) = \underbrace{\vec c^{\mathsf T}\big(s_k\matn I - \mat A\big)^{-1}\mat N \cdots \mat N (s_1\matn I - \mat A)^{-1}\vec b}_{G_k^{\square}(s_1,\ldots,s_k)} \, U(s_1)U(s_2 - s_1)\cdots U(s_k - s_{k-1})\, .
\end{aligned}
\end{equation}
}}\\[1ex]
In order to derive the regular output~$y^{\square}_k(t_1,\ldots,t_k)$, the convolution is then applied backwards:
\begin{equation}
\begin{aligned}
\lefteqn{Y^{\square}_k(s_1,\ldots,s_k) = \underbrace{\int_{\tau_k = 0}^{\infty}\!\!\!\!\cdots\int_{\tau_1 = 0}^{\infty}\vec c^{\mathsf T}\eu^{\mat A \tau_k}\mat N\cdots \mat N \eu^{\mat A \tau_1}\vec b \eu^{-s_k \tau_k}\cdots\eu^{-s_1\tau_1}\,\dd\tau_1\cdots\dd\tau_k}_{G_k^{\square}(s_1,\ldots,s_k)}} \\
&\qquad\times
\underbrace{\int_{\tilde{\tau}_1=0}^{\infty}u(\tilde{\tau}_1)\eu^{-s_1\tilde{\tau}_1}\,\dd\tilde{\tau}_1}_{U(s_1)}
\underbrace{\int_{\tilde{\tau}_2=0}^{\infty}u(\tilde{\tau}_2)\eu^{-(s_2 - s_1)\tilde{\tau}_2}\,\dd\tilde{\tau}_2}_{U(s_2 - s_1)}\cdots \underbrace{\int_{\tilde{\tau}_k=0}^{\infty}u(\tilde{\tau}_k)\eu^{-(s_k - s_{k-1})\tilde{\tau}_k}\,\dd\tilde{\tau}_k}_{U(s_k - s_{k-1})} \\
& = \smo{\int_{\ \ \tilde{\tau}_1,\cdots,\tilde{\tau}_k = -\infty}} \cdots \int\limits^{\infty}  \ \int \cdots \smo{\int_{\!\!\!\!\!\!\!\!\!\!\!\! \tau_1,\cdots,\tau_k = -\infty}}^{\infty} \ \ \ \underbrace{\vec c^{\mathsf T} \eu^{\mat A\tau_k}\mat N\!\cdots\! \mat N \eu^{\mat A \tau_1}\vec b}_{g_k^{\square}(\tau_1,\ldots,\tau_k)} u(\tilde{\tau}_1)\cdots u(\tilde{\tau}_k)\\
&\qquad\times
\eu^{-s_1\scriptsize{\overbrace{(\tau_1 + \tilde{\tau}_1 - \tilde{\tau}_2)}^{t_1}}}\cdots
\eu^{-s_{k-1}\scriptsize{\overbrace{(\tau_{k-1} + \tilde{\tau}_{k-1} - \tilde{\tau}_k)}^{t_{k-1}}}} 
\eu^{-s_k\scriptsize{\overbrace{(\tau_k + \tilde{\tau}_k)}^{t_k}}}\sigma(\ast)\,\dd \tau_k\cdots\dd \tau_1\dd \tilde{\tau}_k\cdots\dd \tilde{\tau}_1
\end{aligned}
\end{equation}
To finally come up with well-known~$\eu^{-s_kt_k}$-terms, the change of variables:
\begin{equation*}
\begin{aligned}
\left.\begin{array}{r@{\:=\:}l}
	t_k & \tau_k + \tilde{\tau}_k, \\ 
	t_{k-1} & \tau_{k-1} + \tilde{\tau}_{k-1} - \tilde{\tau}_k, \\
	\multicolumn{2}{l}{\quad\quad \, \, \vdots} \\
	t_1 & \tau_1 + \tilde{\tau}_1 - \tilde{\tau}_2,
\end{array}\right\}
\qquad \text{with} \qquad
\left\{\begin{array}{r@{\:=\:}l}
	\tilde{\tau}_k & t_k - \tau_k, \\
	\tilde{\tau}_{k-1} & \tilde{\tau}_k + t_{k-1} - \tau_{k-1} = t_k + t_{k-1} - \tau_k - \tau_{k-1}, \\
	\multicolumn{2}{l}{\quad\quad \, \, \vdots} \\
	\tilde{\tau}_1 & \tilde{\tau}_2 + t_1 - \tau_1 = t_k + \ldots + t_1 - \tau_k - \ldots - \tau_1,
\end{array}\right.	
\end{aligned}
\end{equation*}
yields
\begin{equation}
\begin{aligned} \label{eq:y_k-reg-aux-SISO}
\lefteqn{Y^{\square}_k(s_1,\ldots,s_k)= \int_{t_1 = - \infty}^{\infty} \!\!\!\!\! \cdots\int_{t_k = - \infty}^{\infty}}\\ 
&\qquad \times
\underbrace{\int_{\tau_1 = - \infty}^{\infty} \!\!\!\!\!\!\!\! \cdots \int_{\tau_k = - \infty}^{\infty} \!\!\!\!\!\! g_k^{\square}(\tau_1,\ldots,\tau_k) u(t_k + \cdots + t_1 - \tau_k - \cdots - \tau_1)\cdots u(t_k - \tau_k)\dd \tau_k \cdots \dd \tau_1}_{:=y_k^{\square}(t_1,\ldots,t_k)}\\
&\qquad\times
\eu^{-s_1t_1}\cdots\eu^{-s_kt_k}\sigma(\ast)\,\dd t_k\cdots\dd t_1 \\
&=\int_{t_1 = - \infty}^{\infty}\!\!\!\!\!\!\!\!\!\cdots\int_{t_k = - \infty}^{\infty}y_k^{\square}(t_1,\ldots,t_k)\eu^{-s_1t_1}\cdots\eu^{-s_kt_k}\sigma(\ast)\dd t_k\cdots\dd t_1\, .
\end{aligned}
\end{equation}
The following ansatz for the heaviside step function~$\sigma(\ast)$ of the~$k$-th Laplace transform fulfills the requirements:
\begin{equation}
\begin{aligned}
\sigma(\ast) &= \sigma(t_k)\sigma(t_k + t_{k-1})\cdots\sigma(t_k + \cdots + t_1) \\
&= 
\underbrace{\sigma(\tau_k + \tilde{\tau}_k)}_{\checkmark}
\underbrace{\sigma(\tau_k + \tau_{k-1} + \tilde{\tau}_k + \tilde{\tau}_{k-1})}_{\checkmark}\cdots
\underbrace{\sigma(\tau_k + \cdots + \tau_1 + \tilde{\tau}_k + \cdots +\tilde{\tau}_1)}_{\checkmark}
> 0\, , \\ &\qquad\text{for } \tau_1,\ldots,\tau_k,\tilde{\tau}_1,\ldots,\tilde{\tau}_k > 0
\end{aligned}
\end{equation}
Thus, we define the Laplace transform of~$y^{\square}_k(t_1,\ldots,t_k)$ as:\\
\fbox{\parbox{14.7cm}{
\begin{equation}
\begin{aligned}
\lefteqn{Y^{\square}_k(s_1,\ldots,s_k)} \\ 
&:= \int_{t_1 = -\infty}^{\infty}\!\!\!\!\!\!\!\cdots \int_{t_k = - \infty}^{\infty} y_k^{\square}(t_1,\ldots,t_k)\eu^{-s_1t_1}\cdots\eu^{-s_kt_k}\sigma(t_k)\cdots\sigma(t_k + \ldots + t_1)\dd t_k\cdots\dd t_1 \\
&= \int_{t_k = 0}^{\infty}\int_{t_{k-1} = - t_k}^{\infty}\!\!\!\!\!\!\cdots \int_{t_1 = - t_k - \ldots - t_2}^{\infty}y_k^{\square}(t_1,\ldots,t_k)\eu^{-s_1t_1}\cdots\eu^{-s_kt_k}\dd t_k\cdots\dd t_1\, .
\end{aligned}
\end{equation}
}}\\[0.5ex]
\begin{remark}[Response computation from the regular frequency domain representation]
	The Laplace transform~$Y^{\square}_k(s_1,\ldots,s_k)$ of the regular auxiliary output is given in terms of~$G^{\square}_k(s_1,\ldots,s_k)$ in equation~\eqref{eq:Y-reg-SISO}. Please bear in mind the shifted input terms~$U(s_1)U(s_2 - s_1)\ldots U(s_k - s_{k-1})$. The regular auxiliary output~$y^{\square}_k(t_1,\ldots,t_k)$ can then be computed via a $k$-dimensional (regular) inverse Laplace transform. As the limits of the regular Laplace transform  are induced by the heaviside step function~$\sigma(\ast)$, we think, that the same heaviside step function~$\sigma(\ast)$ will result when applying a $k$-dimensional (regular) inverse Laplace transform:~$\mathcal{L}_k^{-1}\{Y^\square_k(s_1,\ldots,s_k)\} = y^\square_k(t_1,\ldots,t_k)\sigma(\ast)$.		
\end{remark}

\subsubsection*{MIMO case}
Similar to the SISO case, the regular auxiliary output~$\vec y^{\square}_k(t_1, \ldots, t_k)$ is:
\begin{equation}
\begin{aligned} \label{eq:y(t1,...,tk)-reg-MIMO}
\lefteqn{\vec y^{\square}_k(t_1, \ldots, t_k) = \sum_{j_1 = 1}^{m} \! \cdots \! \sum_{j_k = 1}^{m}} \\ 
&\qquad \times \int \cdots \smo{\int_{\!\!\!\!\!\!\!\!\!\!\!\!\!\!\!\!\!\!\!\!\!\! \tau_1,\cdots,\tau_k = -\infty}}^{\infty} \, \vec g^{(j_1,\ldots,j_k)}_{k, \square}(\tau_1,\ldots,\tau_k)u_{j_k}(t_k + \ldots + t_1 - \tau_k - \ldots - \tau_1)\cdots u_{j_1}(t_k-\tau_k)\,\dd \tau_k \cdots \dd \tau_1\, .
\end{aligned}
\end{equation}
Due to the linearity property of the Laplace transform,~$\vec Y^{\square}_k(s_1, \ldots, s_k) \in \mathbb{C}^p$ is given by:\\[1ex]
\fbox{\parbox{14.7cm}{
\begin{equation}
	\begin{aligned} \label{eq:Y_k-reg-MIMO}
		\lefteqn{\vec Y^{\square}_k(s_1, \ldots, s_k) =  \sum_{j_1 = 1}^{m}\cdots \sum_{j_k = 1}^{m} \vec Y^{(j_1,\ldots,j_k)}_{k,\square}(s_1, \ldots, s_k)} \\
		&= \sum_{j_1 = 1}^{m} \cdots \sum_{j_k = 1}^{m} \underbrace{\mat C (s_k \matn I - \mat A)^{-1}\mat N_{j_k} \cdots \mat N_{j_2} (s_1 \matn I - \mat A)^{-1}\vec b_{j_1}}_{\vec G_{k,\square}^{(j_1,\ldots,j_k)}(s_1,\ldots,s_k)} U_{j_k}(s_k - s_{k-1}) \cdots U_{j_2}(s_2 - s_1) U_{j_1}(s_1)\, .
	\end{aligned}
\end{equation}
}}\\[0.5em]

\section{Impulse response of bilinear systems} \label{sec:imp-resp}
After having introduced the triangular and regular auxiliary outputs~$y^\triangle_k(t_1,\ldots,t_k)$ and $y^\square_k(t_1,\ldots,t_k)$ in the previous section, at this point one could think that the response of the $k$-th subsystem to an impulse input $u(t) = \delta(t)$ could be calculated from $g_k(t) = g^\triangle_k(t_1 = t,\ldots,t_k = t)$ or $g_k(t) = g^\square_k(t_1 = 0,\ldots,t_{k-1} = 0,t_k = t)$. This, however, is wrong.

Remember that the $k$-dimensional Laplace transform is not unique in terms of~$(k-1)$-di\-men\-sio\-nal discontinuities. We will show that the one-dimensional time line~$t_1,\ldots,t_k = t$ in the triangular case and~$t_1,\ldots,t_{k-1} = 0$,~$t_k = t$ in the regular case  is a discontinuity for the impulse response, and thus, the impulse response cannot be computed via the $k$-dimensional inverse Laplace transform.\footnote{Remember that we did not yet define the triangular kernel~$g^\triangle_k(t_1,\ldots,t_k)$ and the regular kernel~$g^\square_k(t_1,\ldots,t_k)$ on the~$(k-1)$-dimensional surface, especially not for~$t_1,\ldots,t_k = t$ in the triangular case and~$t_1,\ldots,t_{k-1} = 0$,~$t_k = t$ in the regular case.}\\
Hence, in the following we will derive the impulse response of bilinear systems directly in the time domain, where a factor 1/k! -- that has not appeared until now -- will arise in the solution.

\subsection{Derivation of the impulse response} \label{subsec:derivation-imp-direct}
The MIMO bilinear system from equation~\eqref{eq:bilinear-FOM} is for~$\vec u(t) = \vec \mu \delta(t)$ with the Dirac delta function~$\delta(t)$ and the scaling of inputs~$\vec \mu \in \mathbb{R}^m$:
\begin{equation}
\begin{aligned}
\dot{\vec x}(t) &= \underbrace{\left(\mat A + \sum\limits_{j = 1}^{m}\mat N_j \mu_j \delta(t)\right)}_{\mat A(t)}
\vec x(t) + \mat B \vec \mu \delta(t), \quad &&\vec x(0) = \vec x_0\,,\\
\quad \vec y(t)&=\mat C \vec x(t)\, .
\end{aligned}
\end{equation}
The solution of a system of differential equations with distributional dynamic matrix~$\mat A(t)$ is, to the best of the authors' knowledge, not defined~\cite{walter1994einfuehrung}. Therefore, we will derive the solution of the bilinear system for a nascent delta function~$\delta_\varepsilon(t)$ with constant amplitude~$1/\varepsilon$ and duration of action~$\varepsilon$ in the limit of~$\varepsilon \to 0$:
\begin{equation}
\begin{aligned}
\delta(t) := \lim\limits_{\varepsilon\rightarrow 0} \delta_\varepsilon(t)\, , \qquad \text{with} \qquad 
\delta_\varepsilon(t) = \left\{\ \begin{array}{cl}
\D{\tfrac{1}{\varepsilon}}, & 0 \leq t \leq \varepsilon \\
0, & \mathrm{else.}
\end{array} \right.
\end{aligned}
\end{equation}
Then, the system excited by~$\vec u(t) = \vec \mu \delta_{\varepsilon}(t)$ is (with~$\hat{\mat N} = \sum_{j=1}^{m}\mat N_j\mu_j$ and~$\hat{\vec b} = \mat B \vec \mu = \sum_{j=1}^{m} \vec b_j \mu_j$):
\begin{equation}
\begin{aligned} \label{eq:nascent-system}
\dot{\vec x}_\varepsilon(t) &= \left(\mat A + \tfrac{1}{\varepsilon}\hat{\mat N}\right)\vec x_\varepsilon(t) + \tfrac{1}{\varepsilon}\hat{\vec b}, \quad &&0 \leq t \leq \varepsilon, \quad \vec x_\varepsilon(0) = \vec x_0\,, \\
\dot{\vec x}_\varepsilon(t) &= \mat A\vec x_\varepsilon(t), &&t > \varepsilon\, .
\end{aligned}
\end{equation}
The solution~$\vec x_\varepsilon(t)$ in interval~$0 \leq t \leq \varepsilon$ is (with~$\widehat{\mat A} := \mat A + \frac{1}{\varepsilon}\hat{\mat N}, \:\: \det \widehat{\mat A} \not= 0$):
\begin{equation}
\begin{aligned}
\label{eq:im_lo_1}
\vec x_\varepsilon(t) &= \int_{\tau =0}^{t}\eu^{\widehat{\mat A}(t-\tau)}\tfrac{1}{\varepsilon}\hat{\vec b}\, \dd \tau \ + \ \eu^{\widehat{\mat A}t}\vec x_0 = \tfrac{1}{\varepsilon}\eu^{\widehat{\mat A}t} \int_{\tau = 0}^{t}\eu^{-\widehat{\mat A}\tau}\, \dd \tau\: \hat{\vec b} \ + \ \eu^{\widehat{\mat A}t}\vec x_0 \\[0.2em]
&= \tfrac{1}{\varepsilon}\eu^{\widehat{\mat A}t} \left.\left(-\eu^{-\widehat{\mat A}\tau}\hat{\mat A}^{-1}\right)\right|_{\tau = 0}^{t}\: \hat{\vec b} \ + \ \eu^{\widehat{\mat A}t}\vec x_0 = \tfrac{1}{\varepsilon}( \eu^{\widehat{\mat A}t} - \matn I)\widehat{\mat A}^{-1}\hat{\vec b} \ + \ \eu^{\widehat{\mat A}t}\vec x_0 \, .
\end{aligned}
\end{equation} 
The transition condition for the solution~$\vec x_\varepsilon(\varepsilon)$ at the time point~$t = \varepsilon$ is:
\begin{equation}
\begin{aligned}
\vec x_\varepsilon(\varepsilon) &= \tfrac{1}{\varepsilon}(\eu^{(\mat A + \frac{1}{\varepsilon}\hat{\mat N})\varepsilon} - \matn I) \left(\mat A + \tfrac{1}{\varepsilon}\hat{\mat N}\right)^{-1}\hat{\vec b} \ + \ \eu^{(\mat A + \frac{1}{\varepsilon}\hat{\mat N})\varepsilon}\vec x_0 \\
&= (\eu^{\varepsilon\mat A + \hat{\mat N}} - \matn I)(\varepsilon \mat A + \hat{\mat N})^{-1}\hat{\vec b} \ + \ \eu^{\varepsilon \mat A + \hat{\mat N}}\vec x_0\, .
\end{aligned}
\end{equation}
With the exponential series expansion
\begin{equation}
\begin{aligned}
(\eu^{\varepsilon\mat A + \hat{\mat N}} - \matn I)(\varepsilon \mat A + \hat{\mat N})^{-1} &=
\Bigg(\sum_{k = 0}^{\infty}\frac{(\varepsilon\mat A + \hat{\mat N})^k}{k!} - \matn I\Bigg)(\varepsilon \mat A + \hat{\mat N})^{-1} \\
&= \sum_{k = 1}^{\infty}\frac{(\varepsilon\mat A + \hat{\mat N})^k}{k!}(\varepsilon \mat A + \hat{\mat N})^{-1} = \sum_{k = 1}^{\infty}\frac{(\varepsilon\mat A + \hat{\mat N})^{k-1}}{k!}\, , 
\end{aligned}
\end{equation}
the transition condition becomes:
\begin{equation}
\begin{aligned}
\vec x_\varepsilon(\varepsilon)  &= \sum_{k = 1}^{\infty}\frac{(\varepsilon\mat A + \hat{\mat N})^{k-1}}{k!}\hat{\vec b} \ + \ \eu^{\varepsilon \mat A + \hat{\mat N}}\vec x_0\, .
\end{aligned}
\end{equation}
The solution~$\vec x_\varepsilon(t)$ of the system~\eqref{eq:nascent-system} for $t>\varepsilon$ is thus:
\begin{equation}
\begin{aligned}
\label{eq:im_lo_2}
\vec x_\varepsilon(t) &= \eu^{\mat A(t - \varepsilon)}\vec x(\varepsilon) \\
&= \eu^{\mat At}\eu^{-\varepsilon\mat A}\sum_{k = 1}^{\infty}\frac{(\varepsilon\mat A + \hat{\mat N})^{k-1}}{k!}\hat{\vec b} \ + \ \eu^{\mat At}\eu^{-\varepsilon\mat A}\eu^{\varepsilon \mat A + \hat{\mat N}}\vec x_0\, .
\end{aligned}
\end{equation}
Note that the simplification~$\eu^{-\varepsilon\mat A}\eu^{\varepsilon \mat A + \hat{\mat N}} = \eu^{\hat{\mat N}}$ does generally not hold and is only admissible for commuting matrices (i.e.~$\mat A \hat{\mat N} = \hat{\mat N}\mat A$), since then~$\eu^{\mat A + \hat{\mat N}} = \eu^{\mat A} \eu^{\hat{\mat N}}$.

The interesting part of the solution is the interval~$t > \varepsilon$, i.e. the response \emph{after} the action of the nascent delta impulse. On the contrary, the interval~$0 \leq t \leq \varepsilon$, which describes how the system moves from the initial state~$\vec x_0$ to the transition state~$\vec x_\varepsilon(t = \varepsilon)$, is not of interest. A similar thing happens in billiard: the force that acts on billiard balls during the impact are generally unknown, difficult to measure or rather uninteresting. Of more relevance is, however, the velocity of the balls after the impact.

If the duration of action~$\varepsilon$ of the nascent delta function~$\delta_\varepsilon(t)$ approaches zero -- and at the same time the amplitude~$1/\varepsilon \to \infty$ -- then the solution~$\vec x(t)$ for~$t > \varepsilon$ becomes:
\begin{equation}
	\begin{aligned} \label{eq:imp-x(t)}
		\vec x(t) &:= \lim\limits_{\varepsilon\rightarrow 0} \vec x_\varepsilon(t) \\
		&= \eu^{\mat A t}\bigg(\lim\limits_{\varepsilon \rightarrow 0}\underbrace{\eu^{-\varepsilon\mat A}}_{\rightarrow \matn I} \sum_{k = 1}^{\infty} \underbrace{\frac{(\varepsilon\mat A + \hat{\mat N})^{k-1}}{k!}}_{\to \frac{\hat{\mat N}^{k-1}}{k!}}\hat{\vec b} \ + \ \underbrace{\eu^{-\varepsilon\mat A}}_{\to \matn I}\underbrace{\eu^{\varepsilon \mat A + \hat{\mat N}}}_{\rightarrow \eu^{\hat{\mat N}}}\vec x_0\bigg) \\
		&= \eu^{\mat A t}\sum_{k = 1}^{\infty}\frac{\hat{\mat N}^{k-1}}{k!}\hat{\vec b} \ + \ \eu^{\mat A t}\eu^{\hat{\mat N}}\vec x_0\, .
	\end{aligned}
\end{equation}

The output of the system~$\vec y_\varepsilon(t)$ for~$\varepsilon \to 0$ represents the \emph{impulse response}~$\vec g(t;\vec \mu)$ of the bilinear system and is given by:
\\[1ex]
\fbox{\parbox{14.7cm}{		
\begin{equation}
\begin{aligned} \label{eq:impulse-response}
\vec g(t; \vec \mu) &:= \lim\limits_{\varepsilon \rightarrow 0}\vec y(t)\ \ \text{subject to} \ \ \vec u(t) = \vec \mu \delta_\varepsilon(t)\\
&= \underbrace{\mat C\eu^{\mat A t}\sum_{k = 1}^{\infty}\frac{\hat{\mat N}^{k-1}}{k!}\hat{\vec b}}_{\vec g_{\mat B}(t)} \ \ + \ \  
\underbrace{\mat C \eu^{\mat A t}\eu^{\hat{\mat N}}\vec x_0}_{\vec g_{\mathrm{IC}}(t)}, \quad \hat{\mat N} = \sum_{j = 1}^{m}\mat N_j \mu_j, \quad \hat{\vec b} = \mat B \vec \mu = \sum_{j=1}^{m} \vec b_j \mu_j\, ,
\end{aligned}
\end{equation}}}
\\[1ex]
where~$\vec g_{\mathrm{IC}}(t)$ and~$\vec g_{\mat B}(t)$ correspond to the initial condition and the excitation term, respectively.\footnote{Note that in \cite{gebhart2017impulse} the initial condition term is denoted by the abbreviation AW (german: \emph{Anfangswert}).}
As one would expect, setting~$\mat N_1,\ldots,\mat N_m = \matn 0$ implies $\hat{\mat N} = \matn 0$ and therefore
\begin{equation}
\begin{aligned}
\sum_{k = 1}^{\infty}\frac{\matn 0^{k-1}}{k!} = \matn I \ \ \mathrm{and} \ \ \eu^{\matn 0} = \matn I\, ,
\end{aligned}
\end{equation}
so that the above expression boils down to the well-known impulse response of linear systems:
\begin{equation}
\begin{aligned}
\vec g_{\mathrm{lin}}(t;\vec \mu) = \mat C\eu^{\mat A t}\mat B \vec \mu \ + \ 
\mat C \eu^{\mat A t}\vec x_0\, .
\end{aligned}
\end{equation}
\vspace{-2em}
\begin{remark}[Single-variable transfer function]
	The impulse response~$\vec g(t; \vec \mu)$ could be transformed into the frequency domain to obtain the single-variable \enquote{transfer function}~$\vec G(s;\vec \mu)$. This \enquote{transfer function}, however, does not represent the input-output behavior of the bilinear system. Thus, it cannot be used as a starting point for model order reduction with e.g. moment matching. 
\end{remark}

\subsection{Plausibility check of the result}
A \emph{scalar} bilinear initial value problem with~$a(t):=a + nu(t)$ can be solved for~$a(t),u(t) \in \mathcal{C}^0([t_0,t_\mathrm{max}],\mathbb{R})$ by the separation of variables for the homogeneous solution and by the method of variation of parameters for the particular solution.\footnote{Note, however, that the solution~$\vec x(t)$ of a system of bilinear differential equations cannot be analytically given.}
For an impulse excitation~$u(t) = \delta(t)$, where the impulse affects the input and the dynamic simultaneously, the assumptions~$a(t),u(t) \in \mathcal{C}^0$ are not met. Solutions for distributional excitations~$u(t)\not\in \mathcal{C}^0$ are given in the theory of distributions~\cite{walter1994einfuehrung}, solutions for distributional dynamic~$a(t)\not\in \mathcal{C}^0$, however, not. Therefore, our aim is now to estimate the interaction between excitation~$u(t)$ and dynamic~$a(t)$ by computing the response of the system, when the impulse affects the input and the dynamic consecutively, instead of simultaneously.

\subsubsection*{Impulse on the input and then on the dynamic}
If an impulse first acts on the input term~$bu(t)$ (while the bilinear term~$nx(t)u(t)$ is zero), then:
\vspace{-0.2em}
\begin{equation}
\begin{aligned}
	\dot{x}(t) &= ax(t) + b\delta(t), \quad &&x(0) = x_0\, .
\end{aligned}
\end{equation}
\vspace{-0.2em} 
After the action of the first impulse at~$t = 0^{+}$, the state is~$x(0^+) = x_0 + b$, like in the linear case. 
If another impulse acts afterwards on the bilinear term~$nx(t)u(t)$ (while~$bu(t)$ is zero), then:
\vspace{-0.2em}
\begin{equation}
\begin{aligned}
	\dot{x}(t) &= \big(a + n \delta(t)\big)x(t), \quad &&x(0^+) = b + x_0\, .
\end{aligned}
\end{equation}
\vspace{-0.2em}
The state after the action of the second impulse at~$t = 0^{++}$ is~\mbox{$x(0^{++}) = \eu^n(x_0 + b)$}. After both impulses, the system becomes autonomous:
\vspace{-0.2em}
\begin{equation}
\begin{aligned}
	\dot{x}(t) &= ax(t), \quad &&x(0^{++}) = \eu^n(b + x_0)\, .
\end{aligned}
\end{equation}
\vspace{-0.2em}
Thus, once both impulses acted on the system the state is given by~$x(t) = \eu^{at}\eu^n(b + x_0)$ for~$t > 0^{++}$.

\subsubsection*{Impulse on the dynamic and then on the input}
On the contrary, if an impulse first acts on the bilinear term~$nx(t)u(t)$ (while the input term~$bu(t)$ is zero), then:
\vspace{-0.2em}
\begin{equation}
\begin{aligned}
	\dot{x}(t) &= (a + n\delta(t))x(t), \quad &&x(0) = x_0\, .
\end{aligned}
\end{equation} 
\vspace{-0.2em}
After the action of the first impulse at ~$t = 0^{+}$, the state is~$x(0^{+}) = \eu^n x_0$. If another impulse acts afterwards on the input term~$bu(t)$ (while~$nx(t)u(t)$ is zero), then:
\vspace{-0.2em}
\begin{equation}
\begin{aligned}
	\dot{x}(t) &= ax(t) + b\delta(t), \quad &&x(0^+) = \eu^n x_0\, .
\end{aligned}
\end{equation}
\vspace{-0.2em}
The state after the action of the second impulse at~$t = 0^{++}$ is~$x(0^{++}) = b + \eu^n x_0$. After both impulses, the system becomes autonomous:
\vspace{-0.2em}
\begin{equation}
\begin{aligned}
	\dot{x}(t) &= ax(t), \quad &&x(0^{++}) = b + \eu^n x_0\, .
\end{aligned}
\end{equation}
\vspace{-0.2em}
Thus, once both impulses acted on the system the state is given by~$x(t) = \eu^{at}(b + \eu^n x_0)$ for~$t > 0^{++}$.\\

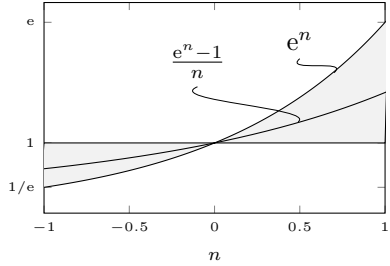
\begin{figure}[tp!]
\pgfplotsset{footnotesize,samples=5}
\setlength\figurewidth{0.3\textwidth}
\setlength\figureheight{0.61803\figurewidth}
\begin{center}
	\centering
%
%
\begin{tikzpicture}

\begin{axis}[%
width=\figurewidth,
height=\figureheight,
at={(0.1\figurewidth,0.1\figureheight)},
xlabel style={font=\color{white!15!black}, font=\scriptsize},
ylabel style={font=\color{white!15!black}, font=\footnotesize},
x tick label style={font=\tiny},
y tick label style={font=\tiny},
scale only axis,
xlabel={$n$},
ytick={0.36788,1,2.71828}, 
yticklabels={$1/\eu$, 1, $\eu$},
xtick={-1,-0.5,0,0.5,1}, 
xmin=-1,
xmax=1,
ymin=0,
ymax=3,
axis background/.style={fill=white}
]

\addplot[area legend, draw=black, fill=black!5!white, forget plot]
table[row sep=crcr] {%
	x	y\\
	-1	1\\
	1	1\\
	1.01995075129403	2.77305819107383\\
	0.999949746218403	2.71814522795009\\
	0.979948741142771	2.66431966844765\\
	0.959947736067139	2.61155997946086\\
	0.939946730991507	2.55984505428938\\
	0.919945725915875	2.50915420419429\\
	0.899944720840243	2.45946715012157\\
	0.879943715764611	2.41076401458933\\
	0.859942710688979	2.36302531373577\\
	0.839941705613347	2.31623194952462\\
	0.819940700537716	2.27036520210487\\
	0.799939695462084	2.22540672232182\\
	0.779938690386452	2.18133852437649\\
	0.75993768531082	2.13814297863026\\
	0.739936680235188	2.09580280455215\\
	0.719935675159556	2.05430106380562\\
	0.699934670083924	2.01362115347235\\
	0.679933665008292	1.97374679941021\\
	0.65993265993266	1.93466204974268\\
	0.639931654857028	1.8963512684773\\
	0.619930649781396	1.85879912925038\\
	0.599929644705764	1.82199060919569\\
	0.579928639630132	1.78591098293447\\
	0.5599276345545	1.75054581668451\\
	0.539926629478868	1.71588096248583\\
	0.519925624403236	1.68190255254081\\
	0.499924619327604	1.64859699366625\\
	0.479923614251973	1.61595096185542\\
	0.459922609176341	1.5839513969477\\
	0.439921604100709	1.55258549740388\\
	0.419920599025077	1.52184071518481\\
	0.399919593949445	1.49170475073151\\
	0.379918588873813	1.4621655480447\\
	0.359917583798181	1.4332112898618\\
	0.339916578722549	1.40483039292929\\
	0.319915573646917	1.3770115033689\\
	0.299914568571285	1.34974349213536\\
	0.279913563495653	1.32301545056425\\
	0.259912558420021	1.29681668600791\\
	0.239911553344389	1.27113671755783\\
	0.219910548268757	1.24596527185173\\
	0.199909543193125	1.22129227896363\\
	0.179908538117493	1.19710786837538\\
	0.159907533041861	1.17340236502791\\
	0.139906527966229	1.15016628545068\\
	0.119905522890598	1.1273903339678\\
	0.0999045178149656	1.10506539897923\\
	0.0799035127393336	1.08318254931571\\
	0.0599025076637017	1.06173303066575\\
	0.0399015025880698	1.04070826207349\\
	0.0199004975124378	1.02009983250582\\
	0	1\\
	-0.0199004975124378	0.980296210365564\\
	-0.0399015025880696	0.960884079086308\\
	-0.0599025076637016	0.941856352884637\\
	-0.0799035127393335	0.923205419651324\\
	-0.0999045178149656	0.904923818014498\\
	-0.119905522890598	0.887004234354705\\
	-0.139906527966229	0.869439499879063\\
	-0.159907533041861	0.85222258775336\\
	-0.179908538117493	0.835346610290949\\
	-0.199909543193125	0.818804816197302\\
	-0.219910548268757	0.802590587869129\\
	-0.239911553344389	0.786697438746986\\
	-0.259912558420021	0.771119010720302\\
	-0.279913563495653	0.755849071583792\\
	-0.299914568571285	0.740881512544247\\
	-0.319915573646917	0.726210345776684\\
	-0.339916578722549	0.711829702028899\\
	-0.359917583798181	0.697733828273449\\
	-0.379918588873813	0.683917085406136\\
	-0.399919593949445	0.670373945990061\\
	-0.419920599025077	0.657098992044356\\
	-0.439921604100709	0.644086912876698\\
	-0.45992260917634	0.631332502958749\\
	-0.479923614251973	0.61883065984367\\
	-0.499924619327604	0.606576382124863\\
	-0.519925624403236	0.594564767435142\\
	-0.539926629478868	0.582791010485528\\
	-0.5599276345545	0.571250401142872\\
	-0.579928639630132	0.559938322545548\\
	-0.599929644705764	0.548850249256468\\
	-0.619930649781396	0.537981745452658\\
	-0.639931654857028	0.527328463150693\\
	-0.65993265993266	0.516886140467274\\
	-0.679933665008292	0.506650599914243\\
	-0.699934670083924	0.496617746727367\\
	-0.719935675159556	0.486783567228208\\
	-0.739936680235188	0.477144127218443\\
	-0.75993768531082	0.46769557040597\\
	-0.779938690386452	0.458434116862187\\
	-0.799939695462083	0.449356061509815\\
	-0.819940700537715	0.440457772640673\\
	-0.839941705613347	0.431735690462796\\
	-0.859942710688979	0.423186325676331\\
	-0.879943715764611	0.414806258077629\\
	-0.899944720840243	0.406592135190978\\
	-0.919945725915875	0.39854067092744\\
	-0.939946730991507	0.39064864427023\\
	-0.959947736067139	0.382912897986146\\
	-0.979948741142771	0.3753303373625\\
	-0.999949746218403	0.36789792896907\\
	-1.01995075129403	0.360612699444566\\
}--cycle;

\draw[ultra thin] (0.49,1.3) .. controls (0.6,1.5) and (-0.3,1.6) .. (-0.1,1.8) node[anchor=south] {$\frac{\eu^n - 1}{n}$};

\draw[ultra thin] (0.7,2.0) .. controls (0.8,2.2) and (0.4,2.0) .. (0.5,2.2) node[anchor=south] {$\eu^n$};

\addplot [color=black, forget plot]
  table[row sep=crcr]{%
-2	0.432332358381694\\
-1.83999195939494	0.457165746080954\\
-1.69998492386552	0.480777041190856\\
-1.5599778883361	0.50632722950462\\
-1.41997085280667	0.53401021310847\\
-1.29996482235288	0.55959869625814\\
-1.17995879189909	0.587061686832146\\
-1.0599527614453	0.616563159211469\\
-0.939946730991507	0.648282860760623\\
-0.839941705613347	0.676552081816491\\
-0.739936680235188	0.706622453985345\\
-0.639931654857028	0.738628153900139\\
-0.539926629478868	0.77271422955589\\
-0.439921604100709	0.809037528063353\\
-0.339916578722549	0.847767705400198\\
-0.259912558420021	0.880607657710115\\
-0.179908538117493	0.915206089893969\\
-0.0999045178149656	0.951670495638581\\
-0.0199004975124377	0.990115429130408\\
0.0599025076637019	1.0305583701492\\
0.139906527966229	1.0733329433129\\
0.219910548268757	1.11847873504973\\
0.299914568571285	1.16614372486621\\
0.379918588873813	1.21648574610338\\
0.459922609176341	1.26967316956538\\
0.519925624403236	1.31153865194371\\
0.579928639630132	1.35518567152626\\
0.639931654857028	1.40069843658157\\
0.699934670083924	1.44816537427817\\
0.75993768531082	1.49767935007033\\
0.819940700537716	1.54933789879166\\
0.879943715764612	1.60324346809326\\
0.939946730991507	1.65950367489865\\
0.999949746218403	1.71823157558443\\
1.0599527614453	1.77954595063511\\
1.11995577667219	1.84357160456181\\
1.17995879189909	1.91043968191815\\
1.23996180712599	1.98028800029306\\
1.29996482235288	2.05326140120862\\
1.35996783757978	2.12951211990278\\
1.41997085280667	2.20920017503102\\
1.47997386803357	2.29249377937844\\
1.53997688326047	2.37956977273426\\
1.59997989848736	2.47061407814502\\
1.65998291371426	2.56582218282993\\
1.71998592894115	2.66539964511367\\
1.77998894416805	2.76956262880702\\
1.83999195939494	2.87853846654553\\
1.89999497462184	2.99256625368046\\
1.95999798984874	3.11189747440528\\
2	3.19452804946533\\
};
\end{axis}
\end{tikzpicture}%
\end{center} \vspace{-1ex}
\caption{Scaling of the excitation term of the impulse response. -- The scaling of the excitation term of the impulse response~$(\eu^n - 1)/n$ lies for all~$n \in \mathbb{R}$ between the limiting cases~$1$ and~$\eu^n$. These cases result, when the impulse~$u(t) = \delta(t)$ acts on the input term~$bu(t)$ and on the bilinear term~$nx(t)u(t)$ consecutively, instead of simultaneously.}
\label{fig:zusatz1}
\vspace{-1em}
\end{figure}

If both impulses act on the system simultaneously and mutually influence each other, then the solution is likely to lie between the two described cases~$x(t) = \eu^{at}\eu^n(b + x_0)$ and~$x(t) = \eu^{at}(b + \eu^n x_0)$. That means:
\begin{equation}
\begin{aligned}
	x(t) = \eu^{at}(f(n)b + \eu^n x_0), \quad \min (1,\eu^n)\leq f(n) \leq \max (1, \eu^n)\, .
\end{aligned}
\end{equation}
The function~$f(n) = (\eu^n - 1)/n = \sum_{k = 1}^{\infty}n^{k-1}/k!$ fulfills the above condition for all~$n \in \mathbb{R}$, as it can be seen in Figure~\ref{fig:zusatz1}. This calculation makes the result obtained in equation~\eqref{eq:imp-x(t)} plausible:
\begin{equation}
x(t) = \eu^{at}\sum_{k = 1}^{\infty} \frac{n^{k-1}}{k!} b \ + \ \eu^{at} \eu^n x_0\, .
\end{equation}

\subsection{Derivation of the impulse response of the~$k$-th subsystem} \label{subsec:derivation-imp-k-subs}
In section~\ref{subsec:derivation-imp-direct} we have derived the impulse response~$\vec g(t; \vec \mu)$ of bilinear systems by applying a nascent delta function~$\delta_{\varepsilon}(t)$ and~$\varepsilon \to 0$ directly to~\eqref{eq:bilinear-FOM}. Now we want to show an alternative derivation of the impulse res\-pon\-se~$\vec g(t; \vec \mu)$ by making use of the Volterra theory and the interpretation of a bilinear system as a sequence of interconnected subsystems (cf.~\eqref{eq:bilinear-subsystems}). Towards this aim, the impulse response of the~$k$-th subsystem~$\vec g_k(t; \vec \mu)$ will be first derived and then~$\vec g(t; \vec \mu) = \sum_{k=1}^{\infty} \vec g_k(t; \vec \mu)$. Similar as before, the impulse response~$\vec g_k(t; \vec \mu)$ is obtained by computing the solution~$\vec y_{\varepsilon,k}(t)$ of the~$k$-th subsystem for a nascent delta function~$\delta_{\varepsilon}(t)$ with constant amplitude~$1/\varepsilon$ and duration of action~$\varepsilon$ in the limit~$\varepsilon \to 0$. For the sake of brevity, we will focus only on the scalar case and consider solely up to the second subsystem. The interested reader is referred to~\cite[\S3.5]{gebhart2017impulse} for a general and extensive derivation.

At this point, it should be stressed that the factor~$1/k!$ only arises in the solution, if the limit~$\varepsilon \to 0$ is built \emph{once}. If the limit~$\varepsilon \to 0$ is built after each subsystem instead, then the factor~$1/k!$ does not pop up. This fact will be extensively discussed in the following.


\subsubsection*{Successive action of Dirac impulses}
First, we discuss the case where the limit~$\varepsilon \to 0$ is built after each subsystem. This corresponds to the case where nascent delta impulses act \emph{successively} on each subsystem.

The solution~$x_1(t)$ of the first subsystem when excited by~$u(t) = \delta_{\varepsilon_1}(t)$ in the limit $\varepsilon_1 \rightarrow 0$ (using L'Hospital's (L.H.) rule) is -- for initial condition $x_0 = 0$ and~$t > \varepsilon_1$ -- given by:
\begin{equation}
\begin{aligned}
	x_1(t)&=  \lim\limits_{\varepsilon_1 \rightarrow 0} \int_{\tau = 0}^{t} \eu^{a(t - \tau)}b \underbrace{\varepsilon_1^{-1}\sigma(\tau)\sigma(\varepsilon_1 - \tau)}_{\delta_{\varepsilon_1}(t)}\,\dd \tau\, =\eu^{at}b\lim\limits_{\varepsilon_1 \rightarrow 0}\varepsilon_1^{-1}\int_{\tau = 0}^{\varepsilon}\eu^{-a\tau}\,\dd \tau \\
	&= \eu^{at}b \lim\limits_{\varepsilon_1 \rightarrow 0}\frac{1 - \eu^{-a\varepsilon_1}}{a\varepsilon_1} \overset{\mathrm{L.H.}}{=}\eu^{at}b\lim\limits_{\varepsilon_1 \rightarrow 0}\frac{a\eu^{-a\varepsilon_1}}{a} \\
	&=\eu^{at}b\, .
\end{aligned}
\end{equation}
The solution~$x_2(t)$ of the second subsystem when excited by~$u(t) = \delta_{\varepsilon_2}(t)$ in the limit~$\varepsilon_2 \rightarrow 0$ would be for~$t > \varepsilon_2$:
\begin{equation}
\begin{aligned} \label{eq:x2-suc-limit}
	x_2(t) &= \lim\limits_{\varepsilon_2 \rightarrow 0}\int_{\tau_1 = 0}^{\varepsilon_2}\eu^{a(t-\tau_1)}n\varepsilon_2^{-1} \underbrace{\lim\limits_{\varepsilon_1 \rightarrow 0}\int_{\tau_2 = 0}^{\varepsilon_1}\eu^{a(\tau_1 - \tau_2)}b\varepsilon_1^{-1}\,\dd \tau_2}_{x_1(\tau_1)} \,\dd \tau_1 \\
	&=  \lim\limits_{\varepsilon_2 \rightarrow 0}\int_{\tau_1 = 0}^{\varepsilon_2}\eu^{a(t-\tau_1)}n\varepsilon_2^{-1}\eu^{a\tau_1}b\,\dd \tau_1 = \eu^{at}nb \lim\limits_{\varepsilon_2 \rightarrow 0} \varepsilon_2^{-1} \int_{\tau_1 = 0}^{\varepsilon_2} \dd \tau_1 \\[0.2em]
	&= \eu^{at}nb\, .
\end{aligned}
\end{equation}
In the calculation above, the limit~$\varepsilon_1 \rightarrow 0$ is initially built for the first subsystem, i.e. a Delta impulse~$\delta_{\varepsilon_1}(t)$ acts foremost on the first subsystem. Once the nascent impulse~$\delta_{\varepsilon_1}(t)$ has affected the first subsystem (for $t > \varepsilon_1$), the solution~$x_1(t)$ and the nascent impulse~$\delta_{\varepsilon_2}(t)$ affect the second subsystem. Thus, building multiple limits~$\varepsilon_k \rightarrow 0$ corresponds to the successive action of Dirac impulses on the subsystems.

Interestingly enough, applying the sifting property of the Dirac impulse 
\begin{equation}
\begin{aligned}  
	\int_{a}^{b}f(\tau)\delta(\tau)\,\dd\tau = \left\{ \begin{array}{ll}
		f(0), & 0 \in [a, b] \\ 0, & \mathrm{else}
	\end{array}  \right.
\end{aligned}
\end{equation}
successively results in the exact same procedure: the impulse~$\delta(t)$ acts initially on the first subsystem, and afterwards the solution~$x_1(t)$ and the impulse~$\delta(t)$ affect the second subsystem. Therefore, the solution~$x_2(t)$ of the second subsystem for~\mbox{$u(t) = \delta(t)$} when applying the sifting property would be again: 
\begin{equation}
\begin{aligned} \label{eq:x2-suc-sifting}
x_2(t)
&= \int_{\tau_1 = 0}^t\int_{\tau_2 = 0}^{\tau_1}  \eu^{a(t - \tau_1)} n \delta(\tau_1) 
\eu^{a(\tau_1 - \tau_2)} b \delta(\tau_2) \, \dd \tau_2\dd \tau_1 \\
&= \eu^{at}nb \int_{\tau_1 = 0}^t\int_{\tau_2 = 0}^{\tau_1} \eu^{-a\tau_2} \delta(\tau_2)\dd\tau_2 \delta(\tau_1) \,\dd \tau_1 \\
&= \eu^{at}nb \int_{\tau_1 = 0}^t \delta(\tau_1) \,\dd \tau_1\, \\
&= \eu^{at}nb\, .
\end{aligned}
\end{equation}
Note that the factor~$1/k!$ does neither arise in~\eqref{eq:x2-suc-limit} nor in~\eqref{eq:x2-suc-sifting}. This factor will, however, appear again in the following. 

\subsubsection*{Simultaneous action of Dirac impulse}
The input of the bilinear system~$u(t)$ actually acts on all subsystems simultaneously. Therefore, we believe that the nascent delta function~$\delta_{\varepsilon}(t)$ affects all subsystems at the same time and, consequently, only a single limit~$\varepsilon \rightarrow 0$ should be performed. In such case, the solution~$x_2(t)$ of the second subsystem for~$u(t) = \delta_\varepsilon(t)$ in the limit~$\varepsilon \rightarrow 0$ for $t > \varepsilon \geq \tau_1 \geq \tau_2$ is given by:
\begin{equation}
\begin{aligned}
\label{eq:zusatz01}
	x_2(t) &= \eu^{at}nb\lim\limits_{\varepsilon \rightarrow 0}\int_{\tau_1 = 0}^{\varepsilon}\varepsilon^{-1}\int_{\tau_2 = 0}^{\tau_1}\varepsilon^{-1}\eu^{-a\tau_2}\,\dd \tau_2 \dd \tau_1 = \eu^{at}nb \lim\limits_{\varepsilon \rightarrow 0}\varepsilon^{-2}a^{-1}\int_{\tau_ 1 = 0}^{\varepsilon}1 - \eu^{-a\tau_1}\,\dd \tau_1 \\[0.3em]
	&= \eu^{at}nb \lim\limits_{\varepsilon \rightarrow 0} \frac{a \varepsilon + \eu^{-a\varepsilon} - 1 }{a^2 \varepsilon^2} \overset{\mathrm{L.H.}}{=}\eu^{at}nb \lim\limits_{\varepsilon \rightarrow 0}\frac{a - a \eu^{a\varepsilon}}{2 a^2 \varepsilon} \overset{\mathrm{L.H.}}{=}\eu^{at}nb\lim\limits_{\varepsilon \rightarrow 0}\frac{a^2\eu^{-a\varepsilon}}{2a^2} \\[0.3em]
	&= \tfrac{1}{2}\eu^{at}nb\, .
\end{aligned}
\end{equation}
With the simultaneous action of the Dirac impulse on all subsystems, the factor~$1/k!$ appears in the solution. The derivation for higher subsystems and the matrix case is extensively performed in~\cite[\S3.5]{gebhart2017impulse}. Here, we simply state the final result for the impulse response of the~$k$-th subsystem: 
\\[1ex]
\fbox{\parbox{14.7cm}{
\begin{equation}
\begin{aligned} \label{eq:imp-res-kth}
\vec g_k(t; \vec \mu) &:= \lim\limits_{\varepsilon \rightarrow 0}\vec y_k(t)\ \ \text{subject to} \ \ \vec u(t) = \vec \mu \delta_\varepsilon(t)\\
&= \underbrace{\mat C \eu^{\mat A t}\frac{\hat{\mat N}^{k-1}}{k!}\hat{\vec b}}_{\vec g_{\mat B, k}(t)} \ \ + \ \ \underbrace{\mat C \eu^{\mat A t}\frac{\hat{\mat N}^{k-1}}{(k-1)!}\vec x_0}_{\vec g_{\mathrm{IC}, k}(t)}, \quad \hat{\mat N} = \sum_{j = 1}^{m}\mat N_j \mu_j, \quad \hat{\vec b} = \mat B \vec \mu = \sum_{j=1}^{m} \vec b_j \mu_j\, .
\end{aligned}
\end{equation}
}}\\[1ex]

In \cite[\S3.4]{gebhart2017impulse} an error analysis for the impulse response of bilinear systems has been performed. It shows that the output~$y_{\mathrm{blin}}(t)$ of the (scalar) bilinear system subject to~$\delta_\varepsilon(t)$ converges to the impulse response~$g_{\mathrm{blin}}(t)$ of the bilinear system like the output~$y_{\mathrm{lin}}(t)$ of the linear system subject to~$\delta_\varepsilon(t)$ converges to the impulse response~$g_{\mathrm{lin}}(t)$ of the linear system.
\section{Adjustment of the multidimensional Volterra kernels} \label{sec:5}
In this section, we will adjust the definitions of the triangular and regular kernels especially along lines of equal time arguments, such that the kernels and the previously derived impulse response are compatible to each other:~$g_k(t) = g_k^\triangle(t,\ldots,t) = g_k^\square(0,\ldots,0,t)$.

\subsection{Multidimensional symmetric kernels} \label{subsec:sym-kernels}
Until now, the triangular and regular kernels have been used to describe a bilinear system. However, there is another form of interest, the so-called \emph{symmetric} kernel, that will be discussed in the following.

\subsubsection*{Second subsystem}
According to \cite[p. 12]{rugh1981nonlinear}, the symmetric output of the second subsystem $y^{\mathrm{sym}}_2(t_1,t_2)$ can be written by summing over the triangular outputs~$y^{\triangle}_2(t_1,t_2)$ and~$y^{\triangle}_2(t_2,t_1)$:
\begin{equation}
\begin{aligned} \label{eq:y2-sym}
y^{\mathrm{sym}}_2(t_1,t_2) := \tfrac{1}{2}\big(y^{\triangle}_2(t_1,t_2) + y^{\triangle}_2(t_2,t_1) \big)\, .
\end{aligned}
\end{equation}
The output is called \emph{symmetric}, since~$y^{\mathrm{sym}}_2(t_1,t_2) \!=\! y^{\mathrm{sym}}_2(t_2,t_1)$. Then the single-variable output of the second subsystem~$y_2(t)$ can be calculated from:
\begin{equation}
\begin{aligned}
	y_2(t) &= y_2^\triangle(t,t) = \frac{1}{2}(y_2^\triangle(t,t) + y_2^\triangle(t,t)) \\
		   &= y_2^{\mathrm{sym}}(t,t)\,.
\end{aligned}
\end{equation}
The Laplace transform of the triangular output~$y^{\triangle}_2(t_2,t_1)$ with permuted time variables is given by the expression~$\mathcal{L}\{y^{\triangle}_2(t_2,t_1)\}(s_1,s_2) \!=\! Y^\triangle_2(s_2,s_1)$, since the time variables~$t_1,t_2$ are linked to the frequency variables~$s_1,s_2$ by the exponential terms~$\eu^{-s_1t_1},\eu^{-s_2t_2}$. Consequently, the Laplace transform of the two-dimensional symmetric output becomes:
\begin{equation}
\begin{aligned}
Y_2^{\mathrm{sym}}(s_1,s_2) &:= \mathcal{L}\{y^{\mathrm{sym}}_2(t_1,t_2)\}(s_1,s_2) \\
&= \tfrac{1}{2}Y_2^\triangle(s_1,s_2) + \tfrac{1}{2}Y_2^\triangle(s_2,s_1)\\
&= \tfrac{1}{2}\vec c^{\mathsf T} \big((s_1 + s_2)\matn I - \mat A \big)^{-1}\mat N (s_1 \matn I - \mat A)^{-1}\vec b \, U(s_1)U(s_2) \\
&+ \tfrac{1}{2}\vec c^{\mathsf T} \big((s_2 + s_1)\matn I - \mat A \big)^{-1}\mat N (s_2 \matn I - \mat A)^{-1}\vec b \, U(s_2)U(s_1) \\
&= G^{\mathrm{sym}}_2(s_1,s_2) \, U(s_1)U(s_2)\, ,
\end{aligned}
\end{equation}
with the two-dimensional symmetric transfer function
\begin{equation}
\begin{aligned}
G_2^{\mathrm{sym}}(s_1,s_2) &= \tfrac{1}{2}\big(G_2^\triangle(s_1,s_2) + G_2^\triangle(s_2,s_1)\big) \\
&= \tfrac{1}{2}\vec c^{\mathsf T} \big((s_1 + s_2)\matn I - \mat A \big)^{-1}\mat N \left[(s_1 \matn I - \mat A)^{-1}\vec b + (s_2 \matn I - \mat A)^{-1}\vec b \right]. 
\end{aligned}
\end{equation}
The two-dimensional transfer function $G_2^{\mathrm{sym}}(s_1,s_2) = Y_2^{\mathrm{sym}}(s_1,s_2)$ s.t.~$U(s) = 1$ corresponds to the Laplace transform of the symmetric impulse response~$g^{\mathrm{sym}}_2(t_1,t_2) = y^{\mathrm{sym}}_2(t_1,t_2)$ s.t.~$u(t)=\delta(t)$, as the Laplace transform~$\mathcal{L}\{\delta(t)\} = 1$:
\begin{equation}
\begin{aligned}
G^{\mathrm{sym}}_2(s_1,s_2) = \mathcal{L}\{g^{\mathrm{sym}}_2(t_1,t_2)\}(s_1,s_2)\, .
\end{aligned}
\end{equation}
According to the definition of the triangular kernel~$g^{\triangle}_2(t_1,t_2)$ from equation~\eqref{eq:tri-SISO-kernels}, the two-dimensional symmetric impulse response becomes for~$t_1, t_2 > 0$ (cf. Figure~\ref{fig:defdsym}):
\begin{equation}
\begin{aligned}
g^{\mathrm{sym}}_2(t_1,t_2) &= \tfrac{1}{2}\big(g^{\triangle}_2(t_1,t_2) +  g^{\triangle}_2(t_2,t_1)\big) \\[0.5em]
&= \frac{1}{2}\underbrace{\left\{\begin{array}{ll} \vec c^{\mathsf T} \eu^{\mat A t_2}\mat N \eu^{\mat A (t_1 - t_2)}\vec b, & t_1 > t_2 \\
	\text{not yet defined}, & t_1 = t_2 \\
	0, & t_1 < t_2 \\
	\end{array} \right.}_{g^{\triangle}_2(t_1,t_2)} \ \ + \ \ \frac{1}{2}\underbrace{\left\{\begin{array}{ll} \vec c^{\mathsf T} \eu^{\mat A t_1}\mat N \eu^{\mat A (t_2 - t_1)}\vec b, & t_2 > t_1 \\
	\text{not yet defined}, & t_2 = t_1 \\
	0, & t_2 < t_1 \\
	\end{array} \right.}_{g^{\triangle}_2(t_2,t_1)} \\[0.5em]
&= \left\{\begin{array}{ll} \frac{1}{2}\vec c^{\mathsf T} \eu^{\mat A t_2}\mat N \eu^{\mat A (t_1 - t_2)}\vec b, & t_1 > t_2 \\
\text{not yet defined}, & t_1 = t_2 \\
\frac{1}{2}\vec c^{\mathsf T} \eu^{\mat A t_1}\mat N \eu^{\mat A (t_2 - t_1)}\vec b, & t_1 < t_2\, . \\
\end{array} \right.
\end{aligned}
\end{equation}
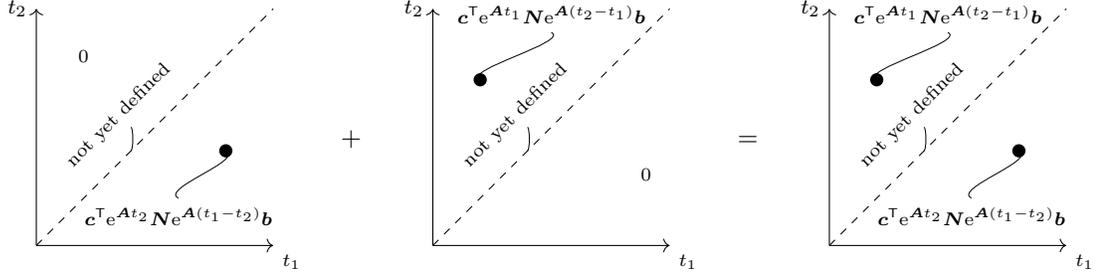
\begin{figure}[htp!]
	\begin{center}
		\begin{tabular}{m{0.273\linewidth}m{0.02\linewidth}m{0.273\linewidth}m{0.02\linewidth}m{0.273\linewidth}}
			\begin{tikzpicture}
[xscale=1/(10 - (-0.5))*0.0352778*0.44*\textwidth,
yscale=1/(10- (-0.5))*0.0352778*0.44*\textwidth,
xtick/.style={below,fill=white,font=\tiny},
ytick/.style={left=0.2,fill=white,font=\tiny},
xlabel/.style={below right,font=\scriptsize},
ylabel/.style={left,font=\footnotesize},
beziernode/.style={circle, inner sep=0pt,minimum size=5pt,fill=black}
]
\def\tickin{0.1}
\def\tickout{0.1}
\def\delta{0.5}
\def\xmax{5}
\def\xmin{0}
\def\ymax{5} 
\def\ymin{0}

\begin{scope}
\draw [->] (0,\ymin) -- (0,\ymax) node [ylabel] {$t_2$};
\draw [->] (\xmin,0) -- (\xmax,0) node [xlabel] {$t_1$};
\draw [dashed] (\xmin,\xmin) -- (\ymax,\ymax);

\draw [ultra thin]  (0.2*\ymax,0.8*\ymax) node[font=\scriptsize] {$0$};

\draw [ultra thin] (0.4*\ymax,0.4*\ymax) .. controls +(+0.05,+0.05) and ++(+0.05,-0.05) .. (0.4*\ymax,0.5*\ymax)  node[anchor=south,font=\scriptsize, rotate = 45] {not yet defined};

\node[beziernode] at (0.8*\ymax,0.4*\ymax) {};
\draw[ultra thin] (0.8*\ymax,0.4*\ymax) .. controls +(+0.4,-0.2) and ++(-0.4,+0.2) .. (0.6*\ymax,0.2*\ymax) node[anchor=north, font=\scriptsize] {$\vec c^{\mathsf T} \eu^{\mat A t_2}\mat N \eu^{\mat A (t_1 - t_2)}\vec b$};
\end{scope}

\end{tikzpicture}%
 & + & \begin{tikzpicture}
[xscale=1/(10 - (-0.5))*0.0352778*0.44*\textwidth,
yscale=1/(10- (-0.5))*0.0352778*0.44*\textwidth,
xtick/.style={below,fill=white,font=\tiny},
ytick/.style={left=0.2,fill=white,font=\tiny},
xlabel/.style={below right,font=\scriptsize},
ylabel/.style={left,font=\footnotesize},
beziernode/.style={circle, inner sep=0pt,minimum size=5pt,fill=black}
]
\def\tickin{0.1}
\def\tickout{0.1}
\def\delta{0.5}
\def\xmax{5}
\def\xmin{0}
\def\ymax{5} 
\def\ymin{0}

\begin{scope}[xshift = 7cm]
\draw [->] (0,\ymin) -- (0,\ymax) node [ylabel] {$t_2$};
\draw [->] (\xmin,0) -- (\xmax,0) node [xlabel] {$t_1$};
\draw [dashed] (\xmin,\xmin) -- (\ymax,\ymax);

\node[beziernode] at (0.2*\ymax,0.7*\ymax) {};
\draw [ultra thin] (0.2*\ymax,0.7*\ymax) .. controls +(-0.4,0.2) and ++(0.4,-0.2) .. (0.5*\ymax,0.9*\ymax)  node[anchor = south, font=\scriptsize] {$\vec c^{\mathsf T} \eu^{\mat A t_1}\mat N \eu^{\mat A (t_2 - t_1)}\vec b$};
\draw [ultra thin] (0.4*\ymax,0.4*\ymax) .. controls +(+0.05,+0.05) and ++(+0.05,-0.05) .. (0.4*\ymax,0.5*\ymax)  node[anchor=south,font=\scriptsize, rotate = 45] {not yet defined};
\draw [ultra thin] (0.9*\ymax,0.3*\ymax) node[font=\scriptsize] {$0$};
\end{scope}

\end{tikzpicture}%
 & = & \begin{tikzpicture}
[xscale=1/(10 - (-0.5))*0.0352778*0.44*\textwidth,
yscale=1/(10- (-0.5))*0.0352778*0.44*\textwidth,
xtick/.style={below,fill=white,font=\tiny},
ytick/.style={left=0.2,fill=white,font=\tiny},
xlabel/.style={below right,font=\scriptsize},
ylabel/.style={left,font=\footnotesize},
beziernode/.style={circle, inner sep=0pt,minimum size=5pt,fill=black}
]
\def\tickin{0.1}
\def\tickout{0.1}
\def\delta{0.5}
\def\xmax{5}
\def\xmin{0}
\def\ymax{5} 
\def\ymin{0}

\begin{scope}[xshift = 14cm]
\draw [->] (0,\ymin) -- (0,\ymax) node [ylabel] {$t_2$};
\draw [->] (\xmin,0) -- (\xmax,0) node [xlabel] {$t_1$};
\draw [dashed] (\xmin,\xmin) -- (\ymax,\ymax);
\node[beziernode] at (0.2*\ymax,0.7*\ymax) {};
\draw [ultra thin] (0.2*\ymax,0.7*\ymax) .. controls +(-0.4,0.2) and ++(0.4,-0.2) .. (0.5*\ymax,0.9*\ymax)  node[anchor = south, font=\scriptsize] {$\vec c^{\mathsf T} \eu^{\mat A t_1}\mat N \eu^{\mat A (t_2 - t_1)}\vec b$};
\draw [ultra thin] (0.4*\ymax,0.4*\ymax) .. controls +(+0.05,+0.05) and ++(+0.05,-0.05) .. (0.4*\ymax,0.5*\ymax)  node[anchor=south,font=\scriptsize, rotate = 45] {not yet defined};
\node[beziernode] at (0.8*\ymax,0.4*\ymax) {};
\draw[ultra thin] (0.8*\ymax,0.4*\ymax) .. controls +(+0.4,-0.2) and ++(-0.4,+0.2) .. (0.6*\ymax,0.2*\ymax) node[anchor=north, font=\scriptsize] {$\vec c^{\mathsf T} \eu^{\mat A t_2}\mat N \eu^{\mat A (t_1 - t_2)}\vec b$};
\end{scope}

\end{tikzpicture}%
		\end{tabular}
	\end{center}
	\caption{Domain of the symmetric kernel~$g^{\mathrm{sym}}_2(t_1,t_2)$ of the second subsystem. -- Due to the same left- and right-handed limit values, the symmetric kernel is continuously extended at the discontinuity $t_1 = t_2 = t$.}
	\label{fig:defdsym}
\end{figure}
The two-dimensional symmetric impulse response~$g^{\mathrm{sym}}_2(t_1,t_2)$ possesses for~$t_1 > t_2$ and~$t_1 < t_2$ the same (left- and right-handed) limit value as~$t_1,t_2 \rightarrow t$, since the transformation~$g_2^\triangle(t_1,t_2) \rightarrow g_2^\triangle(t_2,t_1)$ can be interpreted as mirroring~$g_2^\triangle(t_1,t_2)$ at the line~$t_1 = t_2$: 
\begin{equation}
\begin{aligned}
g^{\mathrm{sym}}_2(t_1 \rightarrow t^+,t) = g^{\mathrm{sym}}_2(t_1 \rightarrow t^-,t) = \frac{1}{2}\vec c^{\mathsf T} \eu^{\mat A t}\mat N\vec b\, .
\end{aligned}
\end{equation}
Therefore it seems reasonable to continuously extend the symmetric impulse response at~$t_1 = t_2 = t$. This finally yields:
\begin{equation}
\begin{aligned}
g^{\mathrm{sym}}_2(t_1,t_2)
&:= \left\{\begin{array}{ll} \frac{1}{2}\vec c^{\mathsf T} \eu^{\mat A t_2}\mat N \eu^{\mat A (t_1 - t_2)}\vec b, & t_1 > t_2 > 0 \\[0.2em]
\frac{1}{2}\vec c^{\mathsf T} \eu^{\mat A t}\mat N \vec b, & t_1 = t_2 = t > 0\\[0.2em]
\frac{1}{2}\vec c^{\mathsf T} \eu^{\mat A t_1}\mat N \eu^{\mat A (t_2 - t_1)}\vec b, & t_2 > t_1 > 0 \\[0.2em]
\end{array} \right. \\[0.5em]
&:= \tfrac{1}{2}\vec c^{\mathsf T} \eu^{\mat A t_{i_2}}\mat N \eu^{\mat A (t_{i_1} - t_{i_2})}\vec b,\:\: t_{i_1} \geq t_{i_2} > 0,\:\: (i_1,i_2) \in \{(1,2), (2,1)\}\, .
\end{aligned}
\end{equation} 
This definition of the two-dimensional symmetric kernel is consistent with the impulse response of the second subsystem $g_2(t) = g^{\mathrm{sym}}_2(t,t) = 1/2 \, \vec c^{\mathsf T} \eu^{\mat A t}\mat N \vec b$ according to equation~\eqref{eq:imp-res-kth}.

We will now adjust the definition of the two-dimensional triangular kernel by defining~$g^{\triangle}_{2}(t_1,t_2)$ by the equal share of the symmetric kernel~$g^{\mathrm{sym}}_{2}(t_1,t_2)$ along the time line~$t_1 = t_2 = t$:
\begin{equation}
\begin{aligned} 
&& g^{\mathrm{sym}}_{2}(t,t) &= \tfrac{1}{2}\big( g^{\triangle}_{2}(t,t) + g^{\triangle}_{2}(t,t)\big)\:\:\text{and}\:\: g^{\mathrm{sym}}_{2}(t,t) = \tfrac{1}{2}\vec c^\intercal \eu^{\mat A t}\mat N \vec b \\
&\implies& g^{\triangle}_{2}(t,t) &= \tfrac{1}{2}\vec c^\intercal \eu^{\mat A t}\mat N \vec b\,,
\end{aligned}
\end{equation}
since both triangles border the intersection line~$t_1 \!=\! t_2$ ($n\!=\!1$-simplex). Consequently, the two-dimensional triangular impulse response is finally defined by
\begin{equation}
\begin{aligned} \label{eq:g2-tri}
g^{\triangle}_{2}(t_1,t_2) 
&:= \left\{\begin{array}{ll}
\vec c^{\mathsf T} \eu^{\mat A t_2}\mat N \eu^{\mat A (t_1 - t_2)}\vec b, & t_1 > t_2 > 0, \ (n=2)\\[0.2em]
\frac{1}{2}\vec c^{\mathsf T} \eu^{\mat A t}\mat N \vec b, & t_1 = t_2 = t > 0, \ (n=1)\\[0.2em]
0, & \text{else} 
\end{array}  \right. \\
&:=\tfrac{1}{(2+1-n)!}\vec c^{\mathsf T} \eu^{\mat A t_2} \mat N \eu^{\mat A(t_1 - t_2)} \vec b\,, \ \ t_1 \geq t_2 > 0\, ,
\end{aligned}
\end{equation}
where~$n$ is the dimension of the corresponding region, i.e.~$n \!=\! 1$ for the time line~$t_1 \!=\! t_2 \!=\! t$ (1-simplex) and~$n \!=\! 2$ for the triangle~$t_1 > t_2 > 0$ (2-simplex).

Respectively, the two-dimensional regular kernel~$g^{\square}_{2}(t_1,t_2)$ can be defined along the time line $t_1 = 0, t_2 = t$ by means of the triangular kernel according to the transformation from regular to triangular from section~\ref{subsubsec:reg-kernel-rep}. Thus, the two-dimensional regular kernel becomes
\begin{equation}
\begin{aligned}
g^{\square}_{2}(t_1,t_2) &= g^{\triangle}_{2}(t_1 + t_2,t_2) \\
&:= \left\{\begin{array}{ll}
\vec c^{\mathsf T} \eu^{\mat A t_2}\mat N \eu^{\mat A t_1}\vec b, & t_1, t_2 > 0, \ (n=2)\\[0.2em]
\frac{1}{2}\vec c^{\mathsf T} \eu^{\mat A t} \mat N \vec b, & t_1 = 0, t_2 = t > 0, \ (n=1) \\[0.2em]
0, & \mathrm{else}
\end{array}  \right. \\
&:= \tfrac{1}{(2+1-n)!}\vec c^{\mathsf T} \eu^{\mat A t_2} \mat N \eu^{\mat A t_1} \vec b\,, \ \ t_1 \geq 0,t_2 > 0\, .
\end{aligned}
\end{equation}

\subsubsection*{Third subsystem}
The symmetric output of the third subsystem $y^{\mathrm{sym}}_3(t_1,t_2,t_3)$ can be obtained by summing over all permutations~$\pi(\cdot)$ of the time variables~$(t_1,t_2,t_3)$ of the triangular output~$y^{\triangle}_3(t_1,t_2,t_3)$:
\begin{equation}
\begin{aligned}
y^{\mathrm{sym}}_3(t_1,t_2,t_3) = \tfrac{1}{6}\big(y^\triangle_3(t_1,t_2,t_3) + y^\triangle_3(t_1,t_3,t_2) + \ldots \big)
&= \tfrac{1}{6}\sum_{\pi(\cdot)}y^{\triangle}_3(t_{\pi(1)},t_{\pi(2)},t_{\pi(3)})\,,
\end{aligned}
\end{equation}
where~$(\pi(1),\pi(2),\pi(3)) \in \{(1,2,3), (1,3,2), (2,1,3), (2,3,1), (3,1,2), (3,2,1) \}$. The three-\-di\-men\-sio\-nal symmetric output is symmetric in each time argument $t_1, t_2, t_3$, since~$y^{\mathrm{sym}}_3(t_1,t_2,t_3) \!=\! y^{\mathrm{sym}}_3(t_1,t_3,t_1) \!=\! \ldots \!=\! y^{\mathrm{sym}}_3(t_3,t_2,t_1)$. The single-variable output of the third subsystem~$y_3(t)$ can be calculated from:~$y_3(t) \!=\! y^\triangle_3(t,t,t) \!=\! y^{\mathrm{sym}}_3(t,t,t)$.

\begin{figure}[tp!]	
	\definecolor{mycolor1}{rgb}{0.97690,0.98390,0.08050}%
	\definecolor{mycolor2}{rgb}{0.24220,0.15040,0.66030}%
	\definecolor{mycolor3}{rgb}{0.09640,0.75000,0.71204}%
	\begin{center}
				\definecolor{mycolor1}{rgb}{0.97690,0.98390,0.08050}%
		\definecolor{mycolor2}{rgb}{0.24220,0.15040,0.66030}%
		\definecolor{mycolor3}{rgb}{0.09640,0.75000,0.71204}%

\begin{tikzpicture}[
helpline/.style={dashed},
]

\begin{axis}[
legend columns=-1,
legend entries={$t_2 = t_3$,$t_1 = t_2$,$t_1 = t_3$},
legend to name=named,
  width=0.44*\textwidth,
  tick label style={font=\tiny, fill=white},
  label style={font=\color{white!15!black}, font=\scriptsize},
  axis on top,
  scale only axis,
  tick align=outside,
  axis lines=center,
  y dir=reverse,
  name=myplot,
  view={22}{35.6},
  enlargelimits=.1,
  ymin=0, ymax=1, xmin=0, xmax=1, zmin=0, zmax=1,
  xlabel=$t_1$, ylabel=$t_2$, zlabel=$t_3$,
  ticks=none,
  every axis x label/.append style={at=(ticklabel* cs:1.03)},
  every axis y label/.append style={at=(ticklabel* cs:0.2)},
  every axis z label/.append style={at=(ticklabel* cs:1.05)}
  ]

\draw [helpline] (1,0,0) -- (1,1,0);
\draw [helpline] (0,1,0) -- (1,1,0);
\draw [helpline] (0,0,1) -- (1.0,0,1);
\draw [helpline] (0,0,1) -- (0,1,1);
\draw [helpline] (1,0,0) -- (1,0,1);
\draw [helpline] (0,1,0) -- (0,1,1);

\draw (0,0,0) -- (1,1,1);

%
%

\addplot3[area legend, draw=black, fill=mycolor2, fill opacity=0.5]
table[row sep=crcr] {%
	x	y	z\\
	0	0	0\\
	1	0	0\\
	1	1	1\\
}--cycle;
\addlegendentry{$t_2 = t_3$}

\addplot3[area legend, draw=black, fill=mycolor1, fill opacity=0.5]
table[row sep=crcr] {%
	x	y	z\\
	0	0	0\\
	1	1	0\\
	1	1	1\\
}--cycle;
\addlegendentry{$t_1 = t_2$} 

\addplot3[area legend, draw=black, fill=mycolor3, fill opacity=0.5]
table[row sep=crcr] {%
	x	y	z\\
	0	0	0\\
	1	1	1\\
	0	1	0\\
}--cycle;
\addlegendentry{$t_1 = t_3$}

\addplot3[area legend, draw=black, fill=mycolor3, fill opacity=0.5, forget plot]
table[row sep=crcr] {%
	x	y	z\\
	0	0	0\\
	1	0	1\\
	1	1	1\\
}--cycle;
\addlegendentry{türkis oben}

\addplot3[area legend, draw=black, fill=mycolor1, fill opacity=0.5, forget plot]
table[row sep=crcr] {%
	x	y	z\\
	0	0	0\\
	1	1	1\\
	0	0	1\\
}--cycle;
\addlegendentry{gelb oben} 

\addplot3[area legend, draw=black, fill=mycolor2, fill opacity=0.5, forget plot]
table[row sep=crcr] {%
	x	y	z\\
	0	0	0\\
	1	1	1\\
	0	1	1\\
}--cycle;
\addlegendentry{lila oben}

\end{axis}

\end{tikzpicture}
		\ref{named}			
	\end{center}
	\caption{Domain of the symmetric kernel~$g^{\mathrm{sym}}_3(t_1,t_2,t_3)$ of the third subsystem. -- It consists of six tetrahedra. Due to the same left- and right-handed limit values, the symmetric kernel is continuously extended at the discontinuities $t_2 = t_3$,~$t_1 = t_2$ and~$t_1 = t_3$.}
	\label{fig:3-vis}
\end{figure}
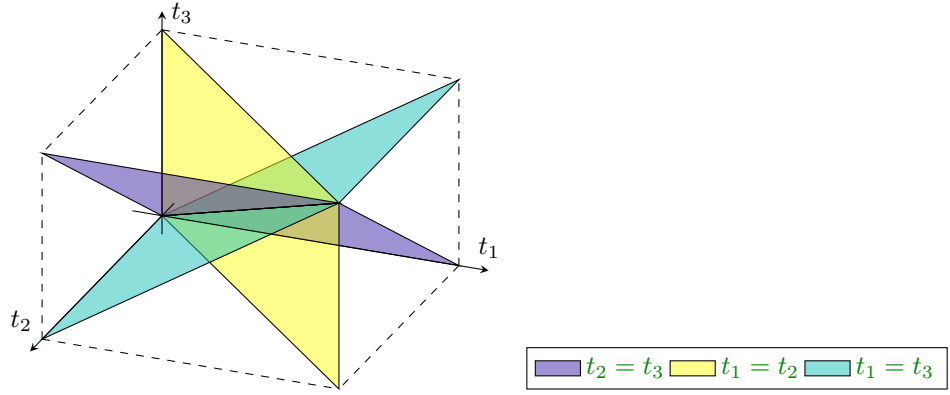

The Laplace transform of the three-dimensional symmetric output~$Y^{\mathrm{sym}}_3(s_1,s_2,s_3)$ becomes
\begin{equation}
\begin{aligned}
Y^{\mathrm{sym}}_3(s_1,s_2,s_3) &= \mathcal{L}\{ y^{\mathrm{sym}}_3(t_1,t_2,t_3)\}(s_1,s_2,s_3) \\
&= \mathcal{L}\left\{ \tfrac{1}{6}\big(y^\triangle_3(t_1,t_2,t_3) + y^\triangle_3(t_1,t_3,t_2) + \ldots \big)\right\}(s_1,s_2,s_3) \\
&= \tfrac{1}{6}\big(Y^\triangle_3(s_1,s_2,s_3) + Y^\triangle_3(s_1,s_3,s_2) + \ldots \big)\\
&= G^{\mathrm{sym}}_{3}(s_1,s_2,s_3) \, U(s_1)U(s_2)U(s_3)\,,
\end{aligned}
\end{equation}
with the three-dimensional symmetric transfer function
\begin{equation}
\begin{aligned}
\lefteqn{G^{\mathrm{sym}}_3(s_1,s_2,s_3) = \tfrac{1}{6}\big(G^\triangle_3(s_1,s_2,s_3) + G^\triangle_3(s_1,s_3,s_2) + \ldots \big)} \\
&= \tfrac{1}{6}\vec c^{\mathsf T} \big((s_1 + s_2 + s_3)\matn I - \mat A \big)^{-1}\mat N \left[ \big((s_1 + s_2)\matn I - \mat A\big)^{-1} \mat N \left[(s_1 \matn I - \mat A)^{-1}\vec b + (s_2 \matn I - \mat A)^{-1}\vec b \right] \right.  \\
& \qquad \qquad \qquad \qquad \qquad \qquad \qquad + \big((s_2 + s_3)\matn I - \mat A\big)^{-1} \mat N \left[(s_2 \matn I - \mat A)^{-1}\vec b + (s_3 \matn I - \mat A)^{-1}\vec b \right] \\
& \qquad \qquad \qquad \qquad \qquad \qquad \qquad + \left.\big((s_1 + s_3)\matn I - \mat A\big)^{-1} \mat N \left[(s_1 \matn I - \mat A)^{-1}\vec b + (s_3 \matn I - \mat A)^{-1}\vec b \right] \right].
\end{aligned}
\end{equation}
The transfer function $G_3^{\mathrm{sym}}(s_1,s_2,s_3) = Y_3^{\mathrm{sym}}(s_1,s_2,s_3)$ s.t.~$U(s) = 1$ corresponds to the Laplace transform of the symmetric impulse response $g^{\mathrm{sym}}_3(t_1,t_2,t_3) = y^{\mathrm{sym}}_3(t_1,t_2,t_3)$ s.t.~$u(t)=\delta(t)$:
\begin{equation}
\begin{aligned}
G^{\mathrm{sym}}_3(s_1,s_2,s_3) = \mathcal{L}\{g^{\mathrm{sym}}_3(t_1,t_2,t_3)\}(s_1,s_2,s_3)\, .
\end{aligned}
\end{equation}
The symmetric impulse response of the third subsystem~$g^{\mathrm{sym}}_3(t_1,t_2,t_3)$ is the sum of six triangular impulse responses~$g^\triangle_3(t_1,t_2,t_3)$,~$g^\triangle_3(t_1,t_3,t_2)$, etc, which are non-zero on a corresponding tetrahedron~$t_1 > t_2 > t_3$,~$t_1 > t_3 > t_2$, etc (cf. Figure~\ref{fig:3-vis}). The symmetric impulse response can be continuously extended at the surface intersections~$t_1 = t_2$,~$t_1 = t_3$ and~$t_2 = t_3$ of the six tetrahedra, as well as at the intersection line~$t_1 = t_2 = t_3$ of the surfaces, since the limit values are equal, as the surface can be interpreted as mirror face. Therefore: 
\begin{equation}
\begin{aligned}
g^{\mathrm{sym}}_3(t_1,t_2,t_3)
&= \tfrac{1}{6}\vec c^{\mathsf T} \eu^{\mat A t_{i_3}}\mat N \eu^{\mat A (t_{i_2} - t_{i_3})}\mat N \eu^{\mat A (t_{i_1} - t_{i_2})} \vec b, \\
&\qquad t_{i_1} \geq t_{i_2} \geq t_{i_3} > 0,\:\: (i_1,i_2,i_3) \text{ permutations of } (1,2,3)\, .
\end{aligned}
\end{equation}
This definition of the three-dimensional symmetric kernel is consistent with the impulse response of the third subsystem $g_3(t) = g^{\mathrm{sym}}_3(t,t,t) = 1/6 \, \vec c^{\mathsf T} \eu^{\mat A t}\mat N^2 \vec b$ according to equation~\eqref{eq:imp-res-kth}.

All six tetrahedra border the intersection line~$t_1 \!=\! t_2 \!=\! t_3$ ($n\!=\!1$-simplex), whereas two tetrahedra border each cutting triangle ($n\!=\!2$-simplex). Therefore, the three-dimensional triangular impulse response is
\begin{equation}
\begin{aligned}
g^{\triangle}_{3}(t_1,t_2,t_3) &:= \left\{\begin{array}{ll}
\vec c^{\mathsf T} \eu^{\mat A t_3}\mat N\eu^{\mat A(t_2 - t_3)}\mat N \eu^{\mat A(t_1 - t_2)}, & t_1 > t_2 > t_3 > 0, \ (n=3) \\[0.3em]
\tfrac{1}{6}\vec c^{\mathsf T} \eu^{\mat A t}\mat N^2 \vec b, & t_1 = t_2 = t_3 = t > 0, \ (n=1) \\[0.3em]
\tfrac{1}{2}\vec c^{\mathsf T} \eu^{\mat A t_3}\mat N\eu^{\mat A(t - t_3)}\mat N\vec b, & t_1 = t_2 = t > t_3 > 0, \ (n=2) \\[0.3em]
\tfrac{1}{2}\vec c^{\mathsf T} \eu^{\mat A t}\mat N^2\eu^{\mat A(t_1 - t)} \vec b, & t_1 > t_2 = t_3 = t > 0, \ (n=2) \\[0.3em]
0, & \mathrm{else}
\end{array}  \right. \\
&:= \tfrac{1}{(3+1-n)!}\vec c^{\mathsf T} \eu^{\mat A t_3} \mat N \eu^{\mat A(t_2 - t_3)} \mat N \eu^{\mat A(t_1 - t_2)} \vec b\,, \ \ t_1 \geq t_2 \geq t_3 > 0\,.
\end{aligned}
\end{equation}

Similarly, the regular impulse response of the third subsystem is given by
\begin{equation}
\begin{aligned}
g^{\square}_{3}(t_1,t_2,t_3) &= g^{\triangle}_{3}(t_1 + t_2,t_2+t_3,t_3) \\
&:= \left\{\begin{array}{ll}
\vec c^{\mathsf T} \eu^{\mat A t_3} \mat N \eu^{\mat A t_2}\mat N \eu^{\mat A t_1}\vec b, & t_1, t_2, t_3 > 0, \ (n=3)\\[0.3em]
\tfrac{1}{6}\vec c^{\mathsf T} \eu^{\mat A t}\mat N^2 \vec b, & t_1 = t_2 = 0, t_3 = t > 0, \ (n=1) \\[0.3em]
\tfrac{1}{2}\vec c^{\mathsf T} \eu^{\mat A t_3} \mat N \eu^{\mat A t_2} \mat N \vec b, & t_1 = 0, t_2, t_3 > 0, \ (n=2) \\[0.3em]
\tfrac{1}{2}\vec c^{\mathsf T} \eu^{\mat A t_3} \mat N^2 \eu^{\mat A t_1} \vec b, & t_2 = 0, t_1, t_3 > 0, \ (n=2) \\[0.3em]
0, & \text{else}
\end{array}  \right. \\
&:= \tfrac{1}{(3+1-n)!}\vec c^{\mathsf T} \eu^{\mat A t_3} \mat N \eu^{\mat A t_2} \mat N \eu^{\mat A t_1} \vec b\,, \ \ t_1, t_2 \geq 0, t_3 > 0\,.
\end{aligned}
\end{equation}

\subsubsection*{$k$-th subsystem}
The symmetric output of the~$k$-th subsystem $y^{\mathrm{sym}}_k(t_1,\ldots,t_k)$ can be obtained by summing over all~$k!$ permutations~$\pi(\cdot)$ of the time variables~$(t_1,\ldots,t_k)$ of the triangular output~$y^{\triangle}_k(t_1,\ldots,t_k)$:
\begin{equation}
\begin{aligned}
y^{\mathrm{sym}}_k(t_1,\ldots,t_k) &:= \tfrac{1}{k!}\big(y^{\triangle}_k(t_1,\ldots,t_{k-1},t_k) + y^{\triangle}_k(t_1,\ldots,t_k,t_{k-1}) + \ldots
\big) \\
&:= \tfrac{1}{k!}\sum_{\pi(\cdot)} y^{\triangle}_k(t_{\pi(1)},\ldots,t_{\pi(k)})\, .
\end{aligned}
\end{equation}
The single-variable output~$y_k(t)$ of the~$k$-th subsystem can be calculated from~$y_k(t) = y^\triangle_k(t,\ldots,t) = y^{\mathrm{sym}}_k(t,\ldots,t)$.

The Laplace transform of the $k$-dimensional symmetric output~$Y^{\mathrm{sym}}_k(s_1,\ldots,s_k)$ becomes
\begin{equation}
\begin{aligned}
Y^{\mathrm{sym}}_k(s_1,\ldots,s_k) &= \mathcal{L}\{ y^{\mathrm{sym}}_k(t_1,\ldots,t_k)\}(s_1,\ldots,s_k) \\
&= \mathcal{L}\left\{ \tfrac{1}{k!}\big(y^\triangle_k(t_1,\ldots,t_{k-1},t_k) + y^\triangle_k(t_1,\ldots,t_k,t_{k-1}) + \ldots \big)\right\}(s_1,\ldots,s_k) \\
&= \tfrac{1}{k!}\big(Y^\triangle_k(s_1,\ldots,s_{k-1},s_k) + Y^\triangle_k(s_1,\ldots,s_k,s_{k-1}) + \ldots \big)\\
&= G^{\mathrm{sym}}_k(s_1,\ldots,s_k) \, U(s_1) \ldots U(s_k)\,,
\end{aligned}
\end{equation}
with the~$k$-dimensional symmetric transfer function
\begin{equation}
\begin{aligned}
G^{\mathrm{sym}}_k(s_1,\ldots,s_k) &= \tfrac{1}{k!} \big(G^\triangle_k(s_1,\ldots,s_{k-1},s_k) + G^\triangle_k(s_1,\ldots,s_{k},s_{k-1}) + \ldots \big) \\
&= \tfrac{1}{k!} \sum\limits_{\pi(\cdot)}G^\triangle_k(s_{\pi(1)},\ldots,s_{\pi(k)})\,. 
\end{aligned}
\end{equation}
The transfer function $G_k^{\mathrm{sym}}(s_1,\ldots,s_k) = Y_k^{\mathrm{sym}}(s_1,\ldots,s_k)$ s.t.~$U(s) = 1$ corresponds to the Laplace transform of the symmetric impulse response $g^{\mathrm{sym}}_k(t_1,\ldots,t_k) = y^{\mathrm{sym}}_k(t_1,\ldots,t_k)$ s.t.~$u(t)=\delta(t)$:
\begin{equation}
\begin{aligned}
G^{\mathrm{sym}}_k(s_1,\ldots,s_k) = \mathcal{L}\{ g^{\mathrm{sym}}_k(t_1,\ldots,t_k)\}(s_1,\ldots,s_k)\, .
\end{aligned}
\end{equation}
The symmetric impulse response of the~$k$-th subsystem~$g^{\mathrm{sym}}_k(t_1,\ldots,t_k)$ is the sum of $k!$ triangular impulse responses~$g^\triangle_k(t_1,\ldots,t_{k-1},t_k)$,~$g^\triangle_k(t_1,\ldots,t_k,t_{k-1})$, etc, which are non-zero on a corresponding $k$-simplex~$t_1 > \ldots > t_{k-1} > t_k > 0$,~$t_1 > \ldots > t_k > t_{k-1} > 0$, etc. The symmetric impulse response can be continuously extended at the one-dimensional intersection line, at the surface intersections, etc, since the limit values are equal. Therefore: 
\\[1ex]
\fbox{\parbox{14.7cm}{
\begin{equation}
\begin{aligned}
g^{\mathrm{sym}}_k(t_1,\ldots,t_k)
&= \tfrac{1}{k!}\vec c^{\mathsf T} \eu^{\mat A t_{i_k}}\mat N \eu^{\mat A (t_{i_{k-1}} - t_{i_k})} \mat N \cdots \mat N \eu^{\mat A (t_{i_2} - t_{i_3})} \mat N\eu^{\mat A (t_{i_1} - t_{i_2})}\vec b, \\
&\qquad t_{i_1} \geq \ldots \geq t_{i_k} > 0,\:\: (i_1,\ldots,i_k) \text{ permutations of } (1,\ldots,k)\,.
\end{aligned}
\end{equation}
}}\\[1ex]
This definition of the $k$-dimensional symmetric kernel is consistent with the impulse response of the~$k$-th subsystem $g_k(t) = g^{\mathrm{sym}}_k(t,\ldots,t) = 1/k! \, \vec c^{\mathsf T} \eu^{\mat A t}\mat N^k \vec b$ according to equation~\eqref{eq:imp-res-kth}.

All $k!$ $k$-simplexes border the intersection line $t_1 \!= \cdots =\! t_k > 0$ ($n \!=\! 1$-simplex), $(k-1)!$ $k$-simplexes border each cutting triangle ($n\!=\!2$-simplex), e.g. $t_1 = \ldots = t_{k-1} > t_k > 0$, etc, and $(k-n)!$ $k$-simplexes border each~$n$-simplex of the surface. Thus, the triangular impulse response of the $k$-th subsystem is
\\[1ex]
\fbox{\parbox{14.7cm}{
\begin{equation}
\begin{aligned} \label{eq:tri-imp-resp}
g^{\triangle}_{k}(t_1,\ldots,t_k) := \tfrac{1}{(k+1-n)!}\vec c^{\mathsf T} \eu^{\mat A t_k}\mat N\eu^{\mat A(t_{k-1} - t_k)}\mat N\!\cdots\!\mat N\eu^{\mat A(t_2 - t_3)}\mat N\eu^{\mat A(t_1 - t_2)}\vec b\,, \ \, t_1 \geq \ldots  \geq t_k > 0\,. 
\end{aligned}
\end{equation}
}}\\

Similarly, the $k$-dimensional regular impulse response is given by
\\[1ex]
\fbox{\parbox{14.7cm}{
\begin{equation}
\begin{aligned}	\label{eq:reg-imp-resp}
g^{\square}_{k}(t_1,\ldots,t_k) &= g^{\triangle}_{k}(t_1+t_2,t_2+t_3,\ldots,t_{k-1}+t_k,t_k) \\[0.2em]
&:= \tfrac{1}{(k+1-n)!}\vec c^{\mathsf T} \eu^{\mat A t_k} \mat N \eu^{\mat A t_{k-1}} \mat N \cdots \mat N \eu^{\mat A t_2} \mat N \eu^{\mat A t_1} \vec b\,, \ \ t_1, \ldots, t_{k-1} \geq 0, t_k > 0\,. 
\end{aligned}
\end{equation}
}} 

\subsection{Additional derivation} \label{subsec:add-der}
In the previous section, the symmetric kernels were introduced and continuously extended at the discontinuities. Afterwards, the definitions of the triangular and regular kernels were adjusted such that they correspond to the equal share of the symmetric kernels at the discontinuities, such as the time line~$t_1 = \ldots = t_k$ (for the triangular kernel domain) and~$t_1 = \ldots = t_{k-1}$ (for the regular kernel domain). This seems reasonable, but is not a proof. Hence, we will briefly show that:
\begin{equation}
\begin{aligned}
g^{\triangle}_2(t_1,t_2) = \lim\limits_{\varepsilon \rightarrow 0} y^{\triangle}_{2}(t_1,t_2)\ \ \text{subject to}\ \ u(t) = \delta_\varepsilon(t)\, , \\
g^{\square}_2(t_1,t_2) = \lim\limits_{\varepsilon \rightarrow 0} y^{\square}_{2}(t_1,t_2)\ \ \text{subject to}\ \ u(t) = \delta_\varepsilon(t)\,.
\end{aligned}
\end{equation}
Only the scalar case and the second subsystem will be considered, since the extension to higher subsystems is cumbersome.

\subsubsection*{Triangular impulse response}
\begin{figure}[htp]
	\begin{center}
		\renewcommand{\thesubfigure}{$\triangle$}
		\subfloat[]{\centering
			\begin{tikzpicture}
[xscale=1/(10 - (-3))*0.0352778*0.44*\textwidth,
yscale=1/(7.8 - (-0.5))*0.0352778*0.61803*0.44*\textwidth,
xtick/.style={below,fill=white,font=\tiny},
ytick/.style={left=0.2,fill=white,font=\tiny},
xlabel/.style={below right,font=\scriptsize},
ylabel/.style={left,font=\footnotesize},
ABnode/.style={circle, inner sep=0pt,minimum size=5pt,fill=black}
]
\def\tickin{0.1}
\def\tickout{0.1}
\def\eps{3}
\def\tzA{5} 
\def\tzB{2} 
\def\teAa{0.7} 
\def\teAb{3.4} 
\def\teAc{6.75}
\def\teAd{9.3} 
\def\teBa{1.3} 
\def\teBb{2.5} 
\def\teBc{3.5}
\def\teBd{6} 
\def\delta{0.5}
\def\xmax{10}
\def\xmin{0}
\def\ymax{10} 
\def\ymin{0}
\begin{scope}
\draw [->] (\xmin,0) -- (\xmax,0) node [xlabel] {$t_1$};
\draw [->] (0,\ymin) -- (0,\ymax) node [ylabel] {$t_2$};
\draw (\eps,0+\tickin) -- (\eps,0-\tickout) node[xtick, font=\footnotesize] {$\varepsilon$};
\draw (0-\tickout,\eps) -- (0+\tickin,\eps) node[ytick, font=\footnotesize] {$\varepsilon$};

\draw (0-\tickout,\tzA) -- (0+\tickin,\tzA) node[ytick] {$t_{2,\mathrm{A}}$};
\draw (0-\tickout,\tzB) -- (0+\tickin,\tzB) node[ytick] {$t_{2,\mathrm{B}}$};

\draw (\xmin,\eps) -- (\xmax,\eps);
\draw (\eps,\ymin) -- (\eps,\eps);
\draw (\xmin,\xmin+\eps) -- (\ymax-\eps,\ymax); 
\draw (\xmin,\xmin) -- (\ymax,\ymax); 
\draw (\ymin+\eps,\ymin) -- (\xmax,\xmax-\eps); 

\node[ABnode] (Aa) at (\teAa,\tzA) {};
\node[ABnode] (Ab) at (\teAb,\tzA) {};
\node[ABnode] (Ac) at (\teAc,\tzA) {};
\node[ABnode] (Ad) at (\teAd,\tzA) {};
\node[ABnode] (Ba) at (\teBa,\tzB) {};
\node[ABnode] (Bb) at (\teBb,\tzB) {};
\node[ABnode] (Bc) at (\teBc,\tzB) {};
\node[ABnode] (Bd) at (\teBd,\tzB) {};
\draw [-latex,dashed] (\xmin,\tzA) -- (\xmax,\tzA) node[xlabel] {}; 
\draw [-latex,dashed] (\xmin,\tzB) -- (\xmax,\tzB) node[xlabel] {}; 
\draw [ultra thin] (Aa) .. controls +(-0.2,0.4) and +(0.2,-0.4) .. ++(1,0.5) node[xshift=-6,yshift=8] {~(A1)};
\draw [ultra thin] (Ab) .. controls +(-0.2,0.4) and +(0.2,-0.4) .. ++(1.4,0.5) node[xshift=-6,yshift=8] {~(A2)};
\draw [ultra thin] (Ac) .. controls +(+0.2,-0.4) and +(-0.2,0.4) .. ++(-1.7,-0.5) node[xshift=6,yshift=-8] {~(A3)};
\draw [ultra thin] (Ad) .. controls +(+0.2,-0.4) and +(-0.2,0.4) .. ++(-0.8,-0.5) node[xshift=-6,yshift=-8] {~(A4)};
\draw [ultra thin] (Ba) .. controls +(-0.5,-0.8) and +(0.5,0.8) .. ++(0,-0.2-\tzB) node[xshift=-6,yshift=-8] {~(B1)};
\draw [ultra thin] (Bb) .. controls +(-0.5,-0.8) and +(0.5,0.8) .. ++(2,-0.2-\tzB) node[xshift=-6,yshift=-8] {~(B2)};
\draw [ultra thin] (Bc) .. controls +(-0.5,-0.8) and +(0.5,0.8) .. ++(3,-0.2-\tzB) node[xshift=-6,yshift=-8] {~(B3)};
\draw [ultra thin] (Bd) .. controls +(-0.5,-0.8) and +(0.5,0.8) .. ++(2.5,-0.2-\tzB) node[xshift=-6,yshift=-8] {~(B4)};
\node[rotate=45, above,font=\footnotesize] at (2*\eps,3*\eps) {$t_2 = t_1 + \varepsilon$};
\node[rotate=45, above,font=\footnotesize] at (2.5*\eps,2.5*\eps) {$t_2 = t_1 \phantom{+ \varepsilon}$};
\node[rotate=45, above,font=\footnotesize] at (3*\eps,2*\eps) {$t_2 = t_1 - \varepsilon$};

\end{scope}
\end{tikzpicture}%

			\label{fig:t1t2triangle}}%
		\renewcommand{\thesubfigure}{$\square$}	
		\subfloat[]{\centering
			\begin{tikzpicture}
[xscale=1/(10 - (-3))*0.0352778*0.44*\textwidth,
yscale=1/(7.8 - (-0.5))*0.0352778*0.61803*0.44*\textwidth,
xtick/.style={below,fill=white,font=\tiny},
ytick/.style={left=0.2,fill=white,font=\tiny},
xlabel/.style={below right,font=\scriptsize},
ylabel/.style={left,font=\footnotesize},
ABnode/.style={circle, inner sep=0pt,minimum size=5pt,fill=black}
]
\def\tickin{0.1}
\def\tickout{0.1}
\def\eps{3}
\def\tzA{5} 
\def\tzB{2} 
\def\teAa{0.7} 
\def\teAb{3.4} 
\def\teAc{6.75}
\def\teAd{9.3} 
\def\teBa{1.3} 
\def\teBb{2.5} 
\def\teBc{3.5}
\def\teBd{6} 
\def\delta{0.5}
\def\xmax{6}
\def\xmin{-\eps-\delta}
\def\ymax{10} 
\def\ymin{0}
\begin{scope}
\draw (-\eps,\eps) -- (\xmax,\eps);
\draw (\eps,0) -- (\eps,\ymax);
\draw (-\eps,\eps) -- (-\eps,\ymax);
\draw (-\xmax,\xmax) -- (\ymin,-\ymin) node[sloped, above,font=\tiny,fill=white,midway,xshift=0.1,yshift=0.1] {$t_2 = -t_1$};
\draw (0,\eps) -- (\eps,0); 

\node[ABnode] (Aa) at (\teAa-\tzA,\tzA) {};
\node[ABnode] (Ab) at (\teAb-\tzA,\tzA) {};
\node[ABnode] (Ac) at (\teAc-\tzA,\tzA) {};
\node[ABnode] (Ad) at (\teAd-\tzA,\tzA) {};
\node[ABnode] (Ba) at (\teBa-\tzB,\tzB) {};
\node[ABnode] (Bb) at (\teBb-\tzB,\tzB) {};
\node[ABnode] (Bc) at (\teBc-\tzB,\tzB) {};
\node[ABnode] (Bd) at (\teBd-\tzB,\tzB) {};
\draw [-latex,dashed] (-\tzA,\tzA) -- (\xmax,\tzA) node[xlabel] {}; 
\draw [-latex,dashed] (-\tzB,\tzB) -- (\xmax,\tzB) node[xlabel] {}; 
\draw [ultra thin] (Aa) .. controls +(-0.2,0.4) and +(0.2,-0.4) .. ++(0.6,0.5) node[xshift=-6,yshift=8] {~(A1)};
\draw [ultra thin] (Ab) .. controls +(-0.2,0.4) and +(0.2,-0.4) .. ++(0.6,0.5) node[xshift=-6,yshift=8] {~(A2)};
\draw [ultra thin] (Ac) .. controls +(-0.2,0.4) and +(0.2,-0.4) .. ++(0.6,0.5) node[xshift=-6,yshift=8] {~(A3)};
\draw [ultra thin] (Ad) .. controls +(-0.2,0.4) and +(0.2,-0.4) .. ++(0.6,0.5) node[xshift=-6,yshift=8] {~(A4)};

\draw [ultra thin] (Ba) .. controls +(0.5,-0.8) and +(-0.5,0.8) .. ++(-2,-0.2-\tzB) node[xshift=6,yshift=-8] {~(B1)};
\draw [ultra thin] (Bb) .. controls +(0.5,-0.8) and +(-0.5,0.8) .. ++(-1.3,-0.2-\tzB) node[xshift=6,yshift=-8] {~(B2)};
\draw [ultra thin] (Bc) .. controls +(0.5,-0.8) and +(-0.5,0.8) .. ++(-0.5,-0.2-\tzB) node[xshift=6,yshift=-8] {~(B3)};
\draw [ultra thin] (Bd) .. controls +(-0.5,-0.8) and +(0.5,0.8) .. ++(1,-0.2-\tzB) node[xshift=-6,yshift=-8] {~(B4)};
\draw [->] (\xmin,0) -- (\xmax,0) node [xlabel] {$t_1$};
\draw [->] (0,\ymin) -- (0,\ymax) node [ylabel] {$t_2$};
\draw (\eps,0+\tickin) -- (\eps,0-\tickout) node[xtick, font=\footnotesize] {$\varepsilon$};
\draw (0-\tickout,\eps) -- (0+\tickin,\eps) node[ytick, font=\footnotesize] {$\varepsilon$};

\draw (0-\tickout-\tzA,\tzA) -- (0+\tickin-\tzA,\tzA) node[left=0.2,font=\tiny] {$t_{2,\mathrm{A}}$};
\draw (0-\tickout-\tzB,\tzB) -- (0+\tickin-\tzB,\tzB) node[left=0.2,font=\tiny] {$t_{2,\mathrm{B}}$} ;
\end{scope}
\end{tikzpicture}%
			\label{fig:t1t2square}}\\[1ex]
	\end{center}
	\caption{Integration regions for~$t_1, t_2$ for the Laplace transform in the triangular~$(\triangle)$ and regular~$(\square)$ case.}
	\label{fig:t1t2}
	\vspace{-2em}
\end{figure}
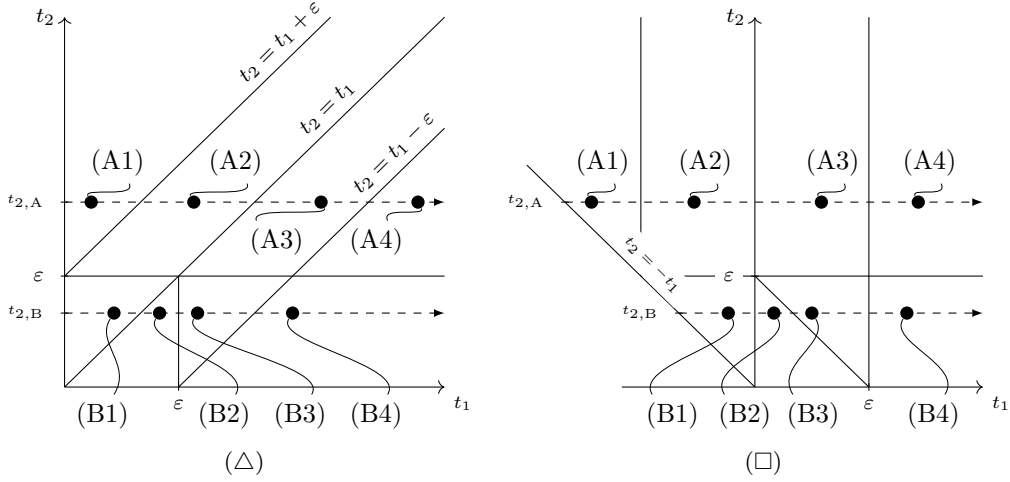

According to equation~\eqref{eq:sy_mu_mu2_1}, the triangular impulse response of the second subsystem $g^{\triangle}_{\varepsilon, 2}(t_1, t_2)$ resembles the triangular output~$y^\triangle_{2}(t_1,t_2)$ for an impulse input~$u(t) = \delta_\varepsilon(t)$ by a nascent delta function~$\delta_\varepsilon(t)$:
\begin{equation}
\begin{aligned}
\lefteqn{g^{\triangle}_{\varepsilon,2}(t_1,t_2) = \int_{\tau_1 = -\infty}^{\infty} \int_{\tau_2 = -\infty}^{\infty} c \, \eu^{a\tau_2}\sigma(\tau_2) \, n \,  \eu^{a\tau_1}\sigma(\tau_1) \, b} \\
&\quad\times \underbrace{\varepsilon^{-1}\sigma(t_2 - \tau_2)\sigma(\varepsilon - (t_2 - \tau_2))}_{\delta_\varepsilon(t_2 - \tau_2)}\underbrace{\varepsilon^{-1}\sigma(t_1 - \tau_2 - \tau_1)\sigma(\varepsilon - (t_1 - \tau_2 - \tau_1))}_{\delta_\varepsilon(t_1 - \tau_2 - \tau_1)}\,\dd \tau_2 \dd \tau_1\, .
\end{aligned}
\end{equation}
The integration limits are taken into account through six heaviside step functions. Consequently, there are altogether eight different cases for~$t_1,t_2,\varepsilon > 0$ that should be considered during the integration. Figure~\ref{fig:A+B} shows the eight integration regions for~$\tau_1, \tau_2$. \\
The cases for~$t_2 > \varepsilon$, i.e.~$t_2 - \varepsilon > 0$, are shown in the upper part of Figure~\ref{fig:A+B}, starting with~$t_1 < t_2 - \varepsilon$ in~\ref{fig:A1}. To obtain the subsequent integration regions, the dashed path~$t_{2\mathrm{A}}$ in Figure~\ref{fig:t1t2triangle} should be followed, increasing the variable~$t_1$ and letting~$t_2$ constant. Doing this first yields~$t_1 > t_2 - \varepsilon$ like in~\ref{fig:A2}, then~$t_1 > t_2$ like in~\ref{fig:A3} and finally~$t_1 - \varepsilon > t_2$ like in~\ref{fig:A4}. \\
The cases for~$t_2 < \varepsilon$, i.e.~$t_2 - \varepsilon < 0$, are shown in the lower part of Figure~\ref{fig:A+B}, starting with~$t_1 < t_2$ in~\ref{fig:B1}. To obtain the subsequent integration regions, the dashed path~$t_{2\mathrm{B}}$ in Figure~\ref{fig:t1t2triangle} should be followed. Doing this first yields~$t_1 > t_2$ like in~\ref{fig:B2}, then~$t_1 - \varepsilon > 0$ like in~\ref{fig:B3} and finally~$t_1 - \varepsilon > t_2$ like in~\ref{fig:B4}.

The calculation of the eight different cases and the computation of the limits~$\varepsilon \to 0$ were performed via the~\emph{Symbolic Math Toolbox} from \textsc{Matlab}. The resulting triangular impulse response 
\begin{equation}
\begin{aligned}
g^{\triangle}_{2}(t_1,t_2) &= \lim\limits_{\varepsilon \rightarrow 0} g^{\triangle}_{\varepsilon,2}(t_1,t_2) \\
&= \left\{\begin{array}{ll}
c\eu^{at_2}n\eu^{a(t_1 - t_2)}b, & t_1 > t_2 \\
\frac{1}{2}c\eu^{at}nb, & t_1 = t_2 = t \\
0, & \text{else} 
\end{array}  \right.
\end{aligned}
\end{equation}
confirms the definition from equation~\eqref{eq:g2-tri}.

\begin{figure}[hbp]
	\begin{center}
		\renewcommand{\thesubfigure}{A1}
		\subfloat[]{\centering
			\begin{tikzpicture}
[xscale=1/(10 - (-3))*0.0352778*0.4*\textwidth,
yscale=1/(7.8 - (-0.5))*0.0352778*0.61803*0.4*\textwidth,
xtick/.style={below,fill=white,font=\tiny},
ytick/.style={left=0.2,fill=white,font=\tiny},
xlabel/.style={below right,font=\scriptsize},
ylabel/.style={left,font=\footnotesize},
]
\def\tickin{0.1}
\def\tickout{0.1}
\def\eps{3}
\def\tz{5}
\def\delta{0.5}
\def\xmax{10}
\def\xmin{-3}
\def\ymax{7.8} 
\def\ymin{-\delta}

\begin{scope}
\def\te{0.7}
\draw [->] (0,\ymin) -- (0,\ymax) node [ylabel] {$\tau_1$};
\draw [->] (\xmin,0) -- (\xmax,0) node [xlabel] {$\tau_2$};
\draw (\te-\eps,0+\tickin) -- (\te-\eps,0-\tickout) node[xtick] {$t_1 - \varepsilon$};
\draw (\te,0+\tickin) -- (\te,0-\tickout) node[xtick] {$t_1$};
\draw (\tz-\eps,0+\tickin) -- (\tz-\eps,0-\tickout) node[xtick] {$t_2 - \varepsilon$};
\draw (\tz,0+\tickin) -- (\tz,0-\tickout) node[xtick] {$t_2$};
\draw [dashed] (\xmin,\te-\eps-\xmin) -- (\te-\eps,0);
\draw [dashed] (\xmin,\te-\xmin) -- (\te,0);
\draw [dashed] (\tz-\eps,0) -- (\tz-\eps,\ymax);
\draw [dashed] (\tz,0) -- (\tz,\ymax);
\draw[arrows={-latex}] (\tz + 1,0) --++ (0,1) node[xshift = 0.75cm, yshift = 0.3cm,fill=white, font=\tiny] {$\sigma(\tau_1)\, ,\: \tau_1 > 0$};
\draw[arrows={-latex}] (\tz,3.5) --++ (-1,0) node[xshift = 1.35cm,yshift = 0.3cm,fill=white,font=\tiny] {$\sigma(t_2 - \tau_2)\, ,\: \tau_2 < t_2$};
\draw[arrows={-latex}] (\tz-\eps,5) --++ (1,0) node[xshift = 1.35cm,yshift = 0.3cm,fill=white,font=\tiny] {$\sigma(\varepsilon - (t_2 - \tau_2)), \tau_2 > t_2 - \varepsilon$};
\draw[arrows={-latex}] (0,6.5) --++ (1,0) node[xshift = 0.5cm, yshift = 0.3cm,fill=white, font=\tiny] {$\sigma(\tau_2)\, ,\: \tau_2 > 0$};
\draw[arrows={-latex}] (\te-\eps,\eps) --++ (-\delta,-\delta) node[xshift = 1.07cm,yshift = 0.65cm,fill=white,font=\tiny,fill=white,align=center] {$\sigma(t_1 - \tau_1 - \tau_2)\,,$\\$\tau_1 < t_1 - \tau_2$};
\draw[arrows={-latex}] (\te-\eps-\delta,\delta) --++ (\delta,\delta) node[xshift = 1.5cm,yshift = 0.2cm,fill=white,font=\tiny,fill=white,align=center] {$\sigma(\varepsilon - (t_1 - \tau_2 - \tau_1))\,,$\\$\tau_1 > t_1 - \varepsilon - \tau_2$};
\end{scope}

\end{tikzpicture}%

			\label{fig:A1}}%
		\renewcommand{\thesubfigure}{A2}
		\subfloat[]{\centering
			\begin{tikzpicture}
[xscale=1/(10 - (-3))*0.0352778*0.4*\textwidth,
yscale=1/(7.8 - (-0.5))*0.0352778*0.61803*0.4*\textwidth,
xtick/.style={below,fill=white,font=\tiny},
ytick/.style={left=0.2,fill=white,font=\tiny},
xlabel/.style={below right,font=\scriptsize},
ylabel/.style={left,font=\footnotesize},
]
\def\tickin{0.1}
\def\tickout{0.1}
\def\eps{3}
\def\tz{5}
\def\delta{0.5}
\def\xmax{10}
\def\xmin{-3}
\def\ymax{7.8} 
\def\ymin{-\delta}

\begin{scope}
\draw [->] (0,\ymin) -- (0,\ymax) node [ylabel] {$\tau_1$};
\draw [->] (\xmin,0) -- (\xmax,0) node [xlabel] {$\tau_2$};
\def\te{3.4}
\draw (\te-\eps,0+\tickin) -- (\te-\eps,0-\tickout) node[xtick] {$t_1 - \varepsilon$};
\draw (\te,0+\tickin) -- (\te,0-\tickout) node[xtick] {$t_1$};
\draw (\tz-\eps,0+\tickin) -- (\tz-\eps,0-\tickout) node[xtick] {$t_2 - \varepsilon$};
\draw (\tz,0+\tickin) -- (\tz,0-\tickout) node[xtick] {$t_2$};
\draw [dashed] (\xmin,\te-\eps-\xmin) -- (\te-\eps,0);
\draw [dashed] (\xmin,\te-\xmin) -- (\te,0);
\draw [dashed] (\tz-\eps,0) -- (\tz-\eps,\ymax);
\draw [dashed] (\tz,0) -- (\tz,\ymax);
\draw[fill=black!5]  (\tz - \eps, 0) -- (\tz - \eps, \te - \tz + \eps) -- (\te, 0) -- cycle;
\end{scope}
\end{tikzpicture}%

			\label{fig:A2}}\\[1ex]	
		\renewcommand{\thesubfigure}{A3}
		\subfloat[]{\centering
			\begin{tikzpicture}
[xscale=1/(10 - (-3))*0.0352778*0.4*\textwidth,
yscale=1/(7.8 - (-0.5))*0.0352778*0.61803*0.4*\textwidth,
xtick/.style={below,fill=white,font=\tiny},
ytick/.style={left=0.2,fill=white,font=\tiny},
xlabel/.style={below right,font=\scriptsize},
ylabel/.style={left,font=\footnotesize},
]
\def\tickin{0.1}
\def\tickout{0.1}
\def\eps{3}
\def\tz{5}
\def\te{6.75}
\def\delta{0.5}
\def\xmax{10}
\def\xmin{-3}
\def\ymax{7.8} 
\def\ymin{-\delta}

\begin{scope}
\draw [->] (0,\ymin) -- (0,\ymax) node [ylabel] {$\tau_1$};
\draw [->] (\xmin,0) -- (\xmax,0) node [xlabel] {$\tau_2$};
\draw (\te-\eps,0+\tickin) -- (\te-\eps,0-\tickout) node[xtick] {$t_1 - \varepsilon$};
\draw (\te,0+\tickin) -- (\te,0-\tickout) node[xtick] {$t_1$};
\draw (\tz-\eps,0+\tickin) -- (\tz-\eps,0-\tickout) node[xtick] {$t_2 - \varepsilon$};
\draw (\tz,0+\tickin) -- (\tz,0-\tickout) node[xtick] {$t_2$};
\draw [dashed] (\xmin,\te-\eps-\xmin) -- (\te-\eps,0);
\draw [dashed] (\te-\ymax,\ymax) -- (\te,0);
\draw [dashed] (\tz-\eps,0) -- (\tz-\eps,\ymax);
\draw [dashed] (\tz,0) -- (\tz,\ymax);
\draw[fill=black!5]  (\te - \eps, 0) -- (\tz - \eps, \te - \tz) -- (\tz - \eps,\te - \tz + \eps) -- (\tz,\te - \tz) -- (\tz,0) -- cycle;
\end{scope}
\end{tikzpicture}%

			\label{fig:A3}}%
		\renewcommand{\thesubfigure}{A4}
		\subfloat[]{\centering
			\begin{tikzpicture}
[xscale=1/(10 - (-3))*0.0352778*0.4*\textwidth,
yscale=1/(7.8 - (-0.5))*0.0352778*0.61803*0.4*\textwidth,
xtick/.style={below,fill=white,font=\tiny},
ytick/.style={left=0.2,fill=white,font=\tiny},
xlabel/.style={below right,font=\scriptsize},
ylabel/.style={left,font=\footnotesize},
]
\def\tickin{0.1}
\def\tickout{0.1}
\def\eps{3}
\def\tz{5}
\def\te{9.3}
\def\delta{0.5}
\def\xmax{10}
\def\xmin{-3}
\def\ymax{7.8} 
\def\ymin{-\delta}

\begin{scope}
\draw [->] (0,\ymin) -- (0,\ymax) node [ylabel] {$\tau_1$};
\draw [->] (\xmin,0) -- (\xmax,0) node [xlabel] {$\tau_2$};
\draw (\te-\eps,0+\tickin) -- (\te-\eps,0-\tickout) node[xtick] {$t_1 - \varepsilon$};
\draw (\te,0+\tickin) -- (\te,0-\tickout) node[xtick] {$t_1$};
\draw (\tz-\eps,0+\tickin) -- (\tz-\eps,0-\tickout) node[xtick] {$t_2 - \varepsilon$};
\draw (\tz,0+\tickin) -- (\tz,0-\tickout) node[xtick] {$t_2$};
\draw [dashed] (\te-\eps-\ymax,\ymax) -- (\te-\eps,0);
\draw [dashed] (\te-\ymax,\ymax) -- (\te,0);
\draw [dashed] (\tz-\eps,0) -- (\tz-\eps,\te - \tz + \eps + \delta);
\draw [dashed] (\tz,0) -- (\tz,\ymax);
\draw[fill=black!5]  (\tz - \eps,\te - \eps - \tz + \eps) -- (\tz - \eps,\te - \tz + \eps) -- (\tz,\te - \tz) -- (\tz,\te - \eps - \tz) -- cycle;
\end{scope}
\end{tikzpicture}%

			\label{fig:A4}} \\[1.5ex]	
		\renewcommand{\thesubfigure}{B1}
		\subfloat[]{\centering
			\begin{tikzpicture}
[xscale=1/(7 - (-2.5))*0.0352778*0.4*\textwidth,
yscale=1/(7 - (-0.5))*0.0352778*0.61803*0.4*\textwidth,
xtick/.style={below,fill=white,font=\tiny},
ytick/.style={left=0.2,fill=white,font=\tiny},
xlabel/.style={below right,font=\scriptsize},
ylabel/.style={left,font=\footnotesize},
]
\def\tickin{0.1}
\def\tickout{0.1}
\def\eps{3}
\def\tz{2}
\def\te{1.3}
\def\delta{0.5}
\def\xmax{7}
\def\xmin{-2}
\def\ymax{7} 
\def\ymin{-\delta}

\begin{scope}
\draw [->] (0,\ymin) -- (0,\ymax) node [ylabel] {$\tau_1$};
\draw [->] (\xmin,0) -- (\xmax,0) node [xlabel] {$\tau_2$};
\draw (\te-\eps,0+\tickin) -- (\te-\eps,0-\tickout) node[xtick] {$t_1 - \varepsilon$};
\draw (\te,0+\tickin) -- (\te,0-\tickout) node[xtick] {$t_1$};
\draw (\tz,0+\tickin) -- (\tz,0-\tickout) node[xtick] {$t_2$};
\draw (\tz-\eps,0+\tickin) -- (\tz-\eps,0-\tickout);
\draw [dashed] (\xmin,\te-\eps-\xmin) -- (\te-\eps,0);
\draw [dashed] (\xmin,\te-\xmin) -- (\te,0);
\draw [dashed] (\tz,0) -- (\tz,\ymax);
\draw [dashed] (\tz-\eps,0) -- (\tz-\eps,\ymax);
\draw[fill=black!5]  (0,0) -- (0,\te) -- (\te,0) -- cycle;
\end{scope}
\end{tikzpicture}%

			\label{fig:B1}}%
		\renewcommand{\thesubfigure}{B2}	
		\subfloat[]{\centering
			\begin{tikzpicture}
[xscale=1/(7 - (-2.5))*0.0352778*0.4*\textwidth,
yscale=1/(7 - (-0.5))*0.0352778*0.61803*0.4*\textwidth,
xtick/.style={below,fill=white,font=\tiny},
ytick/.style={left=0.2,fill=white,font=\tiny},
xlabel/.style={below right,font=\scriptsize},
ylabel/.style={left,font=\footnotesize},
]
\def\tickin{0.1}
\def\tickout{0.1}
\def\eps{3}
\def\tz{2}
\def\te{2.5}
\def\delta{0.5}
\def\xmax{7}
\def\xmin{-2.5}
\def\ymax{7} 
\def\ymin{-\delta}

\begin{scope}
\draw [->] (0,\ymin) -- (0,\ymax) node [ylabel] {$\tau_1$};
\draw [->] (\xmin,0) -- (\xmax,0) node [xlabel] {$\tau_2$};
\draw (\te-\eps,0+\tickin) -- (\te-\eps,0-\tickout) node[xtick] {$t_1 - \varepsilon$};
\draw (\te,0+\tickin) -- (\te,0-\tickout) node[xtick] {$t_1$};
\draw (\tz,0+\tickin) -- (\tz,0-\tickout) node[xtick] {$t_2$};
\draw (\tz-\eps,0+\tickin) -- (\tz-\eps,0-\tickout);
\draw [dashed] (\xmin,\te-\eps-\xmin) -- (\te-\eps,0);
\draw [dashed] (\xmin,\te-\xmin) -- (\te,0);
\draw [dashed] (\tz,0) -- (\tz,\ymax);
\draw [dashed] (\tz-\eps,0) -- (\tz-\eps,\ymax);
\draw[fill=black!5]  (0,0) -- (0,\te) -- (\tz,\te - \tz) -- (\tz,0) -- cycle;
\end{scope}
\end{tikzpicture}%

			\label{fig:B2}}\\[1ex]	
		\renewcommand{\thesubfigure}{B3}
		\subfloat[]{\centering
			\begin{tikzpicture}
[xscale=1/(7 - (-2.5))*0.0352778*0.4*\textwidth,
yscale=1/(7 - (-0.5))*0.0352778*0.61803*0.4*\textwidth,
xtick/.style={below,fill=white,font=\tiny},
ytick/.style={left=0.2,fill=white,font=\tiny},
xlabel/.style={below right,font=\scriptsize},
ylabel/.style={left,font=\footnotesize},
]
\def\tickin{0.1}
\def\tickout{0.1}
\def\eps{3}
\def\tz{2}
\def\te{3.5}
\def\delta{0.5}
\def\xmax{7}
\def\xmin{-2.5}
\def\ymax{7} 
\def\ymin{-\delta}

\begin{scope}
\draw [->] (0,\ymin) -- (0,\ymax) node [ylabel] {$\tau_1$};
\draw [->] (\xmin,0) -- (\xmax,0) node [xlabel] {$\tau_2$};
\draw (\te-\eps,0+\tickin) -- (\te-\eps,0-\tickout) node[xtick] {$t_1 - \varepsilon$};
\draw (\te,0+\tickin) -- (\te,0-\tickout) node[xtick] {$t_1$};
\draw (\tz,0+\tickin) -- (\tz,0-\tickout) node[xtick] {$t_2$};
\draw (\tz-\eps,0+\tickin) -- (\tz-\eps,0-\tickout) node[xtick] {$t_2 - \varepsilon$};
\draw [dashed] (\xmin,\te-\eps-\xmin) -- (\te-\eps,0);
\draw [dashed] (\xmin,\te-\xmin) -- (\te,0);
\draw [dashed] (\tz,0) -- (\tz,\ymax);
\draw [dashed] (\tz-\eps,0) -- (\tz-\eps,\ymax);
\draw[fill=black!5]  (\te - \eps,0) -- (0,\te - \eps) -- (0,\te) -- (\tz,\te - \tz) -- (\tz, 0) -- cycle;
\end{scope}
\end{tikzpicture}%

			\label{fig:B3}}%
		\renewcommand{\thesubfigure}{B4}
		\subfloat[]{\centering
			\begin{tikzpicture}
[xscale=1/(7 - (-2.5))*0.0352778*0.4*\textwidth,
yscale=1/(7 - (-0.5))*0.0352778*0.61803*0.4*\textwidth,
xtick/.style={below,fill=white,font=\tiny},
ytick/.style={left=0.2,fill=white,font=\tiny},
xlabel/.style={below right,font=\scriptsize},
ylabel/.style={left,font=\footnotesize},
]
\def\tickin{0.1}
\def\tickout{0.1}
\def\eps{3}
\def\tz{2}
\def\te{6}
\def\delta{0.5}
\def\xmax{7}
\def\xmin{-2.5}
\def\ymax{7} 
\def\ymin{-\delta}

\begin{scope}
\draw [->] (0,\ymin) -- (0,\ymax) node [ylabel] {$\tau_1$};
\draw [->] (\xmin,0) -- (\xmax,0) node [xlabel] {$\tau_2$};
\draw (\te-\eps,0+\tickin) -- (\te-\eps,0-\tickout) node[xtick] {$t_1 - \varepsilon$};
\draw (\te,0+\tickin) -- (\te,0-\tickout) node[xtick] {$t_1$};
\draw (\tz,0+\tickin) -- (\tz,0-\tickout) node[xtick] {$t_2$};
\draw (\tz-\eps,0+\tickin) -- (\tz-\eps,0-\tickout) node[xtick] {$t_2 - \varepsilon$};
\draw [dashed] (\xmin,\te-\eps-\xmin) -- (\te-\eps,0);
\draw [dashed] (\te-\ymax,\ymax) -- (\te,0);
\draw [dashed] (\tz,0) -- (\tz,\ymax);
\draw [dashed] (\tz-\eps,0) -- (\tz-\eps,\ymax);
\draw[fill=black!5]  (0,\te - \eps) -- (0,\te) -- (\tz,\te - \tz) -- (\tz,\te - \eps - \tz) -- cycle;
\end{scope}
\end{tikzpicture}%

			\label{fig:B4}} 
	\end{center}
	\caption{Integration regions for~$\tau_1, \tau_2$ in the triangular domain~$(\triangle)$. -- (A) for $t_2 > \varepsilon$, i.e.~$t_2 - \varepsilon > 0$ and (B) for $t_2 < \varepsilon$, i.e.~$t_2 - \varepsilon < 0$.}
	\label{fig:A+B}
\end{figure}
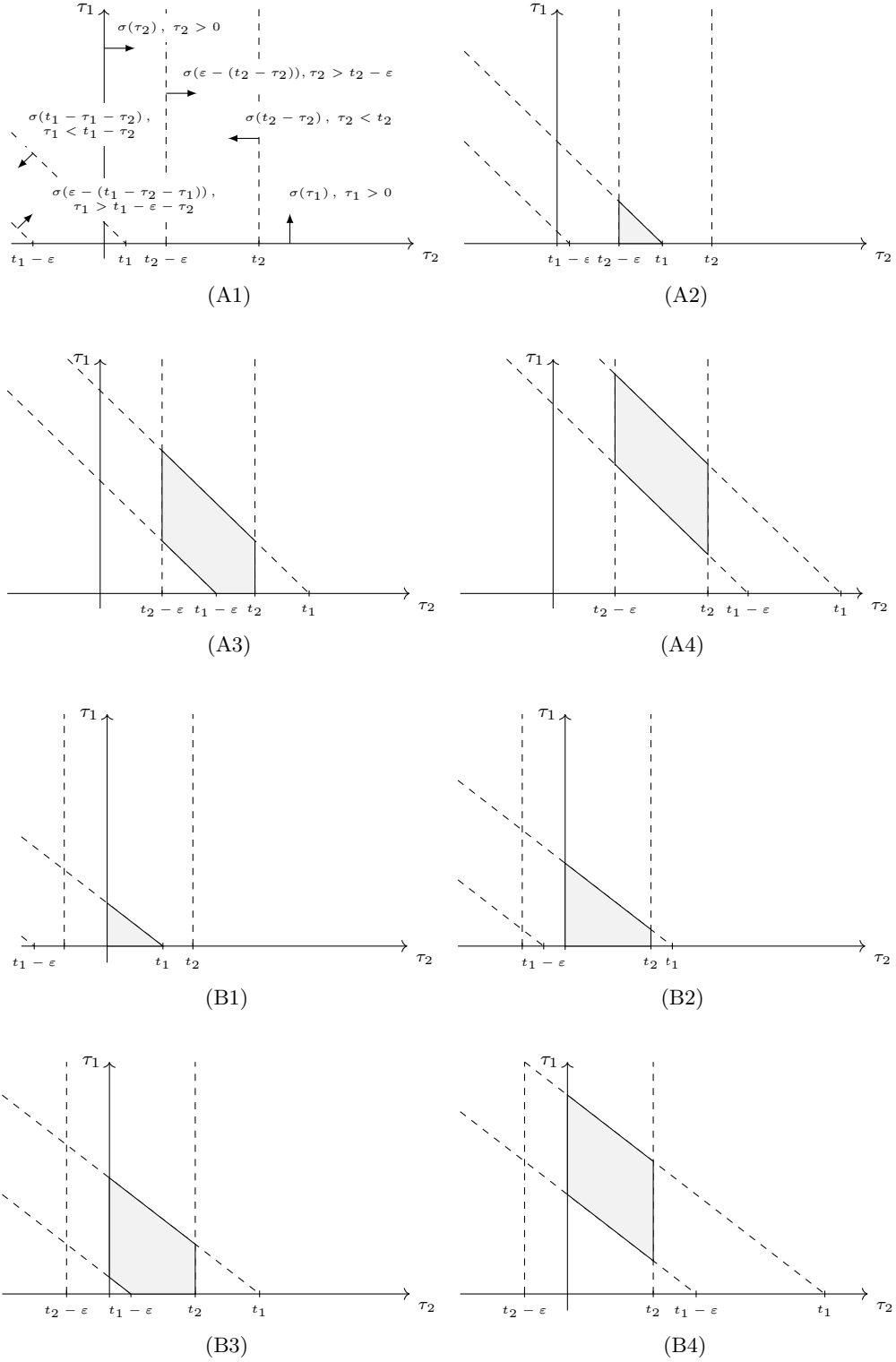	

Turning back to the frequency domain, the Laplace transform of the triangular impulse response~$g^{\triangle}_{\varepsilon,2}(t_1, t_2)$ is:
\begin{equation}
\begin{aligned}
G^{\triangle}_{\varepsilon,2}(s_1,s_2) &= \int_{t_1 = 0}^{\infty}\int_{t_2 = 0}^{\infty}g^{\triangle}_{\varepsilon,2}(t_1,t_2)\eu^{-s_1t_1}\eu^{-s_2t_2}\,\dd t_2 \dd t_1\, .
\end{aligned}
\end{equation}
Eight different integration regions for~$t_1, t_2$ should be considered here again. Figure~\ref{fig:t1t2triangle} shows the eight regions which correspond to the eight cases from Figure~\ref{fig:A+B}. The computation with the \emph{Symbolic Math Toolbox} leads to the Laplace transform
\begin{equation}
\begin{aligned}
G^{\triangle}_{\varepsilon,2}(s_1,s_2) &= c(s_1 + s_2 - a)^{-1}n(s_1 - a)^{-1}b \frac{1 - \eu^{-s_1 \varepsilon}}{s_1 \varepsilon}\frac{1 - \eu^{-s_2 \varepsilon}}{s_2 \varepsilon} \\
&= G^{\triangle}_{2}(s_1,s_2)\mathcal{L}\{\delta_\varepsilon(t)\}(s_1)\mathcal{L}\{\delta_\varepsilon(t)\}(s_2)\, ,
\end{aligned}
\end{equation}
which confirms the triangular output equation~\eqref{eq:Y-tri-SISO} derived from the convolution property. Building the limit~$\varepsilon \to 0$ in the frequency domain finally yields the triangular transfer function of the second subsystem (L'Hospital's rule):
\begin{equation}
\begin{aligned}
G^{\triangle}_{2}(s_1,s_2) &= \lim\limits_{\varepsilon \rightarrow 0}c(s_1 + s_2 - a)^{-1}n(s_1 - a)^{-1}b \frac{1 - \eu^{-s_1 \varepsilon}}{s_1 \varepsilon}\frac{1 - \eu^{-s_2 \varepsilon}}{s_2 \varepsilon} \\
&\!\overset{\mathrm{L.H.}}{=} c(s_1 + s_2 - a)^{-1}n(s_1 - a)^{-1}b\, .
\end{aligned}
\end{equation}
Thus, it can be concluded that building first the limit in the time domain~$g^\triangle_{\varepsilon,2}(t_1,t_2) \overset{\varepsilon \rightarrow 0}{\rightarrow} g^\triangle_2(t_1,t_2)$ and performing the Laplace transform~$G^\triangle_2(s_1,s_2) = \mathcal{L}\{g^\triangle_2(t_1,t_2)\}(s_1,s_2)$ afterwards is equivalent to first computing the Laplace transform~$G^\triangle_{\varepsilon,2}(s_1,s_2) = \mathcal{L}\{g^\triangle_{\varepsilon,2}(t_1,t_2)\}(s_1,s_2)$ and then building the limit in the frequency domain~$G^\triangle_{\varepsilon,2}(s_1,s_2) \overset{\varepsilon \rightarrow 0}{\rightarrow}G^\triangle_{2}(s_1,s_2)$.

\subsubsection*{Regular impulse response}
For the regular form of the impulse response we proceed in a similar way. According to~\eqref{eq:Y2-reg}, the regular impulse response of the second subsystem $g^{\square}_{\varepsilon, 2}(t_1, t_2)$ resembles the regular output~$y^\square_{2}(t_1,t_2)$ for an impulse input~$u(t) = \delta_\varepsilon(t)$ by a nascent delta function~$\delta_\varepsilon(t)$. The resulting eight different integration regions for~$\tau_1, \tau_2$ are obtained by replacing~$t_1 \rightarrow t_1 + t_2$, $t_2 \rightarrow t_2$ in Figure~\ref{fig:A+B}. Using the \emph{Symbolic Math Toolbox} from \textsc{Matlab} yields the regular impulse response
\begin{equation}
\begin{aligned}
g^{\square}_{2}(t_1,t_2) &= \lim\limits_{\varepsilon \rightarrow 0} g^{\square}_{\varepsilon,2}(t_1,t_2) \\
&= \left\{\begin{array}{ll}
c\eu^{at_2}n\eu^{at_1}b, & t_1 > 0, t_2 > 0 \\
\frac{1}{2}c\eu^{at_2}nb, & t_1 = 0, t_2 > 0 \\
0, & \text{else}\, .
\end{array}  \right.
\end{aligned}
\end{equation}
For the Laplace transform of~$g^{\square}_{\varepsilon,2}(t_1, t_2)$ the eight integration regions for~$t_1, t_2$ from Figure~\ref{fig:t1t2square} are used. The computation via the \emph{Symbolic Math Toolbox} yields the Laplace transform
\begin{equation}
\begin{aligned}
G^{\square}_{\varepsilon,2}(s_1,s_2) &= c(s_2 - a)^{-1}n(s_1 - a)^{-1}b \frac{1 - \eu^{-s_1 \varepsilon}}{s_1 \varepsilon}\frac{1 - \eu^{-(s_2 - s_1)\varepsilon}}{(s_2 - s_1) \varepsilon} \\
&= G^{\square}_{2}(s_1,s_2)\mathcal{L}\{\delta_\varepsilon(t)\}(s_1)\mathcal{L}\{\delta_\varepsilon(t)\}(s_2 - s_1)\, ,
\end{aligned}
\end{equation}
which corresponds to the regular output equation~\eqref{eq:Y-reg-SISO} derived from the convolution property. Building the limit~$\varepsilon \to 0$ in the frequency domain finally yields the regular transfer function of the second subsystem (L'Hospital's rule):
\begin{equation}
\begin{aligned}
G^{\square}_{2}(s_1,s_2) &= \lim\limits_{\varepsilon \rightarrow 0}c(s_2 - a)^{-1}n(s_1 - a)^{-1}b \frac{1 - \eu^{-s_1 \varepsilon}}{s_1 \varepsilon}\frac{1 - \eu^{-(s_2 - s_1) \varepsilon}}{(s_2 - s_1) \varepsilon} \\
&\!\overset{\mathrm{L.H.}}{=} c(s_2 - a)^{-1}n(s_1 - a)^{-1}b\, .
\end{aligned}
\end{equation}
To sum up: with this alternative derivation, the definitions for the kernels from the previous section are confirmed again. 

\subsection{Connections and final remarks}
\begin{figure}[h!]
	\begin{center}
		\input{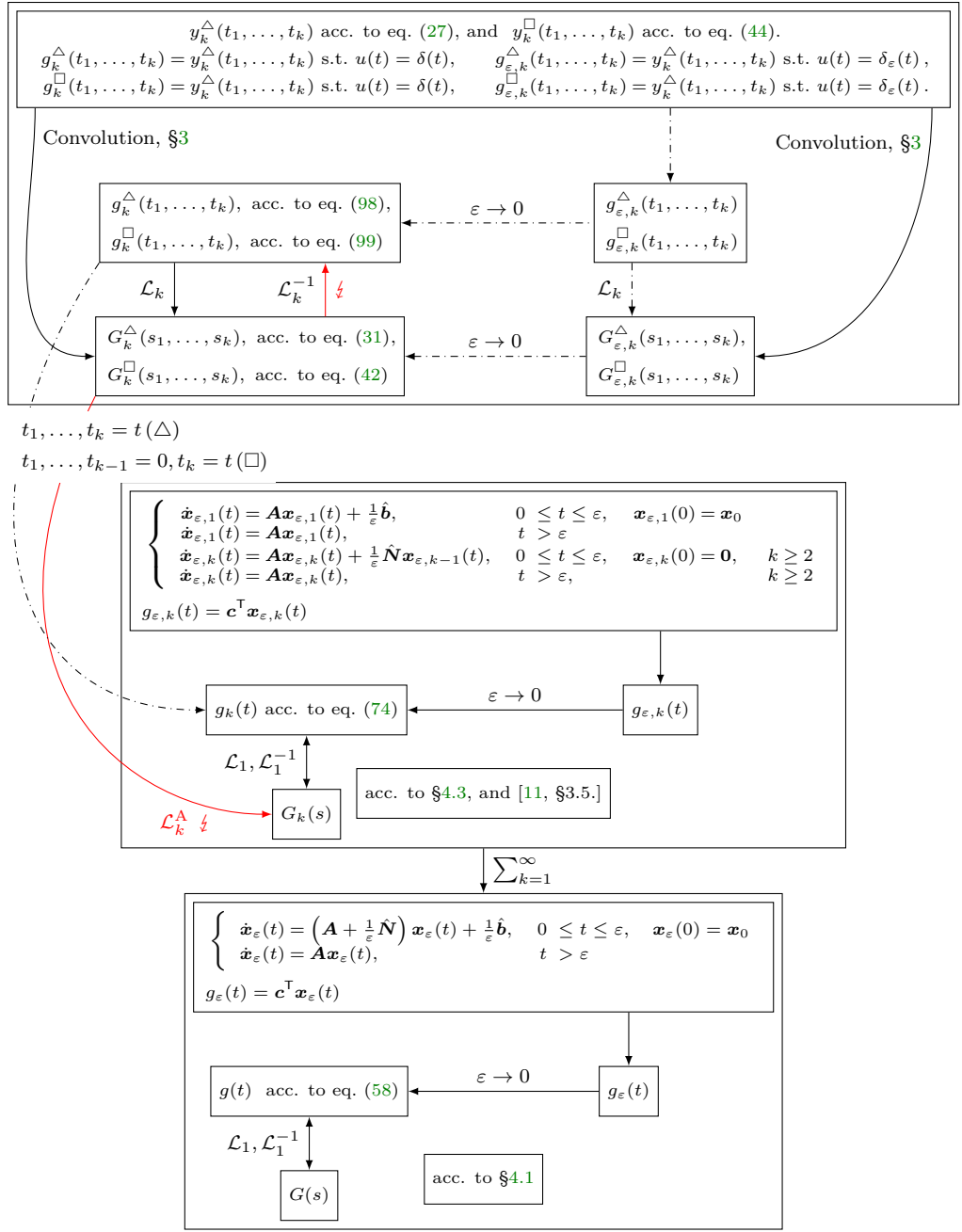}			
	\end{center}
	\caption{Overview and connections.-- The expressions derived in sections \ref{subsec:derivation-imp-direct} (lower box) and \ref{subsec:derivation-imp-k-subs} (middle box) are here considered for the SISO case ($m \!=\! 1$, $p\!=\!1$). Hence, it follows $u(t) \!=\! \mu\delta_\varepsilon(t) \!=\! \delta_\varepsilon(t)$, with the scaling of input~$\mu \!=\! 1$, $\hat{\mat N} \!=\! \sum_{j = 1}^{m}\mat N_j \mu_j = \mat N$ and $\hat{\vec b} = \sum_{j = 1}^{m}\vec b_j \mu_j = \vec b$. The derivations corresponding to dashed arrows are only shown examplarily for the second subsystem in section~\ref{sec:5}. Red arrows corresponds to cases, where the discontinuity is not taken into account.}
	\label{fig:123}
\end{figure}

Figure~\ref{fig:123} gives an overview on the main contents of this paper and the links between the derived statements.
We started with the triangular and regular auxiliary outputs~$y_k^{\triangle}(t_1,\ldots,t_k)$ and~$y_k^{\square}(t_1,\ldots,t_k)$. Using the convolution property of the multidimensional Laplace transform, in section~\ref{sec:inp-out-freq-dom} we derived input-output representations in the frequency domain in terms of the transfer functions~$G^\triangle_k(s_1,\ldots,s_k)$ and~$G^\square_k(s_1,\ldots,s_k)$, which are the multidimensional Laplace transform of the kernels~$g_k^{\triangle}(t_1,\ldots,t_k)$ and~$g_k^{\square}(t_1,\ldots,t_k)$, respectively.

Since the kernels are discontinuous at the time line $t_1,\ldots,t_k\!=\!t$ (for the triangular case) and $t_1,\ldots,t_{k-1} \!=\!0, t_k \!=\! t$ (for the regular case), the impulse response can neither be calculated by the $k$-dimensional inverse Laplace transform $\mathcal{L}^{-1}_k$ nor by the associated Laplace transform $\mathcal{L}^{\mathrm{A}}_k$, as the value at the discontinuity gets lost during the integration. Therefore, in section~\ref{subsec:derivation-imp-direct} we derived the impulse response~$g(t)$ of bilinear systems in the time domain by computing the output~$y(t)$ for a nascent delta function $u(t) \!=\! \delta_{\varepsilon}(t)$ in the limit $\varepsilon \to 0$. Similarly, in section~\ref{subsec:derivation-imp-k-subs} the impulse response~$g_{k}(t)$ of the~$k$-th subsystem was derived in the time domain as the output~$y_{k}(t)$ of the~$k$-th subsystem subject to a nascent delta function $\delta_{\varepsilon}(t)$ in the limit $\varepsilon \to 0$.  Note that $g(t)$ can be also obtained from $g_{k}(t)$ by summing over all subsystems. Particularly remarkable for the impulse response is the fact that a factor $1/k!$ arises in the equations~\eqref{eq:impulse-response} and~\eqref{eq:imp-res-kth}, which is normally not included in the triangular and regular kernels~$g^{\triangle}_{k}(t_1,\ldots,t_k)$ and~$g^{\square}_{k}(t_1,\ldots,t_k)$, since they are usually not defined along lines of equal time arguments. However, when impulse inputs are considered, then the value of a kernel at these precise locations becomes important for the input-output behavior. 

Therefore, we introduced the symmetric kernels in section~\ref{subsec:sym-kernels} and adjusted the definitions of the triangular and regular kernels, especially for the time line $t_1,\ldots,t_k=t$ in the triangular case (cf. equation~\eqref{eq:tri-imp-resp}) and for the time line $t_1,\ldots,t_{k-1}=0,t_k = t$ in the regular case (cf. equation~\eqref{eq:reg-imp-resp}), such that the kernels are consistent with the impulse response. Consequently, the impulse response~$g_k(t)$ according to \eqref{eq:imp-res-kth} can be then obtained from the adjusted triangular and regular kernels directly in the time domain: $g_k(t) \!=\! g^\triangle_k(t,\ldots,t) =  g^\square_k(0,\ldots,0,t)$.

\section{Conclusions} \label{sec:conclusions}
In this paper, we have extensively studied the bilinear systems theory using the Volterra series representation and derived the impulse response of bilinear systems. Afterwards, the triangular and regular Volterra kernels have been adjusted especially along lines of equal time arguments. After the adjustment, the kernels are consistent with the impulse response and characterize the bilinear system even for impulse inputs. Nevertheless, the adjustment does not influence the transfer functions, so that the computation of the output response of a bilinear system over the frequency domain using the auxiliary outputs and the $k$-dimensional inverse Laplace transform or associated Laplace transform can still not be performed (at least for impulse inputs). We therefore advocate not to speak about multidimensional triangular/regular transfer functions, as they should be valid for all Laplace transformable input signals, but they are at least not valid for the Laplace transformable delta function. In fact, similar to the linear setting, they rather should be seen as the scaling of a sum of exponential inputs, as it is regarded by the Growing Exponential Approach. We therefore recommend to use this approach, since derivations become shorter and easier, and the considerations are conducted directly in time domain.

\section*{Acknowledgments} 
We thank Dr. Klaus-Dieter Reinsch from the Chair of Numerical Mathematics at TUM for the valuable remarks, suggestions and discussions concerning some of the derivations. Furthermore, we want to thank Dr. Tobias Breiten from KFU Graz for his interest in the topic and his useful hints.






\bibliographystyle{siamplain}
\bibliography{NLMOR_Bilinear}
\end{document}